\documentclass[12pt]{article}
\usepackage{amssymb}
\usepackage{amsmath,amsthm,bbm}
\usepackage{version}
\usepackage{color}
\usepackage{datetime}
\usepackage[normalem]{ulem}
\usepackage{graphicx}
\usepackage{tikz}
\usetikzlibrary{calc,intersections}

\usepackage{float}

\definecolor{dark_purple}{rgb}{0.4, 0.0, 0.4}
\definecolor{dark_green}{rgb}{0.0, 0.7, 0.0}

\usepackage{latexsym}
\usepackage{amssymb}
\usepackage{euscript}

\def\mcc{M\raise.5ex\hbox{c}C}
\def\mccarthy{M\raise.5ex\hbox{c}Carthy}

\def\K{{\cal K}}
\def\M{{\cal M}}

\def\LL{{\cal L}}

\def\De{\Delta}

\def\a{\alpha}
\def\l{\lambda}
\def\z{\zeta}

\def\th{\theta}

\def\vare{\varepsilon}

\let\i=\infty

\def\={\ = \ }

\def\BB{{\cal B}}
\def\CC{{\cal C}}

\def\I{{\cal I}}

\def\C{\mathbb C}
\def\R{\mathbb R}
\def\T{\mathbb T}
\def\D{\mathbb D}
\def\B{\mathbb B}
\def\Z{\mathbb Z}

\def\be{\setcounter{equation}{\value{theorem}} \begin{equation}}
\def\ee{\end{equation} \addtocounter{theorem}{1}}
\def\beq{\begin{eqnarray*}}
\def\eeq{\end{eqnarray*}}

\def\bp{{\sc Proof: }}
\def\ep{{}{\hfill $\Box$} \vskip 5pt \par}

\def\oec{{}{\hspace*{\fill} $\lhd$} \vskip 5pt \par}
\def\bl{\begin{lemma}}
\def\el{\end{lemma}}
\def\bt{\begin{theorem}}
\def\et{\end{theorem}}
\def\bprop{\begin{prop}}
\def\eprop{\end{prop}}
\def\bd{\begin{definition}}
\def\ed{\end{definition}}
\def\br{\begin{remark}}
\def\er{\end{remark}}
\def\bexer{\begin{exercise}}
\def\eexer{\end{exercise}}
\def\bfig{\begin{figure}}
\def\efig{\end{figure}}

\newtheorem{theorem}{Theorem}[section]
\newtheorem{prop}[theorem]{Proposition}
\newtheorem{lemma}[theorem]{Lemma}
\newtheorem{cor}[theorem]{Corollary}

\newtheorem{question}[theorem]{Question}
\newtheorem{remark}[theorem]{Remark}
\newtheorem{defin}[theorem]{Definition}
\newtheorem{example}[theorem]{Example}

\def\be{\begin{equation}}
\def\ee{\end{equation}}

\def\conv{\operatorname{conv}}

\renewcommand\M{{\mathcal M}}

\renewcommand\O{\Omega}

\newcommand\U{{\mathcal U}}

\renewcommand\o{{\omega}}

\newcommand\car{Carath\'eodory\ }
\newcommand\hol{{\mathcal O}}
\newcommand\om{\omega}
\newcommand\al{\alpha}
\renewcommand\b{\beta}
\newcommand\g{\gamma}
\newcommand\gc{{\gamma_c}}
\renewcommand\CC{{\C}}
\newcommand\omegadotlambda{\omega \cdot \lambda }
\newcommand\taudotlambda{\tau \cdot \lambda }
\newcommand\hatt{\hat T}
\newcommand\hatd{\hat D}
\newcommand\EE{{\mathcal E}}
\newcommand\revp{\overset{\leftarrow}{p}}
\newcommand\revq{\overset{\leftarrow}{q}}
\newcommand\revP{\overset{\leftarrow}{P}}
\newcommand\X{{\mathcal X}}
\newcommand\eith{e^{i \theta}}

\newcommand\dimv{d}
\newcommand\dimpol{N}
\newcommand\dimi{n}

\newtheorem{assumption}{Assumption}

\renewcommand\be{\begin{equation}}
\renewcommand\ee{\end{equation}}
    
\newcommand\wrs{weak Roman surface}
\newcommand\roms{Roman surface}
\newcommand\ovm{\overline{\M}}
\newcommand\pp{\partial_+ \hatt}
\newcommand\ppm{\partial_- \hatt}
\newcommand\vp{\varphi}

\title{Totally geodesic sets and Carath\'eodory geometry in the polydisk}
\author{\L{}ukasz Kosi\'nski
 \thanks{Partially supported by Sheng grant no. 2023/48/Q/ST1/00048 of the
National Science Center. ORCID 0000-0002-9115-8545.}\\
Jagiellonian University\\
 Krakow\\
  lukasz.kosinski@uj.edu.pl
\and
John E. M\raise.5ex\hbox{c}Carthy
\thanks{Partially supported by National Science Foundation Grant  
DMS 2554487. ORCID 0000-0003-0036-7606.
}\\
Washington University in St. Louis\\
mccarthy@wustl.edu
}

\date{\today}
\begin{document}

\bibliographystyle{plain}
\maketitle
\begin{abstract}

A subset $V$ of the polydisk $\D^\dimpol$ is called totally geodesic if every pair of points in $V$ lie on a 
 Kobayashi geodesic of $\D^\dimpol$ that   lies entirely in $V$. It is called a \car set if its intrinsic \car metric agrees with the \car metric inherited from $\D^\dimpol$. 

We prove that for algebraic subvarieties of the polydisk these two notions coincide, and we obtain a complete geometric classification of such sets. The invariant metric geometry of the polydisk is unexpectedly rigid: every such set is, up to automorphism, built from only two basic geometric objects: disks and a distinguished family of two-dimensional algebraic surfaces that we call {\em {\roms}s}.

{\roms}s arise independently in the study of 
the three-point interpolation problem in the polydisk.

  \end{abstract}

\section{Overview}
\label{seca}

Let $\D$ denote the unit disk in the complex plane $\C$, and let $\rho$ be the pseudo-hyperbolic distance on $\D$ defined by
\[
\rho(\l,\mu) = \left| \frac{\l - \mu}{1 - \bar \l \mu} \right| .
\]
The Schwarz-Pick lemma says that every holomorphic function from $\D$ to $\D$ is contractive with respect to $\rho$.
C. \car \cite{cara27} had the idea to port this distance to any set in $\C^N$. Let $\hol(U,G)$ denote the holomorphic maps from $U$ to $G$.
\begin{defin}
\label{defa1}
Let $V \subseteq \C^N$. The \car pseudodistance on $V$ is given by
\be
\label{Eq1}
d_V(\l, \mu) = \sup \{ \rho(f(\l), f(\mu)) : f \in \hol(V, \D) \} .
\ee
\end{defin}

Let $\O$ be a pseudoconvex  domain in $\C^N$, and let $V \subseteq \O$ be an analytic subvariety.
If $\l,\mu$ are points in $V$, then in general $d_V (\l, \mu) \ge d_\Omega (\l,\mu)$. If these distances are always equal,
we shall say that $V$ is a \car set in $\O$. It suggests that there is some special connection between $V$ and $\O$.
\begin{question}
\label{qa1}
For a given domain $\O$, what subvarieties $V$ are  \car sets?
\end{question}

S. Kobayashi \cite{ko67} dualized Carath\'eodory's construction, and used functions in $\hol(\D,\O)$ to define a   Finsler metric on $\O$ by pushing forward the hyperbolic metric on $\D$. A geodesic is a map $\phi \in \hol(\D,\O)$ that is an isometry with respect to these metrics. A set $V \subseteq \O$ is called totally geodesic if any two points of $V$ are contained in the range of a geodesic that lies in $V$. The Kobayashi version of Question \ref{qa1} is:
\begin{question}
\label{qa2}
For a given domain $\O$, what subvarieties $V$ are  totally geodesic?
\end{question}
Question \ref{qa2} can be thought of as a complex analogue of describing the totally geodesic submanifolds of a Riemannian manifold.

When $\O$ is the ball, one can show that the answer to both Questions is that $V$ must be a retract \cite{kmc19}, a very strong restriction.
The object of this paper is to study the much more subtle case when $\Omega$ is the polydisk, which is not strongly convex. We restrict ourselves to considering only subvarieties $V$ that are algebraic, i.e. the intersection of an algebraic set with the polydisk. Then we show that the answer to both questions \ref{qa1} and \ref{qa2} is the same, and that $V$ is built out of 
one-dimensional disks and special two-dimensional complete intersections of quadrics which we call {\roms}s. 

Let us say that a triple $(A,B,C) \in \C^3$ is {\em triangular} if $|A|, |B|, |C|$ form the sides of a triangle (i.e.  $2  \times \max ( |A|,|B|, |C|) <  |A| + |B| + |C|$). A {\em \roms\ of height $1$} is a set of the form
\be
\label{Eq2}
\K_{(A,B,C)} = \{ (x,y,z) \in \D^3 : A x + B y + C z = \bar A yz + \bar B xz + \bar C xy \} 
\ee
where $(A,B,C)$ is triangular. 
If $(A, B, C) \in \C^3 \setminus \{0\}$ then the set \eqref{Eq2} is a retract of $\D^3$ 
if and only if $(A,B,C)$ is not triangular.
Retracts are always totally geodesic and  \car sets. The interesting case is when the triple is triangular;
in this case it was shown in \cite{kz21} that $\K_{(A,B,C)} $ is still totally geodesic and a \car set in $\D^3$.

A {\em \roms\ of height $m$} is given by $m$ triangular triples $(A_j, B_j, C_j)$ that satisfy a compatibility condition:
there exists $\tau = (\tau_1, \tau_2, \tau_3)$ in $\T^3$, a triple of unimodular complex numbers, so that
\be
\label{Eq3}
\notag
A_j \tau_1 + B_j \tau_2 + C_j \tau_3 = 0 \quad \forall\ 1 \le j \le m .
\ee
Given this data, the associated \roms\ is
\[
 \{ (x,y, z_1, \dots, z_m) \in \D^{m+2} : A_j x + B_j y + C_j z_j = \bar A_j y z_j + \bar B_j x z_j + \bar C_j xy , \ 1 \le j \le m \} .
\]

If $U \subset \D^m$ and $N \ge m$ we shall say that $V \subseteq \D^N$ is a graph over $U$ if there exists $F \in \hol (U, \D^{N-m})$ so that $V = \{ (x, F(x)) : x \in U \}$. We allow $N = m$, so $U$ is a graph over itself.
Here is our main theorem.

\bt
\label{thma1}
Let $V$ be an algebraic subvariety of $\D^N$, and assume $V$ is not a single point. Then the following are equivalent.

(i) $V$ is a \car set in $\D^N$.

(ii) $V$ is totally geodesic in $\D^N$.

(iii) After an automorphism of $\D^N$, $V$ is a graph over a product of copies of $\D$ and {\roms}s.
\et

This can be thought of as a rigidity theorem. Condition (i) is a metric condition, condition (ii) is a geometric condition, and condition (iii) says the set has a very specific form of algebraic structure.
The  assumption that $V$ is an algebraic subvariety is used in the proof: it ensures that the local parametrizing functions satisfy polynomial equations, and this gives control over their branching and analytic continuation.
 A natural question is whether  the theorem is true for arbitrary  analytic subvarieties. Nevertheless, the algebraic case already contains the main geometric difficulties, and the theorem gives a complete classification in this case.

When $N=2$, Theorem \ref{thma1} is the main result of \cite{agmcvn}, though the situation is much simpler since Roman surfaces occur for $N\geq 3$. 
 In the bidisk, $(iii)$ simply says that $V$ is a retract of  $\D^2$. 

The Cartesian product of  sets $V_1$ in $\D^{N_1}$ and $V_2$ in $\D^{N_2}$ is a \car set in $\D^{N_1 + N_2}$ if and only if
both $V_1$ and $V_2$ are \car sets
(see Theorem \ref{thmc1}), and a graph over a  set is a \car set if and only if the base set is a \car set. So if we define a {\em basic \car set} in the polydisk as one that cannot be written either as a Cartesian product, or as a graph over a \car set in a lower dimensional polydisk,
 then it follows that
every basic \car set in the polydisk is either a disk or an automorphic copy of a \roms.

{\roms}s\footnote{The Romans were famous for building straight roads, and we wanted a name for a two-dimensional object (a surface) with the property that between any two points there is a geodesic.} are totally geodesic  two-dimensional surfaces that are not built from one-dimensional objects; in this way
they are analogous to the two-dimensional flex locus studied in \cite{mmw17}.
We further show that {\roms}s arise independently in the study of the three-point interpolation problem in the polydisk; this is explained in  Subsection \ref{subsecbc}.

Understanding sets that are totally geodesic or \car is not just of geometric interest.
Questions \ref{qa1} and \ref{qa2} naturally arise when considering von Neumann or And\^o style inequalities, and when studying Nevanlinna-Pick interpolation. See Section \ref{sec:significance}
 where the significance of the main result is explained.

\section{Introduction}
\label{secb}

\subsection{Cartan pairs and retracts}

In this section we shall flesh out the details of the overview.
First let us define what we mean by a holomorphic function on a set that is not open.
Let $\B(a,r)$ be the open ball centered at $a$ of radius $r$.
\begin{defin}
\label{DEFb1} Let $V \subseteq \C^N$.
We say $\phi: V \to \C$ is holomorphic if for every $\mu \in V$ there exists $\vare >0$
and $f \in \hol(\B(\mu,\vare))$ so that $\phi |_{V \cap \B(\mu,\vare)} = f |_{V \cap \B(\mu,\vare)}$.
\end{defin}

\begin{defin}
\label{DEFb2}
Let $\O$ be a pseudoconvex domain in $\C^N$. We say $V \subseteq \O$ is an analytic subvariety if
for every $\mu \in V$ there exists $\vare >0$ and functions $f_1, \dots , f_k \in \hol(\B(\mu,\vare))$ such that
$V \cap \B(\mu,\vare) = \{ z \in \B(\mu,\vare) : f_j(z) = 0, \ 1\leq j \leq k \}$.
If $V$ is an analytic subvariety of $\O$, we shall call $(\O,V)$ a Cartan pair.
We say that $V$ is an algebraic subvariety of $\O$ if there is an algebraic set ${\mathfrak V} \subseteq \C^N$
so that $V = {\mathfrak V} \cap \O$.
\end{defin}

With Definition \ref{DEFb1} in hand, the definition of $d_V$ in \eqref{Eq1} 
and the definition of \car set make sense.
Let $(\O,V)$ be a Cartan pair, and let $\phi \in \hol(V)$.
A deep theorem of H. Cartan \cite{car51} says that if $\phi: V \to \C$ is a holomorphic 
function, then it extends to a holomorphic function $\Phi$ on all of $\O$.
But if there are  restrictions on the range of  $\phi$, it becomes an interesting geometric question
of when $\Phi$ can be chosen with the same restrictions.
We say $\O$ is an {\em extension set} for $V$ if every holomorphic function from 
$V$ to the unit disk $\D$ extends to a holomorphic function from $\O$ to $\D$.
Given $\Omega$, the question of which sets $V$ have  $\O$ as an extension set has been
extensively studied---see e.g. \cite{si75, agmcvn, ghw08, kn11b, dey17, maci19, kmc20, kz21} when $\Omega$ is the polydisk, and 
\cite{bs18, aly19, kmc19, br20, maci20, akmc22, akmc23} for other domains.

There is a finite version of being an extension set. For $n$ a positive integer, say  $\O$ is an $n$-extension set for $V$
if, for any $n$ points $\mu_1, \dots, \mu_n$ in $V$ and every $\phi \in \hol(V, \D)$ there is a $\Phi \in \hol (\O, \D)$ satisfying
$\Phi(\mu_j) = \phi(\mu_j) $ for each $j$. By a normal families argument, being an extension set is the same as being an $n$-extension set for every $n$.
A \car set is the same as a
 $2$-extension set.
To the best of our knowledge, there are no known examples of \car sets that are not extension sets.

\begin{defin}
\label{DEFb4}
A complex geodesic in $\O$ is a holomorphic function $f: \D \to \O$ that has a left inverse 
$\phi \in \hol(\O,\D)$, i.e. satisfies $\phi (f(z)) = z $ for every $z \in \D$.
\end{defin}
An equivalent definition of a complex geodesic is that $f$ is an isometry between $\D$ and $\O$ equipped with
their respective Kobayashi metrics (see Subsection \ref{subsecca}).
\begin{defin}
\label{DEFb5}
If $V$ is an analytic subvariety of $\O$, we say that it is totally geodesic if, for any pair of points $\mu,\nu$ in $V$, there is a complex geodesic  $f: \D \to \O$ passing through $\mu$ and $\nu$ whose range lies in $V$.
\end{defin}

Note that whether $V$ is a  \car set depends a lot on the choice of $\O$. Indeed, in \cite{akmc22} it was shown that if $V$ is fixed, one can always choose $\O$ so that it is an extension set, so a fortiori a \car set. 
Saying that $V$ is totally geodesic is more subtle. It means that any two points in $V$ lie on a Kobayashi geodesic $\phi$ of $\Omega$ that is entirely contained in $V$.

\begin{defin}
\label{DEFb6}
Let $(\O,V)$ be a Cartan pair.
We say $V$ is a retract of $\O$ if there exists $r \in \hol(\O,V)$ such that $r |_V = {\rm id}$.
\end{defin}

Consider the following 4 properties that a Cartan pair  $(\O,V)$ may have:

(R) \quad 
$V$ is a retract of $\O$.

(E) \quad $\O$ is an extension set for $V$. 

(TG) \quad
$V$ is totally geodesic in $\O$.

(C) \quad
$V$ is a \car set in $\O$.

It is straightforward to see that
(R) implies (E),  (TG) and (C); and of course (E) implies (C). If $\O$ is bounded and convex, then  (TG) implies (C) \cite[Prop. 6]{kz21}.
Sometimes the converses hold. For example, if $\O$ is the unit ball in $\C^N$ and
$V$ is an algebraic subvariety (or more generally  $V$ is relatively polynomially convex), it follows from the proof of \cite[Thm. 3.1]{kmc19} that (C) implies (R).
If $V$ is an algebraic subvariety of $\O$   and $\O$ is strictly convex, using techniques developed by Lempert \cite{lem81} and Huang \cite{hu94}, 
 the proof of \cite[Prop. 2.4]{kmc19} shows that (C) implies (TG) (and in dimension 2 this in turn implies (R)).
On the other hand, once one relaxes the requirement of strict convexity, there are examples of sets $\O$ in $\C^2$
with algebraic subvarieties $V$ that have (E) but are not retracts, or even smooth \cite{aly19, akmc22}.

Retracts of the polydisk were characterized by Heath and Suffridge \cite{hs81}; they are all graphs over lower dimensional polydisks.
\bt
\label{thmb1} \cite{hs81}
Let $V$ be a $\dimv$-dimensional retract of $\D^N$. Then, after a permutation of the coordinates, $V$ can be
written as
\[
V \= \{ (z, f(z) ) : z \in \D^\dimv \}
\]
where the function $f $ is in $\hol(\D^\dimv, \D^{N-\dimv})$.
\et

For the rest of the paper,  we restrict ourselves  to algebraic subvarieties $V$ of the polydisk $\D^N$,
and when we talk about \car sets or totally geodesic sets we mean algebraic subvarieties of the polydisk with those properties.

\subsection{Roman Surfaces}
\label{subsecb1}

For $b = (b_1, b_2, b_3) \in  \C^3 \setminus \{0\}$, define the set 
\begin{eqnarray}
\label{Eq4} \notag
\mathcal K_b & := & \{z\in \mathbb D^3:\ b_1 z_1 + b_2 z_2 + b_3 z_3 = \bar b_1 z_2 z_3 + \bar b_2 z_1 z_3 + \bar b_3 z_1 z_2\}
\\
\notag
&=& \{(z_1,z_2)\in D_{\kappa_\beta} :\ z_3 = \kappa_{b}(z_1, z_2)\}
\end{eqnarray}
where
\begin{eqnarray}
\label{eq4.1}
\kappa_b(z_1, z_2) &=& \frac{\bar b_3 z_1 z_2 - b_1 z_1 - b_2 z_2}{b_3 - \bar b_1 z_2 - \bar b_2 z_1} 
\\
\label{eq4.2} \notag
D_{\kappa_\beta} &=& \{ (z_1, z_2) \in \D^2: |\kappa_b(z)| < 1 \}.
\end{eqnarray}
We shall say that $b \in \C^3$ is {\em triangular} if the triple $(|b_1|,|b_2|,|b_3|)$ form the sides of a nondegenerate triangle (i.e. $2 \times \max ( |b_1|, |b_2|, |b_3|) < |b_1| + |b_2| + |b_3|$).
We shall write $\K$ for $\K_{(1,1,1)}$:
 \[
 \K = \{ (x,y,z) \in \D^3 : x + y + z = yz + xz + xy \} .
 \]

There are two possibilities for $\K_b$. If $b$ is not triangular, then after permuting the variables so
$|b_3| \geq |b_2| + |b_1|$ we get $\kappa_b : \D^2 \to \D$, and then $\K_b$ is a
 two-dimensional retract of $\D^3$.
Otherwise, they are all automorphic copies of $\K$, as proved by K. Maciaszek \cite{maci23}.
\bt
\label{thmb2} \cite{maci23}
If $b$ is triangular, then there is an automorphism $\Phi : \D^3 \to \D^3$ such 
that $\Phi(\K_b) = \K$.
\et

Let $\mathcal M$ denote the family of functions 
\[
\M \= \{\kappa_b: \text{ $b$ is triangular}\} .
\]
Note that every function $\a$ in $\M$ can be uniquely written in the form
\be
\label{Eq5}
\a(x,y) \= \tau \frac{ Ax + B y +xy}{1 + \bar{B}x + \bar{A} y } ,
\ee
where $|\tau| = 1$ and $(1,|A|,|B|)$ form the sides of a triangle.
We shall use the short-hand notation $[\tau; A,B]$ to denote a function in $\M$ that has the normalized form of
\eqref{Eq5}.
By  $\overline{\M}$ we mean the set of functions of the form \eqref{Eq5}  where $|\tau| = 1$ and $(1,|A|,|B|)$ form the sides of a possibly degenerate triangle, i.e. we allow the maximum to equal the sum of the other two.
If $1 = |A|+|B|$, then $\a \in O(\D^2, \D)$. If $1 + |A| = |B|$, and one changes variables from $(x_1, x_2, \a(x_1,x_2))$ to $(x_1, \beta(x_1, x_3), x_3)$, then $\b \in O(\D^2, \D)$. Note that unimodular multiples of the coordinate functions are in $\overline{\M}$.

We associate with $\a \in \M$ the set
\be
\label{Eq6}
\hatt(\a) = \{ \omega \in \T : | A+ B \o| < 1 \}.
\ee
The triple $(1,A,B)$ is triangular if and only if $\hatt(\a)$ is a non-empty arc in $\T$ whose closure is not all of $\T$
(Lemma \ref{lemx1}). It turns out that $\hatt(\a)$ plays a key r\^ole in the analysis of \car sets.

\begin{defin}
\label{defb7}
For $ j \in \I$, let $\a_j$  be given by 
\[
\a_j(x,y) \= \tau_j \frac{ A_jx + B_j y +xy}{1 + \bar{B_j}x + \bar{A_j} y }  .
\]
We say $\{ \a_j : j \in \I \}$ form an edge family if

(i) each $\a_j \in \M$ 

(ii)
  there are unimodular $\om$ and $\eta$ so that
\be
\label{eq6.1}
 A_j + B_j \om \= \eta   \qquad \forall\  j \in \I .
\ee
We shall call  $\o$ the anchor point of the family.
If $|\I| = 2$, we call the two functions an edge pair.
\end{defin}

\begin{remark}{\rm
It is important to note that \eqref{eq6.1} refers to $A_j$ and $B_j$, not $\frac{\partial \a_j}{\partial x} (0)  = \tau_j A_j$ 
and $\frac{\partial \a_j}{\partial y} (0) = \tau_j B_j$. Multiplying $\a_j$ by a unimodular constant will rotate
$\frac{\partial \a_j}{\partial x} (0)  + \frac{\partial \a_j}{\partial y} (0 ) \o$, but it will not change $A_j + B_j \om$.
}
\end{remark}

\begin{defin}
\label{defb8}
A \roms\ $\LL$ is a subset of $\D^{2+m}$ of the form
\be
\label{Eq7}
\LL \= \{ (x_1, x_2, \a_3(x_1, x_2),\dots,  \a_{m+2}(x_1, x_2) ) : (x_1, x_2) \in D_{\a_3} \cap 
\dots \cap D_{\a_{m+2}} \}
\ee
for some edge family $\{\a_3,\dots,  \a_{m+2}\}$. We call $m$ the height of $\LL$.
\end{defin}

Sometimes we shall need to work with limits of  {\roms}s.

\begin{defin}
For $ j \in \I$, let $\a_j$  be given by 
\[
\a_j(x,y) \= \tau_j \frac{ A_jx + B_j y +xy}{1 + \bar{B_j}x + \bar{A_j} y }  .
\]
We say  $\{ \a_j : j \in \I \}$ form a weak edge family if

(i) each $\a_j \in \overline{\M}$ 

(ii)
  there are unimodular $\om$ and $\eta$ so that, for all $j \in \I$ for which  $A_j B_j \neq 0$,
\be
\label{eq6.2}
 A_j + B_j \om \= \eta  .
\ee
A \wrs\ is a set of the form \eqref{Eq7} for some weak edge family 
where at least one of the functions $\a_j$ is in $\M$.
\end{defin}

Note that by definition, unimodular multiples of  coordinate functions are in any weak edge family.
That is because some unimodular multiple will satisfy \eqref{eq6.2}, and we should be free to multiply by any $\tau_j$  without changing this. The reason that {\wrs}s are required to contain at least one function in $\M$ is that otherwise they would be graphs over a bidisk, and hence two-dimensional retracts.

The {\roms}s of height $1$ are precisely the sets $\K_b$ where $b$ is triangular.
If some function $\a_m \in \ovm \cap O(\D^2, \D)$, then $\LL$ is a graph over the \roms\ $\LL'$ of height $m-1$ obtained by omitting $\a_m$. Since $\a_m: \D^2 \to \D$, the surface $\LL'$ is a \car set if and only if $\LL$ is.
So when it comes to classifying \car sets, it is natural to look at {\roms}s, not {\wrs}s.
However, permutation automorphisms of $\D^{m+2}$ preserve {\wrs}s not {\roms}s (see Example \ref{exx2.1}),
so we must consider them also. (Note that by definition, all Roman surfaces are normalized to  contain the point $0$.)

Let $\C_\infty = \C \cup \{ \infty \}$.
If a surface in $({\C_\infty})^\dimpol$ is the graph of functions $\a_j(x_1, x_2)$ where each $\a_j \in \ovm$, then
typically it cannot be written as a graph over a different pair of variables. 
Being able to write it as a graph over any pair of variables is equivalent to every pair of functions being an edge pair, and so this property is invariant under automorphisms 
of $\D^\dimpol$,
since every automorphism of $\D^\dimpol$ extends to an automorphism of $({\C_\infty})^\dimpol$.
\bt
\label{thmb3} Let $N \ge 4$.
Let $\a_3, \dots, \a_N \in \overline{\M}$, and let 
\be
\notag
\LL = \{ (x_1, x_2, \a_3(x_1, x_2), \dots, \a_N(x_1, x_2)) : 
 (x_1, x_2) \in D_{\a_3} \cap 
\dots \cap D_{\a_N} \}.
\ee
Let
\[
\LL^\sharp = \{ x \in (\C_\infty)^N : x_j = \a_j(x_1, x_2), \ 3 \le j \le N \}.
\]
Assume that no two of the $\dimpol$ functions $x_1, x_2, \a_3(x_1,x_2), \dots , \a_\dimpol(x_1, x_2)$ are unimodular multiples of each other.
Then the following are equivalent:

(i) For any pair $j,k \in \{ 1, 2, \dots, N \}$, the projection of $\LL^\sharp$ onto the $x_j-x_k$ plane is injective.

(ii) For any pair  $j,k \in \{ 1, 2, \dots, N \}$, there are functions $\g_n$ in $\overline{\M}$ for $n \in \{ 1, \dots , N \} \setminus \{ j,k\}$
so that $\LL$ can be written as 
\[
\LL = \{ (x_1, \dots, x_N) \in \D^N : x_n = \g_n(x_j, x_k), (x_j, x_k) \in D_{\g_1} \cap \dots \cap D_{\g_N} \} .
\]
(iii) For any pair $j,k \in \{ 3, \dots, N \}$, the functions $\a_j$ and $\a_k$ form a weak edge pair.

Moreover, if these conditions hold, then  either $\LL$ is a 
retract, and the functions $\{\a_3, \dots, \a_N \}$  form a weak edge family,  or $\LL$ is a \wrs\ and a graph over a \roms,
or $\LL$ is neither,
 in which case it is a graph over a \roms\ of height one.
\et
We prove Theorem \ref{thmb3} in Section \ref{sece}.
An example where $\LL$ is not a \wrs\ is given in Example \ref{exx8}.

In Theorem \ref{thme3}, proved in Section \ref{sece}, we prove that {\roms}s are  \car sets.
We do this in steps. First, we show that they are totally geodesic in the infinitesimal sense---at every point, and for every vector in the tangent space at that point, there is a geodesic tangent to that vector. Then we use a topological and geometric argument to show that the {\roms} is totally geodesic---given any pair of points, there is a geodesic connecting them. From this, we conclude that they are also \car sets.
\begin{theorem}
\label{thme3}
{\roms}s are totally geodesic, and are \car sets.
\end{theorem}
We define Lempert sets in Definition \ref{defc1.1}; 
they are sets for which the conclusion of Lempert's Theorem \ref{thmc0} holds.
Lempert domains that are not biholomorphic to convex domains are hard to come by.
A consequence of Theorem \ref{thme3} is that {\roms}s are Lempert sets and their
 two-dimensional bases are Lempert domains.
\begin{theorem}
\label{thme5}
Let $\LL$ be a \roms\ as in \eqref{Eq7}. Then $\LL$ is a Lempert set, and 
 $ D_{\a_3} \cap 
\dots \cap D_{\a_{m+2}}$ is a two-dimensional Lempert domain.
\end{theorem}

\begin{remark}
It follows from Theorem \ref{thma1}
that all basic \car sets are either disks or {\roms}s. However, it is not true that all {\roms}s are
basic---they may be a graph over a \roms\ in a lower dimensional polydisk, as in Example \ref{exx2.1}.
\end{remark}

\subsection{More on {\roms}s}
\label{subsecbc}

{\roms}s appeared implicitly in \cite{kosthree}. Although they were not identified or defined there as geometric objects in their own right, they already appeared naturally and can now be recognized as the same geometric objects. We explain this connection below. 

Let $\mu_1, \dots, \mu_n$ be points in $\D^N$, and $w_1, \dots, w_n$ be points in $\D$.
We say that they are the data for an $n$ point {\em extremal Pick problem} if there exists $F \in \hol(\D^N, \D)$
that interpolates, i.e.  satisfies $F(\mu_j ) = w_j$ for each $1 \le j \le n$, and if for no $r < 1$ is there a function in $\hol(\D^N, r \D)$ that interpolates. We say the problem is {\em nondegenerate} if no $n-1$ point subproblem is extremal.  An extremal problem is called {\em strictly $k$-dimensional} if there is a solution that depends on $k$ variables, but not one that depends on fewer variables.

By applying automorphisms, one can always normalize the problem to assume that $\mu_1 = 0$ and $w_1 = 0$.

\bt
\label{thmb5}
Let $\dimpol \ge 3$.
A strictly 3-dimensional nondegenerate 3 point extremal Pick problem in $\D^\dimpol$ never has a  unique solution. After
an automorphism of $\D^\dimpol$, the subset of $\D^\dimpol$  on which  all solutions coincide is
the Cartesian product of a weak \roms\ of height $m$ and $\D^{\dimpol-m-2}$. Moreover, every 
set of this form
is the uniqueness set for some 
strictly 3-dimensional nondegenerate 3 point extremal Pick problem.
\et
We prove Theorem \ref{thmb5} in Section \ref{secn}.
The case $N=3$ was proved in \cite{kosthree} and \cite{maci23}.

 \subsection{Context and Motivation}\label{sec:significance}
Our main result is that, among algebric subvarieties of the polydisk,  every basic \car  set is either a disk or a \roms. Moreover, we show that every \car set is totally geodesic.
Although this is a geometric classification theorem,  it is actually at the nexus of many areas of analysis: von Neumann inequalities, spectral sets, extension sets, Pick interpolation, invariant distances and complex geodesics.

Von Neumann inequalities are inequalities of the form 
\be
\label{eqb1}
\sup_{T \in {\mathcal C}} \| p(T) \| \le \| p \|_V ,
\ee
where $p$ is a polynomial in $\dimv$ variables, $T$ is a commuting $d$-tuple of operators that lie in a certain class ${\mathcal C}$,  $V$ is some  subset of $\C^d$, and the right-hand side means the supremum norm on $V$. The original inequality of von Neumann \cite{vonN51} was when ${\mathcal C}$ consisted of the contractions on a Hilbert space (the operators of norm at most $1$) and $V$ was the unit disk. 
When an inequality of the form \eqref{eqb1} holds, we say $\overline{V}$ is a spectral set for $T$.

The importance of the inequality is that a matrix or operator norm is bounded by a scalar norm.
And\^o extended von Neumann's inequality to the class of commuting pairs of contractions, with $V$ becoming the bidisk \cite{and63}.
The impetus behind \cite{agmcvn} was trying to replace $V$ by some subset of the bidisk, and considering ${\mathcal C}$ to be pairs of commuting contractions that in some sense ``lived on'' $V$; this was developed further in \cite{agmc_dv}. The connection between von Neumann inequalities and describing extension sets for the polydisk is discussed in \cite[Sec. 8]{kmc20}.
The problem on the symmetrized bidisk has been considered in many places, such as  \cite{ps14,dks25}. Finding classes of commuting $d$-tuples  for which
\eqref{eqb1} holds on the polydisk has been done in \cite{gkvw09, bds20}. The sharpest known inequalities that apply to all $d$-tuples of commuting contractions, with a constant on the right-hand side, have been found in \cite{har25}.

In turn, the question of understanding \eqref{eqb1} for subsets of $\D^2$ came out of the following observation: when one solves an extremal Pick problem on $\D^N$, for $N \ge 2$, there is a variety $V$ on which the solution is unique, and has norm $1$. It extends off this variety to a function that has norm one on the whole polydisk. Does this latter extension exist because of the geometry of the variety, or the particular nature of the solution? The answer in the bidisk is that,  not all holomorphic functions from $V$ to $\D$  extend to functions from $\D^2$ to $\D$, since $V$  is never a \car set for strictly two-dimensional problems, and hence it is the special nature of the solution that allows the extension.
For two point problems, however,  $V$ will either be a geodesic disk or the whole bidisk, and therefore will have the extension property. 

This question can be asked again for 3 point nondegenerate strictly 3-dimensional problems in $\D^N$. 
By Theorem \ref{thmb5}, the uniqueness set is a weak \roms\ (times a polydisk factor), and the solution has multiple norm-preserving extensions off  this set. We do not know whether {\roms}s have the extension property for all functions, or just special ones.

The realization that all known proofs that a set did not have the extension property actually proved it was not a \car set, and the identification of \car sets as interesting objects to study, was in \cite{kz21}, where the connection with totally geodesic sets was made. These objects have been intensely studied, starting with 
Kobayashi and Lempert \cite{ko67, lem81, lem84}, and continuing to the present
 \cite{bgnt24, brga18, mark18, nt26}.

\section{Background results}
\label{secc}

\subsection{Kobayashi and \car metrics}
\label{subsecca}

A standard reference for invariant distances and metrics is the book \cite{jp13} by M. Jarnicki and P. Pflug.
\begin{defin}
\label{defc1}
Let $(\O,V)$ be a Cartan pair. For $\mu,\nu \in V$, define the Lempert function
\be
\label{Eq8}
l_V (\mu,\nu)  \ := \ \inf \{ \rho(a, b) : \ \exists\  f \in \hol(\D,V),\ f(a) = \mu, f(b) = \nu \} .
\ee
If no such $f$ exists, we define $l_V (\mu,\nu) = +\i$.
The Kobayashi pseudodistance is the largest pseudodistance  $k_V$ that is less than or equal to $l_V$.
The infinitesimal version, called the Kobayashi-Royden pseudometric,  is defined for $\mu \in V$ and
$X$ a tangent vector at $\mu$ 
\be
\label{Eq9}
\notag
l_V(\mu; X) \ := \ \inf \{ |t| : \ \exists f \in \hol(\D,V), f(0) = \mu, t f'(0) = X \} .
\ee
\end{defin}
The set $\O$ is called Kobayashi hyperbolic if the Kobayashi pseudodistance is a distance (i.e.  distinct points are a non-zero distance apart). Note that it is an immediate consequence of the definitions that if $\phi : V \to U$ is a holomorphic function, then
for any pair of points $(\mu,\nu)$ in $V$ we have
$d_U(\phi(\mu), \phi(\nu)) \le d_V(\mu, \nu)$
and $l_U(\phi(\mu), \phi(\nu)) \le l_V(\mu, \nu)$.

 L. Lempert proved the following celebrated theorem. An operator theoretic proof that used Pick interpolation was found later  by J. Agler \cite{ag91}.
\bt
\cite{lem81,lem84}
\label{thmc0}
Let $\O$ be a bounded convex domain. Then  $\O$ is Kobayashi hyperbolic and 
$d_\O = l_\O = k_\O$.
Moreover, for any $\mu, \nu \in \O$  there exists a complex geodesic 
$f: \D \to \O$ whose range contains $\mu$ and $\nu$.
\et

\begin{defin}
\label{defc1.1}
Let us say that a set $U$ is a Lempert set if $d_U = l_U$.
If in addition $U$ is a domain, i.e.  a connected open set, then we shall call it a Lempert domain.
\end{defin}

For product domains, we will use the following theorem of Jarnicki and Pflug, that the \car metric for the product is the maximum of the metric on each factor.
\bt
\label{thmc1}
Let $V_1$ and $V_2$ be connected complex analytic spaces, and $\mu = (\mu_1, \mu_2)$ and $\nu = (\nu_1, \nu_2)$ be points
in $V_1 \times V_2$.
Then 
\begin{eqnarray}
\label{Eq10}
d_{V_1 \times V_2} (\mu,\nu) &\=&  \max ( d_{V_1}(\mu_1, \nu_1), d_{V_2}(\mu_2, \nu_2))\\
l_{V_1 \times V_2} (\mu,\nu) &\=&  \max ( l_{V_1}(\mu_1, \nu_1), l_{V_2}(\mu_2, \nu_2)) .
\label{Eq11}
\end{eqnarray}
\et
Equality \eqref{Eq10} is proved in \cite{jp99}. Equality \eqref{Eq11} is proved in 
\cite[Prop. 3.7.1]{jp13} for domains, but the proof extends immediately to complex analytic spaces.
It follows that the product of Lempert sets is a Lempert set.

Note that if $\O$ is taut, as long as there is some function $f \in \hol(\D,V)$ that maps
$a$ to $\mu$ and $b$ to $\nu$, then a normal families argument shows that the infimum in 
\eqref{Eq8} is attained, though a priori it might map into $\overline{\O}$.
If $\O$ is the polydisk, looking at each component function separately we
see that each one must map into $\D$, so the infimum is actually attained with a function $f$
whose range is in $V$.
 If  $\phi \in \hol (V, \D)$, then $\phi \circ f$ is a contraction from
$(\D, \rho) $ to itself, so $d_V(\mu,\nu) \le l_V(\mu,\nu)$.
So if $V$ is in addition a Lempert set, then $d_V = k_V = l_V$, so there is an extremal $\phi$
which is a left inverse of $f$. Therefore
$V$ is totally geodesic, as defined in Definition \ref{DEFb5}. 
So we have:
\begin{prop}
\label{PROc2}
Let $(\D^N,V)$ be a Cartan pair. Assume that $V$ is a Lempert set.
Then $V$ is totally geodesic.
\end{prop}

The following result was proved in \cite{kz21}; we include the proof for completeness.
\begin{prop}
\label{propc3}
Let $\O$ be bounded and convex, and let $(\O,V)$ be a Cartan pair. If $V$ is totally geodesic, then it is a \car set.
\end{prop}
\bp
Let $\mu, \nu \in V$, and let $\phi : V \to \D$ be a \car extremal. By Theorem \ref{thmc0} there is a holomorphic $f: \D \to V$ such that $f(a) = \mu, f(b) = \nu$ for some $a,b \in \D$, and $f$ has a left inverse $\psi$.
Therefore
\[
d_V(\mu,\nu) = \rho(\phi(\mu), \phi(\nu)) = \rho(\phi(f(a)), \phi(f(b))) \le \rho(a,b) .
\]
We also have
\[
d_\O(\mu,\nu) \ge \rho(\psi(\mu), \psi(\nu)) = \rho(\psi(f(a)), \psi(f(b))) = \rho(a,b) .
\]
As $d_V(\mu, \nu) \ge d_\O(\mu, \nu)$, they must be equal.
\ep

\subsection{Rudin's theorems}
\label{subseccb}

We shall write $\T$ for the unit circle in $\C$. For a monomial $z_1^{i_1} \cdots z_N^{i_N}$, we define its 
degree to be $i_1 + \dots + i_N$. The degree of a polynomial is the maximum of the degrees of the monomials
with non-zero coefficients. The degree of a rational function is the maximum of the degree of its numerator and denominator,
provided they have no common factors.
If ${\mathcal F}$ is a set of functions 
on a set $X$, a set $E \subseteq X$ is called a determining set for ${\mathcal F}$ if
$f = g $ on $E$ implies $f=g $ on $X$ for any functions $f,g \in {\mathcal F}$.

The following notation will be used in the next theorem.
For each $w \in \C^N$, let $\Delta_w = \{ \l w : \l \in \overline{\D} \}$.
For each $w \in \C^N$, if $f$ is a function defined on a subset of $\C^N$,
the function $f_w(\l) = f(\l w)$ is a function of the single complex variable $\l$.
We write $A(\D)$ to denote the disk algebra, that is the functions that are holomorphic on
$\D$ and continuous on $\overline{\D}$.

W. Rudin proved the following theorem \cite[Thm. 5.2.2]{rud69}.
\bt
\label{thmc2}
Let $\Omega$ be a connected neighborhood of $0$ in $\C^N$.
Let $E$ be a compact subset of $\C^N$ that does not contain $0$ and is a determining set for polynomials.
Let $Y$ be the set of complex valued functions  $f$ defined on $\Omega \cup \cup_{w \in E} \Delta_w$
that are holomorphic on $\O$ and satisfy $f_w \in A(\D)$ for every $w \in E$.

Assume there is an integer $k \geq 0$ such that $f_w$ is a rational function (of one variable) of degree $k$ for every $w \in E$.
Then $f$ is a rational function of $N$ variables of total degree $k$.
\et

Let $p$ be a polynomial in $N$ variables; let the maximum degree in
each variable $z_j$ be $d_j$.
 Define $\revp(z) = z_1^{d_1} \dots z_N^{d_N} \overline{p(\frac{1}{\overline{z}})}$.
An inner function on $\D^N$ is a bounded holomorphic function whose radial boundary values have modulus $1$ a.e. on
$\T^N$.

Rudin characterized rational inner functions on the polydisk in \cite[Thm. 5.2.5]{rud69}.
His theorem can be extended to rational functions with poles inside $\D^N$. This is well-known, but since we have been unable to locate a reference, we include a proof.
\bt
\label{thmc4}
Every rational  function $u$ on $\D^N$ that is unimodular a.e. on $\T^N$  is of the form
\be
\label{Eq13}
\notag
u(z) = z^I \frac{ \revp(z)}{p(z)}
\ee
for some polynomial $p$  and some monomial $z^I$, with $I \in \Z^N$.
\et
\bp
Indeed, write $ u =  \frac{q}{p}$,
where $p,q$ are polynomials with no common factors.
Then 
\be
\label{Eq14}
\frac{1}{\overline{u(\frac{1}{\bar z})}}
= 
z^J \frac{\revp(z)}{\revq(z)}
\ee
where $J_\dimv$ is the maximum degree in $z_\dimv$ of $q$ minus its degree in $p$.
Since $u$ is unimodular on $\T^N$, \eqref{Eq14} and $u$ agree a.e. on $\T^N$, and hence are equal as rational functions.
So there are monomials $z^{I_1}, z^{I_2}$ with non-negative powers so that
\be
\label{Eq15}
z^{I_1} p(z) \revp(z) = z^{I_2} q(z) \revq(z) .
\ee
Factor $p$ and $q$ into irreducible factors. There is a corresponding factoring for $\revp$ and $\revq$.
Comparing both sides of \eqref{Eq15} we see that up to a monomial factor, $q = \revp$.
\ep

\subsection{Automorphisms of the polydisk}
\label{subseccd}

For each $\mu \in \D$, let 
$\vp_\mu$ be the M\"obius map
\beq
\vp_\mu(z) &=& \frac{z-\mu}{1-\bar \mu z}.
\eeq
Its inverse is $\vp_{-\mu}$.
There are three sorts of automorphism of $\D^N$, which we shall call permutations, rotations, and M\"obius automorphisms.
A permutation is a map of the form
\[
P_\sigma :(z_1, z_2, \dots, z_N) \mapsto (z_{\sigma(1)}, z_{\sigma(2)}, \dots , z_{\sigma(N)}) 
\]
for some permutation $\sigma$ of $\{ 1, \dots , N\}$.
A rotation is a map of the form
\be
\label{eq18}
\notag
R_\theta: (z_1, z_2, \dots, z_N) \mapsto (e^{i\th_1} z_1, e^{i\th_2} z_2, \dots, e^{i\th_N} z_N ) .
\ee
A M\"obius automorphism is a map of the form
\be
\label{eq19}
\notag
\Phi_a : 
(z_1, z_2, \dots, z_N) \mapsto (\vp_{a_1}(z_1), \vp_{a_2}(z_2), \dots, \vp_{a_N}(z_N)) ,
\ee
where $a \in \D^N$.

Every automorphism of $\D^N$ can be uniquely written as the composition of a M\"obius automorphism, a permutation and a rotation \cite[Cor. to 7.3.3]{rud69}.
Since the property of being a \car set is preserved under automorphisms, we shall frequently apply automorphisms to simplify our notation.
By a {\em small perturbation} we will mean an automorphism $R_\theta \Phi_a$ where $a$ is close to $0$.

Suppose that $a \in V$ and that near $a$ there is a parametrization of $V$ of the form
$\{ (x_1, \dots, x_d, \a_{d+1}(x), \dots , \a_N(x) )\}$ where $x = (x_1, \dots, x_d)$ lies in a 
neighborhood $U$ of $(a_1, \dots, a_d)$ and the functions $\a_j$ are holomorphic on $U$.
Let $\tilde V = \Phi_a (V)$. Then $0 \in \tilde V$, and near $0$ the set $\tilde V$ has a parametrization
$\{ (y_1, \dots, y_d, \tilde \a_{d+1}(y), \dots , \tilde \a_N(y)) \}$ where $y_i = \vp_{a_i}(x_i)$ and
\be
\label{eq20}
\tilde \a_j(y_1, \dots, y_d) \ = \vp_{a_j} \circ \a_j(\vp_{-a_1}(y_1), \dots, \vp_{-a_d}(y_d) ) .
\ee
By the chain rule, 
\be
\label{eq21}
\frac{\partial \tilde \al_j}{\partial y_i} \Big|_0 \=  \frac{1-|a_i|^2}{1-|a_j|^2}\  \frac{\partial \alpha_j}{\partial x_i}\Big |_a
\ee

\subsection{Theorems of Chevalley and Tworzewski}
\label{subsecce}

We shall use a theorem of Chevalley; for a proof in this form see e.g. \cite{mum76}.
If $X \subset \C^N$ we shall write $\pi_m(X)$ for the projection of $X$ onto the first $m$ coordinates.
\bt
\label{thmc7}
Let $X$ be an irreducible algebraic set in $\C^N$ of dimension $k$. Then for each $m < N$ the set
$\overline{\pi_m(X)}$ is an irreducible algebraic subset of $\C^m$ of dimension less than or equal to $k$.
\et

\begin{defin}
\label{defc2}
Let $\O$ be a connected open set in $\C^N$, and $f \in \hol(\O)$.
Say $f$ is a Nash function at $x_0 \in \O$ if there exists an open neighborhood  $U$ of $x_0$ 
in $\O$ and a polynomial $p \in \C[z_1, \dots, z_{N+1}]$ that is not identically zero and satisfies
$p(x,f(x)) = 0$ for all $x \in U$. Say $f$ is a Nash function on $\O$ if it is a Nash function at every point.
A function $F \in \hol(\O, \C^m)$ is a Nash mapping if each component is a Nash function.
\end{defin}

The significance of Nash functions for us is that functions that parametrize algebraic sets are Nash functions.
The following result was proved by P. Tworzewski \cite[Prop. 1.6]{tw90}.
\begin{theorem}
\label{thmc8}
Let $\O$ be an open connected set in $\C^N$, and $F: \O \to \C^m$ a holomorphic map.
Then the following are equivalent.

(i) $F$ is a Nash mapping.

(ii)  There exists an $N$-dimensional irreducible algebraic subset $X \subset \C^{N+m}$ so that
$\{ (x, F(x) ) : x \in \O \} \subset X$.

(iii) For each basis of $\C^m$, all components of $F$ are Nash functions.
\et

\subsection{Extension sets}

Extension sets were studied by N. Sibony in \cite{si75}, where he gave an example of
a pseudoconvex domain $V \subsetneq \D^2$ with the property that $\D^2$ was an extension set for $V$.
In \cite{agmcvn} it was shown that if one added an extra regularity condition to $V$, namely that it be relatively polynomially convex, then if $\D^2$ were an extension set $V$ would have to be a retract.
Extension sets have been widely studied; see e.g.  \cite{ghw08, kn11b, dey17, maci19, kmc20, aly19, br20, maci20, kz21,  akmc22, akmc23}.
In all cases where it was shown that the extension property forced $V$ to be a retract, examination of the proof shows
that only being a \car set was used. 


\section{Outline of the proof of Theorem \ref{thma1}}
\label{secd}

\subsection{Cartoon Version}

The key geometric idea is the following.
Let $V$ be a \car set in $\D^\dimpol$ containing $0$, and suppose that near $0$ it can be parametrized as $\{ (x, \a(x) ) \}$ where $x \in \C^\dimv$ and $\a$ takes values in $\C^{\dimpol - \dimv}$.
Define 
\[
T(\a) = \{ \om \in \T^\dimv : |\sum_{i=1}^\dimv \frac{\partial \a_j}{\partial x_i} \Big|_{0} \om_i | < 1 , \quad
\forall \ j > \dimv \}.
\]
Then for any $\o \in T(\a)$, each function $\l \mapsto \a_j( \om_1 \l, \dots, \om_\dimv \l)$ 
extends from a neighborhood of $0$ to a holomorphic function from $\D$ to $\D$.
This is Lemma \ref{lemd2}.
Combining this with Schwarz's Lemma, if $\om$ is in the boundary of $T(\a)$ and
$|\sum_{i=1}^\dimv \frac{\partial \a_j}{\partial x_i} \Big|_{0} \om_i | = 1$,
then $\a_j( \om_1 \l, \dots, \om_\dimv \l) = \eta \l$ for some $|\eta| = 1$.
If $\alpha_j$ agrees with an affine function  in a lot of directions, this places big restrictions on what $\a_j$ can be,
and if it persists under perturbations, we eventually show this forces $\a_j$ to be in $\M$  (Lemmas \ref{lemd3}, \ref{lemd4}, \ref{lemd5}).

\subsection{Grown-up Version}

$ (iii) \Rightarrow (i)$: We must first prove that {\roms}s are indeed \car sets.
We do this in Theorem \ref{thme3}, by first proving that they are totally geodesic, and then using Proposition \ref{propc3}
to conclude that they are also \car sets.
It then follows from Theorem \ref{thmc1} 
that any  product of disks and {\roms}s  is a \car set.
It is  immediate that any graph over such a set is also a \car set.

$(iii) \Rightarrow (ii)$:
In Theorem \ref{thme5} we strengthen this result to show that {\roms}s are Lempert sets. Using 
 Proposition \ref{PROc2} we conclude that any graph over a product of disks and {\roms}s is 
totally geodesic.

$(ii) \Rightarrow (i)$: This follows from Proposition \ref{propc3}.

$(i) \Rightarrow (iii)$:
 We shall fix $\mathfrak V$ a $\dimv$-dimensional algebraic set in $\C^N$, with $2 \le \dimv \le N-1$, and let $V = {\mathfrak V} \cap \D^N$.
We assume that $V$ is a \car set for $\D^N$, and must show that after an automorphism of $\D^N$ it has the form described in Theorem \ref{thma1}.

Let 
\be
\label{Eq22}
\notag
U_b \ = \D^{m_1}  \times \LL_1 \times \dots \times \LL_{m_2} ,
\ee
where each $\LL_j$ is a \roms\ of height $h_j$.
Let $m_3 = N - m_1 - \sum_{j=1}^{m_2}(2+ h_j)$.
Let $f=(f_j)$ be holomorphic functions from $U_b $ to $\D$, for $ 1 \leq j \leq m_3$.
Let
\be
\label{Eq23}
W = \{ (x, f(x) ) : x \in U_b \} ,
\ee
where if $m_3 = 0$ we understand that $W = U_b$. Our goal is to show that there is an automorphism of $\D^N$ that maps $V$ to some $W$ of this form, for some choice of $m_1,  \LL_1, \dots, \LL_{m_2}$ and functions $f$.

\vskip 10pt

The singular points of $\mathfrak V$ will be a set of lower dimension, so we can choose a point $a \in V$ that is
not singular, and applying the M\"obius automorphism $\Phi_a$ to $\D^N$ we get a new
algebraic subvariety that is still a  \car set and for which $0$ is a regular point.
After applying a permutation automorphism, we can assume that near $0$ we can write $V$ as 
a graph over the first $\dimv$ coordinates. This means there are holomorphic functions $\alpha_j$ 
defined on a neighborhood $U$ of $0$ in $\C^\dimv$ so that for some $\vare > 0$ we have
\be
\label{Eq24}
 V \cap \B(0,\vare) \= \{ (x_1, \dots, x_\dimv, \alpha_{\dimv+1} (x_1, \dots, x_\dimv), \dots , \alpha_N(x_1, \dots, x_\dimv)) : x \in U \} .
 \ee

If $f$ is the germ of a holomorphic function defined on a neighborhood of $0$ in $\C^\dimv$, 
define
\be
\label{Eq25}
\notag
T(f)
 \ := \ \{ \om  \in \T^{\dimv} : | \om_1 f_{x_1}(0) + \dots + \om_{\dimv-1} f_{x_{\dimv-1}}(0) + \om_\dimv  f_{x_\dimv}(0) | < 1 \}.
\ee
Since points in $T(f)$ are invariant under the multiplication action of $\T$, we can write them projectively as
$\om =  (\om_1 : \om_2 : \dots : \om_\dimv)$. It will often be convenient to normalize by choosing one of the coordinates to be $1$.  When $d=2$, we defined $\hatt(f)$ in Section \ref{secb} with the normalization $\om_1 = 1$, so there is only one number to consider. As we saw in that section, $T(f)$ contains a lot of information, especially if we also consider $T(\tilde f)$ after we apply M\"obius automorphisms.

\begin{remark}\label{remd1}
{\rm
Assume $T(f)$ is non-empty. If $T(f)$ is not dense in $\T^d$, then
$bT(f) = b\overline{T(f)}$ is a $(\dimv -1)$-dimensional subset of $\T^d$ (see Lemma \ref{lemh2}).
But if $T(f)$ is dense, then $b\overline{T(f)}$ is empty, and the boundary of $T(f)$ may be  anything from $1$ to $(\dimv -1)$-dimensional---see Example \ref{exx6}.}
\end{remark}

\begin{defin}
\label{defd1}
Suppose that $0$ is a regular point of $V$. We say $V$ has a good parametrization at $0$ if
locally $V$ has the form \eqref{Eq24} and
$\bigcap_{\dimv+1 \leq j \leq N} T(\al_j) \neq \emptyset$.
We say $V$ has a good $(x_{i_1}, \dots, x_{i_d}) $ parametrization at $0$ if, after permuting the variables so that the first $d$ are  $(x_{i_1}, \dots, x_{i_d}) $, then the permuted set has a good parametrization at $0$.

If $a$ is a regular point of $V$, we say $V$ has a good parametrization
at $a$ if $\Phi_a (V)$ has a good parametrization at $0$.
\end{defin}
Notice that having a good parametrization is an open condition, so will be preserved by small M\"obius perturbations.
We shall write 
\[
T(\a) = \bigcap_{\dimv+1 \leq j \leq N} T(\al_j) .
\]
Geometrically,
\be
\label{eqtt}
T(\a) = \pi_{1,\dots, \dimv}  T_0 (V) \cap ( \T^\dimv \times \D^{\dimpol - \dimv}) ,
\ee
where $T_0 (V)$ is the tangent space at $0$ to $V$, and $\pi_{1,\dots, \dimv}$ is the projection onto 
the first $d$ coordinates.

{\bf Step 1:} Assume that $V$ is irreducible as an analytic subvariety of $\D^N$.
In particular, this means $V$ is connected.
Observe that without loss of generality, we can assume that we never have $\alpha_j (x)  = \eta \alpha_i(x)$
for  $i \neq j$ and $\eta \in \C$.
Otherwise,  permuting so that $j = N$ and $|\eta| \leq 1$,  we can consider
the set $V'$ in $\D^{N-1}$ whose germ at $0$ is given by the functions $\alpha_1, \dots, \alpha_{N-1}$.
This will be a \car set in $\D^{N-1}$. If $V'$ has the form \eqref{Eq23}, then so does $V$ (by appending $\eta x_i$ to the functions $f$).
 Similarly, we can assume that we never have $\a_j(x) = \eta x_i$.
So we will make the standing assumption:
\begin{assumption}
\label{ass1}
Let $a$ be a regular point for $V$. Applying a M\"obius automorphism to move $a$ to $0$, and a permutation automorphism so that $V$ 
is locally given by \eqref{Eq24}, then the following hold:
\begin{itemize}
\item
For no pair $i \neq j$ in $\{ \dimv+1, \dots, N \}$ do we have a complex number $\eta$ so that
$\alpha_j = \eta \alpha_i$.
\item
 For no $i \in \{ 1, \dots , \dimv \}$ is there a $j$ in $\{ \dimv+1, \dots, N\}$ and
a complex number $\eta$ so that $\a_j (x) = \eta x_i$.
\item
For no $j$ in $\{ \dimv+1, \dots, N\}$ is $\a_j$ constant.
 \end{itemize}
\end{assumption}

Under this assumption, we show that there is always a good parametrization, after we permute the coordinates and
do a small M\"obius perturbation.  Note that the following lemma does not require $V$ to be a \car set.
We prove it in Section \ref{secf}.
\begin{lemma}
\label{lemd1}
Let $V$ be a $\dimv$-dimensional algebraic set in $\D^N$ that has $0$ as a regular point.
Assume that Assumption \ref{ass1} holds for $V$ at $0$.
After a permutation automorphism and a M\"obius perturbation, we can assume that $V$ has a good parametrization at $0$.
\end{lemma}

{\bf Step 2:}  The importance of Lemma \ref{lemd1} comes from the next lemma, which we prove in Section \ref{secg}.
For $\om \in \T^{\dimv}$  and $\l \in \C$, write
\be
\label{Eq26}
\notag
\omegadotlambda \ := \ ( \om_1 \l, \dots, \om_{\dimv} \l) .
\ee
\begin{lemma}
\label{lemd2}
Suppose $0$ is a regular point for $V$.
If $\om$ is in $\overline{T(\alpha)}$, then  for each $j$ the function $\l \mapsto
\alpha_j(\omegadotlambda) $ extends to a holomorphic map from $\D$ to $\D$ that is continuous from
$\overline{\D}$ to $\overline{\D}$.
Moreover, this extension is continuous on $W = \{ \omegadotlambda : \lambda \in \overline{\D}, \ \omega \in \overline{T(\alpha)} \}$ except possibly on a set of dimension $d-1$.
\end{lemma}

If $\overline{T(\a)} = \T^{\dimv}$
it  follows from Lemma \ref{lemd2} that every $\a_j$ extends to a holomorphic function from $\D^\dimv$ to $\D$, and hence we are in the situation of \eqref{Eq23} with $m_2 =  0$.  To see this, note that since $\alpha_j(\omegadotlambda) $ is holomorphic on $\D$ for every $\om$,
the power series for $\a_j$ converges in $\D^\dimv$. Moreover its limit at each point has modulus less than one. So in this case, 
$V$ is a retract of $\D^N$. The same applies if  $\overline{T(\tilde \a)} = \T^\dimv$ for 
any $a \in V$ and $\tilde \a$  given by
\eqref{eq20}.

If $V$ can be written as $V = \{ (x, f(x)): x \in V' \}$ for some $f \in O(V' ,\D)$, then $V$ will be a \car set if and only if $V'$ is. So we can assume that this does not happen, and 
 we  continue under the following assumption:
\begin{assumption}
\label{ass2}
 $V$ has a good parametrization at $0$ and
$\overline{T(\tilde \a)}$ is a proper non-empty subset of $\T^{\dimv}$ for every $a$ in a neighborhood of $0$ in $V$, where $\tilde \a$ is given by \eqref{eq20}. Moreover, there is no automorphism of $\D^\dimpol$ so that
$\tilde V = \{ (x, f(x)): x \in V' \}$ for some $f \in O(V' ,\D)$ and an algebraic set 
$V'$ in $\D^{\dimpol -1} $.
\end{assumption}

\vskip 10pt
 
{\bf Step 3:}
Let us look at points in the boundary of ${T(\a)} = \cap {T(\a_j)}$. They lie in either the boundary of a single ${T(\a_j)}$ or
of multiple ones. 

First let us address the single case. Assume for some $j$ 
 that there is a point $\om$ in the boundary of  ${T(\a)}$ that is in  the boundary of $ \overline{T(\al_j)}$ and in the interior of
all the other $T(\al_\ell)$ for $\ell > \dimv$. Then we show that $\alpha_j$ is in $\M$,
 and in particular
it depends on only two variables. We prove the next lemma in Section \ref{sech}.
\begin{lemma}
\label{lemd3}
Suppose $\om \in b \overline{T(\alpha_j)} \cap b {T(\a)} $ and moreover $\om \in \cap_{\dimv< \ell \le N, \ \ell \neq j} T(\alpha_\ell)$.
Then $\alpha_j \in \M$.
\end{lemma}

Note that the lemma is false if we just assume  $\om \in b {T(\alpha_j)} \cap b {T(\a)} $---see Example \ref{exx6}.
In Example \ref{exx7} we show the lemma is false even if 
 we assume 
$bT(\tilde \a_j) \cap_{\ell \neq j} T(\tilde \a_\ell)$ is non-empty for every $\tilde \a_j$.

{\bf Step 4:}
If $\om$ is in the boundary of several different $\overline{T(\a_j)}$, and this remains true under perturbation, then
they form an edge family when $\dimv=2$.

\begin{lemma}
\label{lemd4}
Let $\dimv=2$.
Let $\ell \ge \dimv+2$.
Suppose there exists $\vare > 0$ so that, for all $a \in V \cap B(0,\vare)$, if we write as in \eqref{eq20}
\[
\tilde \a_j(y_1, \dots, y_\dimv) \ = \vp_{a_j} \circ \a_j(\vp_{-a_1}(y_1), \dots, \vp_{-a_\dimv}(y_\dimv) ) ,
\]
then $ b \overline{T( \tilde \al_{\dimv+1}) }\cap \dots \cap b \overline{T(\tilde \al_{\ell}) }\cap T(\tilde \al_{\ell+1}) \cap \dots \cap T(\tilde \al_N)$ is non-empty.
 Then $\al_{\dimv+1}, \dots , \al_{\ell}$ are in $\M$,
 and 
  they form an edge family.
\end{lemma}
We prove 
Lemma \ref{lemd4} 
in Section \ref{seci}. It is not exactly true for larger $\dimv$.

\begin{example}
\label{exd1}
{\rm
Let
\[
\kappa(x,y) = \kappa_{(1,1,1)}(x,y) = \frac{xy-x-y}{1-x-y}.
\]
Let $\dimv=4, \dimpol =6$ and
$\a_5(x) = \kappa(x_1, x_2)$ and $\a_6(x) = \kappa(x_3, x_4)$.
Then $bT(\a_5) = \{ \om \in \T^4 : |\om_1 + \om_2 | = 1 \}$, and
$bT(\a_6) = \{ \om : | \om_3 + \om_4 | = 1 \}$.
Their intersection is 2-dimensional, and this is true after perturbation as well,
but $\a_5$ and $\a_6$ depend on different variables so cannot form an edge family.
}\end{example}

 However,
Lemma \ref{lemd4} remains true if two of the functions share a $(\dimv -1)$-dimensional piece of boundary.
\begin{defin}
\label{DEFd2}
Let $\om \in b T(\a)$. We say $\om$ is in the $(\dimv -1)$-dimensional boundary, and write
$\om \in b_{\dimv -1} T(\a)$, if there exists $j$ so that $\om \in b \overline{T(\a_j)}$ and there is a neighborhood
$U$ of $\om$ in
$\T^\dimv$ so that
$U \cap 
 b \overline{T(\a_j)} =U \cap b T(\a)$.
If this occurs, we shall say that $b T(\a_j)$ shares $ b_{\dimv -1} T(\a)$ at $\om$.
\end{defin}
When $\dimv=2$ the $(\dimv -1)$-dimensional boundary is all of $b{T(\a)}$ under Assumption \ref{ass2}, and saying that
$b T(\a_j)$ shares $ b_{1} T(\a)$ at $\om$ just means $\om \in bT(\a) \cap b\overline{T(\a_j)}$.
In Example \ref{exd1}, $bT(\a_5)$ shares $b_3T(\a)$ on the set $\{ \om \in \T^4 : |\om_1 + \om_2| = 1, |\o_3 + \o_4| < 1 \}$.

\begin{lemma}
\label{lemd5}
Let $\dimv \geq 2$.
Let $\ell \ge \dimv+2$.
Suppose there exists $\vare > 0$ so that, for all $a \in V \cap B(0,\vare)$, if we write as in \eqref{eq20}
\[
\tilde \a_j(y_1, \dots, y_\dimv) \ = \vp_{a_j} \circ \a_j(\vp_{-a_1}(y_1), \dots, \vp_{-a_\dimv}(y_\dimv) ) ,
\]
then $ b {T( \tilde \al_{\dimv+1}) }, \dots , b {T(\tilde \al_{\ell}) }$ all share
$b_{\dimv -1}T(\a)$ at some point that lies in $ T(\tilde \al_{\ell+1} )\cap \dots \cap T(\tilde \al_N)$.
 Then $\al_{\dimv+1}, \dots , \al_{\ell}$ are in $\M$, and they all depend on the same pair of variables.
Moreover, they form an edge family.
\end{lemma}
We prove Lemma \ref{lemd5} in Section \ref{secj}.
Both Lemmas \ref{lemd3} and \ref{lemd4} are special cases of Lemma \ref{lemd5}. We prove them first to make the proof easier to follow.

{\bf Step 5:}
We shall now divide the functions $\{ \a_j : \dimv+1 \leq j \leq N \}$ into two subsets, the bounding functions and the non-bounding functions.

\begin{defin}
\label{DEFd4}
We shall say that $\a_j$ is a bounding function at some point $\o \in \T^\dimv$  if:

(i) 
 The set $bT( \a_j)$ shares $b_{\dimv -1}T( \alpha)$ at some point $\om$;

 (ii) moreover   for all small perturbations of $\tilde V$ there is a perturbation $\tilde \om$
  so that the set $bT(\tilde  \a_j)$ shares $b_{\dimv -1}T( \tilde \alpha)$ at $\tilde \om$.
 
If no such $\o$ exists,  we shall say that $\a_j$ is non-bounding.
   \end{defin}
   
  \begin{defin} 
\label{defd5}
If $\a_j$ is bounding at $\o$, and in addition amongst all the bounding functions $ \a_{j'}$ at $\o$
the set $T(\a_j)$ is minimal, we shall say that $\a_j$ is a minimal bounding function.
\end{defin}

\begin{remark}{\rm
Minimal bounding functions are not necessarily unique, even after perturbations. 
Example \ref{exx5} gives an example of an edge pair $(\a,\b)$ with the property that $\hatt(\a)$ has the anchor point as its initial point, and $\hatt(\b)$ has it as the final point, so neither is contained in the other.}
\end{remark}

We shall let $\BB$ denote the bounding $\a_j$'s. 
By Lemmas \ref{lemd3}, \ref{lemd4} and \ref{lemd5},  each function $\a_j $ in $\BB$ is in $\M$.
\begin{lemma}
\label{lemd7}
Let $\a_\ell $ and $\a_{\ell'}$ be in $\BB$.
Then either $\a_\ell$ and $\a_{\ell'}$ depend on the same pair of variables, or on disjoint pairs of variables.
\end{lemma}
Lemma \ref{lemd7} is proved in Section \ref{secka}.

The functions in $\BB$ that depend on the same pair of variables are  an edge family.
\begin{lemma}
\label{lemd6} 
Suppose for some pair $i,j \in \{ 1, \dots, d \}$ the set ${\EE} = \{ \a_k \in \BB : \a_k \ {\rm is\ a\ function\ of\ } x_i,x_j \}$ is non-empty.
Then
the functions in $\EE$ form an edge family.
\end{lemma}
We prove Lemma \ref{lemd6} in Section \ref{seck}.

There are $m_2$ pairs to which Lemma \ref{lemd6} applies. Let $\EE_j $, for $ 1 \leq j \leq m_2$, denote the corresponding edge families.
Let 
\[
\EE'_j = \{ \a_\ell: \a_\ell \ {\rm is\ in\ the\ same\ edge\ family\ as\ } \EE_j\} .
\]
So $\EE_j'$ may include functions that are not bounding (but perhaps become bounding after a M\"obius permutation).
 We have that $\BB$ is the disjoint union $\cup_{j=1}^{m_2} {\EE_j}$.
 Let $\BB' = \cup_{j=1}^{m_2} {\EE'_j}$.

Once Lemmas \ref{lemd7} and \ref{lemd6} are proved, it turns out that $T(\a)$ has a very simple form. If we knew this at the outset, various technical steps in the proof could be omitted.
\begin{cor}
\label{cord3}
There is some permutation of the coordinates $\{x_1, \dots , x_\dimv \}$, some $k \le \frac{\dimv}{2}$,
and some triangular triples $(1, A_j, B_j)$ so that
\be
\label{eqd6}
\notag
\overline{T(\a)} = \{ \o \in \T^\dimv: |A_j \o_{2j-1} + B_j \o_{2j}| \le 1, \ 1\le j \le k \}.
\ee
\end{cor}

{\bf Step 6.}

Permute the variables so that
\[
\BB \= \bigcup_{j=1}^{m_2} \EE_j \= \{ \a_\ell : \dimv+1 \le \ell \le m_3 \} .
\]
Define
\be
\label{Eq27}
\notag
G \= \bigcap_{\a_\ell \in \BB'} \{ x \in \D^\dimv : |\a_\ell(x) | < 1 \} .
\ee

\begin{lemma}
\label{lemd9}
Let $j' > m_3$. Then $\a_{j'} \in O(G, \D)$.
\end{lemma}
We prove Lemma \ref{lemd9} in Section \ref{secl}.

For each $j$, let 
\[
\LL_j \= \{ (x_1, x_2, \a_\ell(x_1, x_2) ) : \a_\ell \in \EE'_j \} .
\]
This is a \roms\ of height $\# \EE'_j$.
We have proved that if $V$ is irreducible, then after an automorphism of $\D^N$ it is a graph over
$\D^{m_1} \times \LL_1 \times \dots \times \LL_{m_2}$.

{\bf Step 7:}
The last step is to prove that $V$ is in fact irreducible; we do this in Section \ref{secm}.

\begin{lemma}
\label{lemd10}
If $V$ is an algebraic set that is a  \car set for $\D^N$, then $V$ is irreducible as an analytic subvariety 
of $\D^N$.
\end{lemma}


\section{Examples}
\label{secx}

\subsection{Examples of $V$ in $\D^4$}

Let us consider some examples of two-dimensional algebraic sets in $\D^4$.
They will all be  of the form
\begin{eqnarray}
\notag
V = \{ (x,y,z,w) \in \D^4 : Ax + By - z &=& \bar{B} xz + \bar{A} yz - xy, \\
Cx + Dy - w &=& \bar{D} xw + \bar{C} yw - xy \} .
\label{Eq28}
\end{eqnarray}
Writing
\beq
\a(x,y) &=& \frac{Ax + By + xy}{1 + \bar B x + \bar A y} \\
\b(x,y)  &=& \frac{Cx + Dy + xy}{1 + \bar D x + \bar C y}
\eeq
then $V$ is the intersection of the graph of $(\a,\b)$ with the 4-disk.

Recall from \eqref{Eq6} that  $\hatt(\a) = \{ \o \in \T: |A + B \o | < 1 \}$.
We shall let $j$  denote $e^{i\pi/3}$.

\begin{lemma}
\label{lemx1}
The triple $(A,B,1)$ is triangular if and only if $\overline{\hatt(\a)}$ is a proper subarc of $\T$.
\end{lemma}
\bp
The triple $(A,B,1)$ is triangular if and only if there is some $\theta \in (0,\pi)$ such that the points $|A|$ and $|B|e^{i\theta}$ are distance $1$ apart. This means there are unimodular constants $\om_1$ and $\eta_1$ so that
$A + B \om_1 = \eta_1$. Replacing $\theta$ with $-\theta$, we also get $A + B \om_2 = \eta_2$,
where $\om_1 \neq \om_2$.
The set $\hatt(\a)$ is then one of the circular arcs joining $\om_1$ and $\om_2$.
\ep
We shall assume that $\a$ and $\b$ are in $\M$, so by Lemma \ref{lemx1} we have $\hatt(\a)$ and $\hatt(\b)$ are both proper subarcs of $\T$.
Recall that $D_\a = \{ (x,y) \in \D^2: |\a(x,y)|  < 1 \}$.
\begin{lemma}
\label{lemx2}
If $D_\beta \subseteq D_\a$, then $\hatt(\b) \subseteq \hatt(\a)$.
\end{lemma}
\bp
Let $\om \in \T$. Then
\be
\label{Eq29}
\a(x,\om x) = \om x \frac{x + (A \bar \om + B) }{1 + (\bar A \om + \bar B )x} .
\ee
If $\om \in \hatt(\a)$, then \eqref{Eq29} is a degree 2 Blaschke product.
If $\om \in b \hatt(\a)$, then \eqref{Eq29} is the unimodular constant $(A+B \om)$ times $x$.
If $\om \in \T \setminus \overline{\hatt (\a)}$, then \eqref{Eq29} has a pole in $\D$.

If there were some point $\om \in \hatt (\b) \setminus \overline{\hatt (\a)}$, then the entire disk
$\{ (\l, \om \l) : \l \in \D \}$ would lie in $D_\beta$, but there would be some points not in $D_\a$.
\ep

\begin{example}
\label{exx01}
{\rm
$A=B=-1$ and $C = 1-j, D=j$, where $j$ denotes $e^{i\pi/3}$. Then the pair $\a, \b$ form an edge family, since
\[
A + B j^2 = -j = C + D j^2 ,
\]
and hence $V$ is a \roms.

We have $\hatt(\a) = \{ \eith: \frac{2\pi}{3} < \theta < \frac{4\pi}{3}\}$
and $\hatt(\b) =  \{ \eith: 0 < \theta < \frac{2\pi}{3}\}$. So $\{x_1, x_2 \}$ do not give a good parametrization.
If we change variables and write $V = \{ ( \g_1(x_3, x_4), \g_2 (x_3,x_4), x_3, x_4) \}$, then
$\g_1, \g_2$ are in $\M$, 
and calculating we get $\g_1 = [-1; \frac12 - \frac{i}{2\sqrt{3}},- \frac{i}{\sqrt{3}}]$
and $\g_2 = [ -1; \frac12 + \frac{i}{2\sqrt{3}}, \frac{i}{\sqrt{3}}]$.
This gives $\hatt(\g_1) = \{ \eith: \frac{2\pi}{3} < \theta < 2\pi \}$ and
$\hatt(\g_2) = \{ \eith: 0 < \theta < \frac{4\pi}{3}\}$. These overlap, though not at the anchor point $1$.  We have $\{ x_3, x_4 \}$ give a good parametrization.

By Lemma \ref{lemx2}, $V$ is not a graph over a \roms\ of height $1$ over $\{ x_3, x_4 \}$.
Since $\a$ has the property that $\a(x_1,x_2) = x_3$ is the same as $\a(x_1, x_3) = x_2$ and $\a(x_2, x_3) = x_1$,
it follows that $V$ is not a graph over any \roms\ of height $1$, and that the good parametrizations are the ones where one of the coordinates is $x_4$.
}
\end{example}

\begin{example} 
\label{exx2}
{\rm
Let $A=\frac{1}{2} - \frac{\sqrt{2}-1}{2}i, B=\frac{1}{2} + \frac{\sqrt{2}-1}{2}i $ and $C = B, D= A$. Then the pair $\a, \b$ form an edge family, since
\[
A + B  = 1 = C + D ,
\]
and hence $V$ is a \roms.
We have $\hatt(\a) = \{ \eith: 0 < \theta < \frac{3\pi}{2}\}$
and $\hatt(\b) =  \{ \eith: \pi/2 < \theta < 2\pi\}$. 
So $\hatt(\a) \cap \hatt(\b) = \{ \eith: \frac{\pi}{2} < \theta < \frac{3\pi}{2}\}$, and hence $\{ x_1, x_2 \}$ form a good parametrization, even though the anchor point $1$ is not in the closure of $\hatt(\a) \cap \hatt(\b)$.
}\end{example}

\begin{example} 
\label{exx2.1}
{\rm
Let $A= -1 - rj^2, B=-1+r$ and $C = -1 , D= -1 $, for $0 < r < 1$. Then the pair $\a, \b$ form an edge family, since
\[
A + B j^2 = -j = C + D j^2 ,
\]
and hence $V$ is a \roms\ of height $2$.
We have $\hatt(\b) = (j^2, j^4) \subsetneq \hatt(\a) = (j^2,  \frac{r+j^4}{1+rj^4}) $.
The anchor point $j^2$ is in  the closure of $\hatt(\a) \cap \hatt(\b)$.

Switch to $x_1, x_3$ coordinates, $V = \{ x_1, \g_2(x_1, x_3), x_3, \g_4(x_1, x_3) \}$.
Then $\g_4 = [1; rj^2, 1-r]$.
As $|rj^2| + |1-r| = 1$, we have $\g_4 \in O(\D^2, \D)$,
so $V$ is also a graph over a \roms\ of height $1$.
}\end{example}

\begin{example}
\label{exx1}
{\rm
Suppose that  $\hatt(\a) $ and $\hatt(\b)$ are  overlapping arcs but 
with neither one contained in the other,  and with the closure of their union not being all of the circle.
Then $\a, \b$ do not form an edge family, since $b \hatt (\a) \cap b \hatt (\b)$ is empty, so $V$ is not a \roms\
by 
Lemma \ref{lemd6}.
By Lemma \ref{lemx2}, it is also not a graph over a \roms\ of height $1$ in the variables $\{ x,y,z \}$
or $\{ x, y, w \}$. 
}
\end{example}

The converse to Lemma \ref{lemx2} is not true.
\begin{example} 
\label{exx3}
{\rm
We shall give an example where $\overline{\hatt(\b)} \subseteq \hatt(\a)$,
but when we change the base point $a$ as described in Subsection \ref{subseccd}, then $\hatt(\tilde \b) \not\subseteq \hatt(\tilde \a)$.
If $D_\beta \subseteq D_\a$, then after all automorphisms of $\D^4$ we would still have $D_{\tilde \b} \subseteq D_{\tilde \a}$, and so by Lemma \ref{lemx2} we would have $\hatt(\tilde \b) \subseteq \hatt(\tilde \a)$.

Let $A = B = 1$, and let $C = 2.5, D= 3$.
Then, approximating to two decimal places,
 \beq
\hatt(\a)  &=& \{ \eith: \frac{2\pi}{3} < \theta < \frac{4\pi}{3} \} =  \{ \eith: 2.09 < \theta < 4.19 \}
\\
\hatt(\b)  & = & \{ \eith: 2.82 < \theta < 3.46 \}
\eeq
Let 
\[
a= (\frac{-1+ 3i}{4}, \frac{i}{2}, \alpha(\frac{-1+ 3i}{4}, \frac{i}{2}), \beta(\frac{-1+ 3i}{4}, \frac{i}{2})) .
\]
By computation, $|a_3|$ and $|a_4|$ are less than $1$, so $a \in V$.
If we apply the M\"obius automorphism $\Phi_a$ we get new functions $\tilde \a$ and $\tilde \b$ in $\M$,
and by \ref{eq21} we get that their partial derivatives at $0$ become (to two decimals)
\[
\tilde A = 0.07 -0.72 i, \tilde B = 0.34-0.64i, \tilde C = -0.99+0.30i, \tilde D = -0.32-1.76i.\]
We then find that
\[
\hatt(\tilde \a) =  \{ \eith: 1.21 < \theta < 4.28 \},
\quad 
\hatt(\tilde \b) = \{ \eith: 0.96 < \theta < 1.95 \}.
\]
So $V$ is not  a \car set.
}
\end{example}

\begin{example} 
\label{exx4}
{\rm
Now we give an example of a set of the form \eqref{Eq28} that is a \car set that is a graph over a  {\roms}\
 in $\D^3$.
Choose positive constants so that $ \frac{1}{2} < C = D < A = B $.

Then $D_\a \subsetneq D_\b$.
Indeed,
$(x,y) \in D_\a$ if and only if it is in $\D^2$ and $|\a(x,y)|^2 < 1$, which is the condition
\be
\label{Eq30}
\frac{2}{A} \Re [ (x+y) (\bar x \bar y - 1) ] < \frac{1 - |xy|^2}{A^2} .
\ee
If \eqref{Eq30} holds for $A$, then it also holds for $C < A$.
Therefore $D_\a \subset D_\b$, and the inclusion must be proper since $T(\a) \subsetneq T(\b)$.
Therefore $V$ is the graph of $\beta$ over the \roms\ $\{ (x,y,z) \in \D^3 : z = \a(x,y) \}$.
If we write $V$ as a graph over $\{x,z\}$, it is still a graph over a \roms, but the fourth coordinate is given by a function that is not in $\M$.
}
\end{example}

\begin{example} 
\label{exx5}
{\rm
Let $C=D=-1$. Let $r \in \R \setminus \{ 0 \}$, let $A = \frac{j^5}{r-j}$ and $B= rA$.
We have $C + D j^2 = j, \ C+ D j^4 = j^5$, and 
$\hatt(\b)$ is the arc running counterclockwise from $j^2$ to $j^4$.
$A + B j^2 = j$ and $ A + B j^4 = \frac{j^5 - r}{r-j}$ which is unimodular, but only equals $j^5$ if $r=1$.
So $\a$ and $\beta$ form an edge pair, 
and $\hatt(\a)$ and $\hatt(\b)$  have the same end-points; the arcs coincide if $r > 0$, and are complementary if $r < 0$.
If we perturb to some point $(p,q, \a(p,q), \b(p,q))$ the arcs will still have one endpoint in common (since the functions form an edge pair), but depending on the perturbation we can have $\hatt(\a)$ and $\hatt(\b)$ go in opposite directions but overlap
(e.g. $r=-1, p = -0.1, q = 0$
which gives $\hatt(\tilde \a) = \{ \eith : -1.93 < \theta < 2.28 \}$
and $\hatt(\tilde \b) = \{ \eith: \theta \notin [-1.93,1.93] \}$),
 or are disjoint 
(e.g. $r=-1, p = 0.1, q=0$, which gives $\hatt(\tilde \a) = \{ \eith : -2.28 < \theta < 1.93 \}$
and $\hatt(\tilde \b) = \{ \eith: \theta \notin [-2.28,2.28] \}$).
}
\end{example}

\subsection{An example of $V$ in $\D^5$}
\begin{example}
\label{exx8}
{\rm
Let $A_3 = -1 = B_3$, let $A_4 = -\frac{1}{2} - \frac{i}{2\sqrt{3}} = B_5$, and $B_4 = \bar{A_4} = A_5$.
Let $\a_j = [ 1; A_j, B_j]$ for $3 \le j \le 5$.
We have
\beq
A_3  + B_3 j^2 &= j^4 &= A_4 + B_4 j^2 \\
A_3 + B_3 j^4 &= j^2 &= A_5 + B_5 j^4 \\
A_4 + B_4 1 &= -1 &= A_5 + B_5 1 .
\eeq
It is easy to check that $\hatt(\a_3) = (j^2, j^4),\ \hatt(\a_4) = ( j^2, 1)$ and $\hatt(\a_5) = (1,j^4)$,
where the arcs are counterclockwise.

Let $\LL = \{ x \in \D^5: x_j = \a_j(x_1, x_2), 3 \le j \le 5 \}$.
This is an example of the situation in Theorem \ref{thmb3} that is not a \wrs. It follows from the proof of that theorem that
$\LL$ is a graph over the \roms
$\{ x \in \D^3 : x_3 = \a_3(x_1, x_2) \}$, and so is a \car set.
}\end{example}

\subsection{Examples of functions}

\begin{example}
\label{exx6}
{\rm
Let  $2 \le k \le \dimv$, and define
$f(x) = \frac{1}{k} (x_1 + \dots + x_k)$.
Then $T(f) = \{ \om \in \T^\dimv : \om_1 = \dots = \om_k \}$.
The boundary of $\overline{T(f)}$ is empty, and the boundary of $T(f)$ has dimension
$\dimv-k+1$.
}
\end{example}

\begin{example}
\label{exx7}
{\rm
In \cite{kn07b}, Knese gives a description of all holomorphic functions $f : \D^\dimv \to \D$ for which equality holds everywhere in the Schwarz inequality, i.e.  for which
\be
\label{Eq30.1}
\sum_{i=1}^\dimv (1-|x_i|^2) \left| \frac{\partial f}{\partial x_i} (x) \right|
\= 1 - |f(x)|^2 
\ee
holds at every point $x \in \D^\dimv$. They are all rational inner functions, of degree $1$ in each variable, and 
with further restrictions which he describes.
Let us choose such a function, say
\be
\label{Eq30.2}
\notag
\a(x_1,x_2,x_3) = \frac{(x_1+x_2+x_3) - 2(x_1x_2+x_1x_3 + x_2x_3) + 3x_1x_2x_3}{3-2(x_1+x_2+x_3)+(x_1x_2+x_1x_3+x_2x_3)} .
\ee
This corresponds to choosing the symmetric unitary $U$ from the realization formula in \cite[Thm. 1.2]{kn07b} to be 
\[
\begin{pmatrix}
0 & \frac{1}{\sqrt{3}} & \frac{1}{\sqrt{3}} & \frac{1}{\sqrt{3}} \\[6pt]
\frac{1}{\sqrt{3}} & \frac{2}{3} & -\frac{1}{3} & -\frac{1}{3} \\[6pt]
\frac{1}{\sqrt{3}} & -\frac{1}{3} & \frac{2}{3} & -\frac{1}{3} \\[6pt]
\frac{1}{\sqrt{3}} & -\frac{1}{3} & -\frac{1}{3} & \frac{2}{3}
\end{pmatrix} .
\]
At $0$, we have $\frac{\partial \a}{\partial x_1} = \frac{\partial \a}{\partial x_2} = \frac{\partial \a}{\partial x_3} = \frac{1}{3}$,
so $T(\a) = \T^3 \setminus \{ \o_1 = \o_2 = \o_3 \}$, and $b\overline{T(\a)}$ is empty.
The set $V = \{ x_1,x_2,x_3, \a(x_1,x_2,x_3) : (x_1,x_2,x_3) \in \D^3 \}$ is a retract of $\D^4$, and hence a \car set.

If we choose a point $a \in V$ and perturb $V$ to $\tilde V = \Phi_a(V)$, then
by \eqref{eq21}
\[
\sum_{i=1}^3
\left|\frac{\partial \tilde \a}{\partial x_i} \right|
\=
\sum_{i=1}^3 \frac{1-|a_i|^2}{1-|\alpha(a_1,a_2,a_3)|^2} \left| \frac{\partial  \a}{\partial x_i} \right|
\]
which is $1$ by \eqref{Eq30.1}.
So for all $a$ we have $b T(\tilde \a)$ is one-dimensional, and hence Lemmas \ref{lemd3} and \ref{lemh6} are not true if we just assume that
$b T(\tilde \a_{\dimv +1})  \cap bT(\tilde \a) \cap_{\ell \ge \dimv +2} T(\a_\ell)$ is non-empty for all perturbations.
}
\end{example}


\section{Proof of Lemma \ref{lemd1}}
\label{secf}

\begin{lemma}
\label{lemf2}
Let $\X$ be a $d$-dimensional complex vector space, and let $\{ v_1, \dots, v_N \}$ be vectors that
span $\X$ and such that no two are unimodular multiples of each other.
Let $L_i(u) := \langle u, v_i \rangle$.  
Then for each $ 1 \leq m \leq d$, there exists $\sigma \in {\mathfrak S}_N$ and $u \in \X$
so that 
\begin{equation}
\label{Eq31}
\notag
| L_{\sigma(1)} (u) | = |L_{\sigma(2)}(u) |\dots =  |L_{\sigma(m)}(u) | > \max_{j > m}  |L_{\sigma(j)}(u) |,
\end{equation}
and so that $\{ L_{\sigma(1)}, \dots, L_{\sigma(m)} \}$ are linearly independent.
\end{lemma}

\begin{proof}
We assume $N > d$, since otherwise it is obvious.
When $m=1$, the claim follows from the fact that the Hausdorff dimension
of 
\begin{equation}
\label{Eq32}
\notag
E_2 = 
\{ u \in \X : {\rm \ for\ some\ }i\neq j, \ |L_i(u) | = |L_j(u)| \} 
\end{equation}
is $\leq 2d -1 $.
Any point  $u$ not in this set then works.

For $2 \le m \le d$ we will prove by induction on $m$.
Let $\I = \{ 1, \dots, N\}$. 
We shall use the following inductive hypothesis:

(IH) {\em  Suppose there is $u_1$ in $\X$ and $i_1 \in \I$ so that $|L_{i_1}(u_1)| > \max_{j \neq i_1} |L_{j}(u_1)|$.
Then there exists distinct $i_2, \dots, i_m$ in $\I\setminus \{ i_1\}$ and $u_m \in \X$ so that
\be
\label{eq32.1}
|L_{i_1}(u_m)|  = |L_{i_2}(u_m)|  = \dots = |L_{i_m}(u_m)|  > |L_{j}(u_m)| \quad \forall j \notin \{ i_1, \dots, i_m \} ,
\ee
and so that $\{ L_{i_1}, \dots, L_{i_m} \}$ is linearly independent.
Moreover, there is an uncountable set of $\omega$ for which $u_m$ can be chosen satisfying \eqref{eq32.1} and
with $L_{i_2}(u_m) = \omega L_{i_1}( u_m)$.
}

For $m=2$, choose $u_0$  in $\X \setminus E_2$ that is  maximized by a different
functional (we can do this because $m \leq d$). 
Let
\begin{equation}
\label{Eq33}
\notag
E_3 = 
\{ u \in \C^d : {\rm \ for\ some\ distinct\ }j_1, j_2, j_3,  \ |L_{j_1}(u) | = |L_{j_2}(u)|  = |L_{j_3}(u)|\} .
\end{equation}

The Hausdorff dimension of $\{ u \in \C^d: |L_{j_1}(u) | = |L_{j_2}(u)|  = |L_{j_3}(u)|\} $ is at most  $2d-2$.  So the dimension of $E_3$ 
is at most $2d-2$, and does not disconnect $\C^d$.

Draw a path connecting $u_1$ to $u_0$ that does not intersect $E_3 $.
At the first point where $|L_{i_1}(u)|$ ceases to be strictly larger than all the other components,
we get $u_2$ satisfying \eqref{eq32.1} and $L_{i_2}(u_2) = \o L_{i_1}(u_2)$ for some particular $\o$.
Moreover, $L_{i_1}$ and $L_{i_2}$ are linearly independent, since otherwise they would be unimodular multiples of each other, which is excluded.

Now, choose $z \in \X$ so that $L_{i_2} (z)  \neq \o L_{i_1} (z) $. 
For $\theta$ small, let 
\be
\label{eq34.1}
\notag
\vare = \frac{(e^{i\theta} - 1)\omega L_{i_1} (u_2) }{L_{i_2} (z)  - \omega e^{i\theta} L_{i_1} (z) }.
\ee
Then we have 
\[
L_{i_2} (u_2 + \vare z) = e^{i\theta} \o L_{i_1} (u_2 + \vare z) .
\]
So replacing $u_2$ by $u_2 + \vare z$ we get \eqref{eq32.1} for $\o e^{i\theta}$ for all $\theta$ sufficiently close to $0$.

Now assume (IH) has been proved  for $m-1 \ge 2$.
We want to prove it for $3 \le m \le d $.

We have $|L_{i_1}(u_1)| > \max_{j \neq i_1} |L_{j}(u_1)|$.
First we apply (IH) for 2 to get (after a normalization)
\be
\label{eq32.2}
\notag
1 = L_{i_1} (u_2) = \bar \o L_{i_2} (u_2)  >  \max_{j \neq i_1, i_2} |L_{j}(u_2)| .
\ee
Now let $Y_\o =  {\rm ker}( L_{i_2} - \om L_{i_1})$.
We want to apply the inductive hypothesis to $Y_\o$.

It is possible that some of the $L_j$'s are unimodular multiples of each other when restricted to $Y_\o$.
If so, we discard all but one element in each set that are.
Then applying (IH) we get $u_{m} \in Y_\o$ and distinct $i_3, \dots, i_m \in \I \setminus \{i_1, i_2 \}$ so that
\be
\label{eq32.3}
|L_{i_1}(u_{m})| = |L_{i_3}(u_{m})| = \dots = |L_{i_m}(u_{m})| > |L_j (u_{m})| 
\ee
whenever $L_j |_{Y_\o}$ is not a unimodular multiple of some $L_{i_k} |_{Y_\o}$ for $1 \le k \le m$.
Moreover, the set of possible ratios of $L_{i_3}(u_{m})/L_{i_1}(u_{m})$ is uncountable.

Now let us analyze possible repeats in \eqref{eq32.3}. We claim they can only happen for finitely many choices of $\o$,
which we can then discard.

For simplicity of notation, let us assume that $L_3 \in
\{ L_{i_k} : k \in \{1, 3, \dots m \}\}$ and $L_4$ is not, and $L_4 |_{Y_\o} = \tau L_3 |_{Y_\o}$ for some unimodular $\tau$.
Note first that $L_3 \neq L_{i_1}$, since there is no $L_4$ with the same modulus at $u_2 \in Y_\o$.

Let $V$ be the span of $v_{i_1} $ and $v_{i_2}$, and let $y_1, y_2$ in $V$ be biorthogonal to $v_{i_1} $ and $v_{i_2}$
(so $L_{i_k}(y_j) = \delta_{jk}$ for $j,k \in \{ 1,2\}$).
Then \[
\X = \{ a y_1 + b y_2 + \xi : a,b \in \C, \xi \in V^\perp \}\]
 and
 \[
 Y_\o = \{ a(y_1 + \o y_2) + \xi : a \in  \C, \xi \in V^\perp \} .
 \]
We are assuming
\be
\label{eq32.4}
\notag
L_4 ( a(y_1 + \o y_2) + \xi) = \tau L_3 ( a(y_1 + \o y_2) + \xi),
\ee
so $L_4 \xi = \tau L_3 \xi$ and 
\be
\label{eq32.5}
L_4 (y_1 + \o y_2)  = \tau L_3 (y_1 + \o y_2)
.
\ee

Case: There is some $\xi \in V^\perp$ with $L_4 \xi \neq 0$.

Now we repeat the procedure with a different choice of $\o$, say $\o'$.
We must have the same $\tau = (L_4 \xi)/(L_3 \xi)$, and 
\be
\label{eq32.6}
L_4 (y_1 + \o' y_2)  = \tau L_3 (y_1 + \o' y_2)
.
\ee
Comparing \eqref{eq32.5} and \eqref{eq32.6} we get that $L_4 = \tau L_3$ on $\X$, a contradiction.

Case: $L_3$ and $L_4$ are in the span of $L_{i_1}$ and $L_{i_2}$. 
So $L_3 = c L_{i_1} + dL_{i_2}$. 
We know $|L_3 (u_2) |= |c + \o d | < 1 = L_{i_1}(u_2)$.
Since $u_{m} = a u_2 + \xi$, this gives
$|L_3(u_{m})| = |a| |c + \o d| < |a| = |L_{i_1}(u_{m})|$, so $L_3$ is not one of the functionals in 
\eqref{eq32.3}.

Since there are only finitely many pairs to consider, we only need to discard finitely many $\o$. For all the rest, we get a $u_m$ so that $L_{i_2}(u_m) = \o L_{i_1} (u_m)$ and
\eqref{eq32.1} holds.

Finally, let us check the linear independence. By (IH) we have $L_{i_1}, L_{i_3}, \dots, L_{i_m}$ are linearly independent on $Y_\o$, and $L_{i_1}$ and $L_{i_2}$ are linearly independent. Suppose $\sum_{j=1}^m c_j L_{i_j} = 0$.
Restricting to $Y_\o$, we get $c_j  =0 $ for all $ j \ge 3$, and $c_1 + c_2 \o = 0$. By linear independence of
$L_{i_1}$ and $L_{i_2}$, we get that each $c_j = 0$.
\end{proof}

{\bf Lemma \ref{lemd1}.}
{\em
Let $V$ be a $\dimv$-dimensional algebraic set in $\D^N$ that has $0$ as a regular point.
Assume that Assumption \ref{ass1} holds for $V$ at $0$.
After a permutation automorphism and a M\"obius perturbation, we can assume that $V$ has a good parametrization at $0$.
}

\bp
 Let $0$ be a regular point of $V$. Define linear functionals $L_j : T_0 V \to \CC$ by
$L_j(u) = \langle u, e_j \rangle$, where $e_j$ is the standard basis in $\CC^n$.
(So $v_j$ is the projection of $e_j$ onto $T_0 V$, the tangent space to $V$ at $0$.)

Case (i): For no pair $(i,j)$ do we have 
$L_i $ is a unimodular multiple of $L_j$.
Now we can apply Lemma \ref{lemf2} with $m=d$, and we can permute the coordinates so that $\sigma$ is the identity.
Parametrize $V$ near $0$ with respect to the first $d$ coordinates, so that a neighborhood of $0$
is $\{ (\alpha_j(x_1, \dots, x_d )) : 1 \leq j \leq N \}$ where for $i \leq d$ we let $\alpha_i(x) = x_i$.
Writing $\alpha_j^i$ for $\frac{\partial \alpha_j}{\partial x_i} (0)$, we have
\begin{equation}
\label{Eq35}
\notag
T_0V =  \{ (x_1, \dots, x_d, \sum_{i=1}^d x_i \alpha_{d+1}^i , \dots,\sum_{i=1}^d x_i \alpha_{N}^i  ) : x \in \CC^d \} 
\end{equation}

Let $u$ be as in Lemma \ref{lemf2}, scaled so that $|L_1(u)| = 1$.
Then $u_1, \dots, u_{d}$ are unimodular constants, call them $\om_1, \dots, \om_{d}$,
and the other coordinates have modulus less than one.
This means for $j > d$, 
$u_j = \sum_{i=1}^d \om_i \a_j^i$ has modulus less than one, so
 $ (\om_1 :\dots : \om_{d}) \in \cap_{j=d+1}^N T(\alpha_j) $.

Case (ii): Suppose there is some pair, which without loss of generality we can take to be $(1,2)$, so
that $|L_1| = |L_2|$, and this remains true if we perturb by a M\"obius automorphism.
(If it is not true for all M\"obius perturbations, we apply one such and we are in Case (i) again.)
By Assumption \ref{ass1}, no coordinate is constant on $V$,
which means no $v_j $ is zero.
So we can choose $i_2, \dots, i_d \in \{ 2, \dots , N \}$ so that $\{ v_1, v_{i_2}, \dots, v_{i_d} \}$ form a basis for $T_0 V$, and hence 
we can locally write $V$ as a graph over $(x_1, x_{i_2}, \dots, x_{i_d})$. 
We have near $0$ a parametrization
\be
\label{Eq36}
\notag
\{
(\alpha_1(x_1, x_{i_2}, \dots, x_{i_d}), \dots, \al_N(x_1, x_{i_2}, \dots, x_{i_d})) \ : \
(x_1, x_{i_2}, \dots, x_{i_d}) \in U \}
\ee
where for $j \in \{ 1, i_2, \dots, i_d\}$, we write $\al_j(x_1, x_{i_2}, \dots, x_{i_d}) = x_j$.

If $L_1 = \tau L_2$, then we have
\be
\label{Eq37}
\notag
\begin{pmatrix}
\frac{\partial \al_1}{\partial x_1}(0)\\
\vdots\\
\frac{\partial \al_1}{\partial x_{i_d}} (0)
\end{pmatrix}
\=
\tau \begin{pmatrix}
\frac{\partial \al_2}{\partial x_1}(0)\\
\vdots\\
\frac{\partial \al_2}{\partial x_{i_d}} (0)
\end{pmatrix}
\ee
If $|L_1| = |L_2|$ on $\Phi_a (V)$ for all $a$ in a neighborhood of $0$ in $V$, we get from \eqref{eq21}
\be
\label{Eq38}
\begin{pmatrix}
\frac{ 1- |a_{1}|^2}{1-|a_1|^2} \frac{\partial \al_1}{\partial x_1}(a)\\
\vdots\\
\frac{ 1- |a_{i_d}|^2}{1-|a_1|^2} \frac{\partial \al_1}{\partial x_{i_d}} (a)
\end{pmatrix}
\= \tau(a)
\begin{pmatrix}
\frac{ 1- |a_{1}|^2}{1-|a_2|^2} \frac{\partial \al_2}{\partial x_1}(a)\\
\vdots\\
\frac{ 1- |a_{i_d}|^2}{1-|a_2|^2} \frac{\partial \al_2}{\partial x_{i_d}} (a)
\end{pmatrix}
\ee
Since $\al_1(x) = x_1$, the left-hand side of \eqref{Eq38} is the vector $(1,0, \dots, 0)^t$.
This means $\al_2$ is a function of $x_1$ alone, and we have
\be
\label{Eq39}
\notag
 \left| \frac{d \al_2}{d x_1} \right| \= 
 \frac{1-|\al_2|^2}{1-|x_1|^2} .
 \ee
By Lemma \ref{LEMf3}, we conclude that $\alpha_2(x) = \tau  x_1$, which violates Assumption \ref{ass1}.
\ep

\begin{lemma}
\label{LEMf3}
Let $y$ be a holomorphic function of $x$ in a neighborhood of $0$ in $\C$, and assume $y(0) = 0$.
If
\be
\label{Eq40}
\left| \frac{d y}{d x} \right| \= 
 \frac{1-|y|^2}{1-|x|^2} 
 \ee
then $y = \o x$ for some $\o \in \T$.
\end{lemma}
\bp
From \eqref{Eq40}, we have that 
\be
\label{Eq41}
\frac{2 |dy|}{1-|y|^2} \= 
\frac{2 |dx|}{1-|x|^2} .
\ee
If we integrate \eqref{Eq41} as $x$ goes from $0$ to $re^{i\theta}$ along a radius, we get that the hyperbolic distance
$y$ travels is the same as the distance $x$ travels; since $x$ has moved along a geodesic we get that
$|y( x)| \leq |x|$.
Since $\frac{d y}{d x} (0) \neq 0$, we can swap the r\^oles of $x$ and $y$ to get the reverse inequality;
so in a neighborhood of $0$ we have $|y|=|x|$. The conclusion follows.
\ep


\section{Edge Families}
\label{sece}

Throughout this section, we shall 
let $\{\a_1,\dots, \a_m \} \subset \overline{\M}$. We shall assume that Assumption \ref{ass1} holds, which means no $\a_j$ is a unimodular multiple of a coordinate function, and no two $\a_j, \a_k$ are unimodular multiples of each other. 
In Lemmas \ref{leme1} and \ref{leme2}, we assume this applies to $\a$ and $\b$, which we write without subscripts for legibility.

We define
\be
\label{Eq42}
\LL = \{ (x_1, x_2, \a_1(x_1, x_2),\dots,  \a_m(x_1, x_2) ) : (x_1, x_2) \in D_{\a_1} \cap 
\dots \cap D_{\a_m} \}
\ee
\be
\label{Eq42.1}
\LL^\sharp = \{ x \in (\C_\infty)^N : x_j = \a_j(x_1, x_2), \ 3 \le j \le N \}.
\ee
In this section we will prove Theorems \ref{thmb3},  \ref{thme3} and \ref{thme5}.

If some $\a_j \in \ovm\setminus \M$, we shall understand $b \hatt(\a_j)$ to mean the singleton \[
b \hatt(\a_j) \ :=  \
\{ \o \in \T :
|\a_j^1 + \a_j^2 \o| = 1 \}, 
\]
 and we shall understand $\overline{\hatt(\a_j)}$ to mean 
 \[
 \overline{\hatt(\a_j)} \ := \ 
 \{ \o \in \T :
|\a_j^1 + \a_j^2 \o| \le 1 \}.
\]

We start with the following lemma.
\begin{lemma}
\label{leme1}
Suppose $\a, \b$ are in $\overline{\M}$, and are given by
\begin{eqnarray}
\label{Eq43}
\notag
\a(x,y) \= \xi  \frac{x y + A x + B y}{1 + \bar B x + \bar A y} \\
\b(x,y) \= \tau \frac{x y + C x + D y}{1 + \bar D x + \bar C y}.
\label{Eq44}
\end{eqnarray}
Then $\a(x, \b(x,y) )$ is in $\overline{\M}$ if and only if
\be
\label{Eq45}
\eta \ :=\  -\ \frac{ \tau\bar A  D + \bar C}{\bar B \bar C  + \tau \bar A}
\ee
is unimodular and
\be
\label{Eq46}
|A|^2 \eta - (B + \eta) (1 + \bar B \eta) \= 0.
\ee
The function is in $\M$ if in addition $B/D$ is not real.
Moreover, if $\a(x, \b(x,y) )$ is in $\overline{\M}$, then
\be
\label{Eq47}
\a(x, \b(x,y) ) \= - (\xi \tau \eta) \frac{xy+  (C + \bar \tau A \bar D)x - \bar \eta BD y }{1 - \eta \bar B \bar D x +(\bar C + \tau \bar A D) y} .
\ee
\end{lemma} 
\bp
We calculate that
\be
\label{Eq48}
\a(x, \b(x,y) ) \= \xi \frac{P(x) y + Q(x)}{R(x) y + S(x) } ,
\ee
where
\beq
P(x) &\=& \tau x^2 + (\tau B + \tau D + A \bar C) x + \tau B D \\
Q(x) &\=& x \left[ (\tau C + A \bar D) x + (\tau B C + A) \right]\\
R(x) &\=& (\bar B\bar C + \tau \bar A) x + \tau \bar A D + \bar C \\
S(x) &\=& \bar B\bar D x^2 + ( \bar B + \bar D + \tau \bar A C)x + 1 .
\eeq
For \eqref{Eq48} to be in $\overline{\M}$, the functions $P,Q,R,S$ must have a common root,
which cannot be at zero since $S$ and $P$ don't vanish at $0$.
The number $\eta$ given by  \eqref{Eq45} is the root of $R$, and $1/\bar \eta$ is the non-zero root of $Q$,
and these are equal if $\eta$ is unimodular.
We have that $\eta$ is a root of $P$ if and only if it is a root of $\bar A P(x) - (B +x) R(x)$.
Dividing by $\bar C$, this gives exactly \eqref{Eq46}. Hence $\eta$ is also a root of $S$,
 since the roots of $P$ are the reflections in the circle of the roots of $S$.
If we factor $(x-\eta)$ out of all 4 polynomials, we get \eqref{Eq47}.

So we have proved that $\a(x,\beta(x,y))$ is a rational function of degree $(1,1)$ if and only if $|\eta|=1$ and \eqref{Eq46} holds. 
We claim if these happen, then $1, |C+\bar \tau A \bar D|, |BD|$ form the sides of a possibly degenerate triangle, so
the function is in $\overline{\M}$.
Indeed, if $|\eta| = 1$, then 
\eqref{Eq46} says $|A| = |B+\eta|$, so there exists a unimodular $\o$ so that
\be
\label{eq48.1}
A + \om B = - \om \eta.
\ee
From \eqref{Eq45}, we have
\be
\label{eq48.2}
\notag
(\bar \tau A + BC) = - \eta (\bar \tau A \bar D + C) ,
\ee
which with \eqref{eq48.1} gives
\be
\label{eq48.3}
\notag
\bar \tau - \bar \o C = - \eta \bar \tau \bar D .
\ee
So
\be
\label{eq48.4}
\notag
C +  \bar \tau A \bar D = \bar \tau \o (1 - B \bar D).
\ee
The three numbers $1, |1- B \bar D|, |BD|$ form the sides of a possibly degenerate triangle.
It is degenerate when one is the sum of the other two, which occurs if and only if $B \bar D$ is real.
\ep

Suppose $\a \in \overline{\M}$. We shall write $\a^{-1} (x,\underline{z}) =y $ to denote the function in $\overline{\M}$  that gives
the unique solution $y$ to $\a(x,y) = z$, and  $\a^{-1} (\underline{z},y) =x $ to give the function  that yields $x$.
The composition $\a(x, \b^{-1}(x,\underline{z}) )$ is in $\overline{\M}$ if  and only if $\a$ and $\b$ form an edge pair.

\begin{lemma}
\label{leme2}
Suppose  that $\alpha = [\xi; A,B]$ and  $\beta = [\tau; C,D]$ as in \eqref{Eq44}.
Then $\a(x, \b^{-1}(x,\underline{z}) )$ is in $\overline{\M}$ if  and only if
there are $\om_1, \eta_1 \in \T$ so that
\be
\label{Eq49}
A + B \om_1 \= \eta_1 \= C + D \om_1 .
\ee
It is in $\M$ if and only if in addition $B \bar D$ is not real.
If  $\a(x, \b^{-1}(x,\underline{z}) )$ is in $\overline{\M}$ then
\be
\label{Eq55}
\a(x, \b^{-1}(x,\underline{y}) ) \= \left( \xi \bar \tau \bar \o_1 \eta_1 \frac{\bar D}{D} \right) \
\frac{xy + \frac{\tau(A-C)}{\bar D} x + \frac{\o_1 \bar \eta_1 B}{\bar D} y}{1 +  \frac{\bar \o_1  \eta_1 \bar B}{ D}x +
\frac{\bar\tau(\bar A-\bar C)}{D} y} .
\ee
\end{lemma}
\bp
If $\a(x,y) = z$, then 
\[
\a^{-1}(x,\underline{z}) \= \frac{\bar \xi \bar B}{B} \frac{xz - (\xi A/\bar B )x + (1/\bar B) z}{1 + x/B - (\bar \xi \bar A/B) z }
\]
and
\[
\b^{-1}(x,\underline{z}) \=  \frac{\bar \tau \bar D}{D} \frac{xz - (\tau C/\bar D )x + (1/\bar D) z}{1 + x/D - (\bar \tau \bar C/D) z }
\]
So by Lemma \ref{leme1}
 $\a(x, \b^{-1}(x,\underline{z}) )$ is in $\overline{\M}$ if and only if
\be
\label{Eq50}
\eta_0 \ :=\  \ \frac{  \bar A  -  \bar C }{\bar B \bar C  -  \bar D  \bar A }
\ee
is unimodular and
\be
\label{Eq51}
|A|^2 \eta_0 - (B + \eta_0) (1 + \bar B \eta_0) \= 0.
\ee
Equation \eqref{Eq55} follows from \eqref{Eq47}.
Equation \eqref{Eq51} is equivalent to
\[
|A| \= |B + \eta_0| .
\]
This means there exists some $\om_1 \in \T$ so that
\be
\label{Eq52}
A \= - \om_1 (B + \eta_0) .
\ee
Let $\eta_1 = - \om_1 \eta_0$, so we have
\be
\label{Eq53}
\notag
A + \om_1 B = \eta_1 .
\ee
From \eqref{Eq50} and \eqref{Eq52}, we get
\beq
\bar A &\=&  (\bar B \bar C  -  \bar D  \bar A ) \eta_0 + \bar C\\
&=& \bar C( \bar B \eta_0 +1) - \bar D \bar A \eta_0 \\
&=& - \bar C( \om_1  \eta_0 \bar A) - \bar D \bar A \eta_0.
\eeq
Therefore
\be
\label{Eq54}
\notag
\bar C + \bar D \bar \om_1 \= - \bar \om_1 \bar \eta_0 \= \bar \eta_1.
\ee
Taking complex conjugates, we get that \eqref{Eq49} is necessary.

Conversely, suppose \eqref{Eq49} holds.
Then 
\beq
BC - AD &\=& B(C-A) + A(B-D)\\
&=& (B-D)(A + B \om_1 ) \\
&=& \eta_1 (C -A) \bar \om_1
\eeq
Therefore 
$
\eta_0 \= - \bar \om_1 \eta_1 $
is unimodular, and $| B + \eta_0| = | - \bar \om_1 A| = |A|$ as desired.
\ep

Suppose we have a set  of the form \eqref{Eq42}. From the proof of Lemma \ref{leme1}, we see that if we try to rewrite it as a graph over the variables $x_1, x_3$ say, then it would become a set of the form
\[
\{ (x_1, \b_2(x_1, x_3), x_3, \b_4(x_1,x_3), \dots ) : (x_1, x_3) \in U \} ,
\]
for an appropriate set $U$. The function $\b_2$ is always in $\overline{\M}$, and will be in $\M$ if $\a_3$ is,  but for $k \ge 4$,
the function $\b_k(x_1,x_3) = \alpha_k(x_1, \a_2^{-1} (x_1, \underline{x_3}))$ is either in $\overline{\M}$ (which happens if and only if the conditions in Lemma \ref{leme2} are satisfied, i.e.  if and only if $\a_3$ and $\a_k$ form an edge pair), or it is a rational function that is degree $2$ in $x_1$ and degree $1$ in $x_3$. In the latter event, there are $x_1' \neq x_1''$ so that
$\b_k(x_1',x_3) = \b_k(x_1'', x_3)$. This means that if we project $\LL^\sharp$ onto the two-dimensional linear subspace spanned by the $x_3$ and $x_k$ coordinates, it is not one-to-one. Repeating this argument as we change to other pairs of variables,  we conclude:

\begin{prop}
\label{prope2}
Let $\LL^\sharp$ be given by \eqref{Eq42.1}.
Then $\LL^\sharp$ is single-sheeted over each pair of coordinates if and only if, when written as a graph over any pair
$x_j, x_k$, each function is in $\overline{\M}$.
\end{prop}

Lemma \ref{leme2} says that every  $\a_j(x, \a_k^{-1} (x,\underline{z}))$ is in $\overline{\M}$ if and only if
every  $\a_j$ and $\a_k$ form an edge pair, which means we can write $\LL$ as a graph over $(x_1, x_k)$ or $(x_2, x_k)$ for any $k$. If we strengthen this to arbitrary base pairs,
we get that they must either form an edge family, or always form edge pairs with two out of three possible $\o-\eta$ combinations.

Recall the notation $[\xi; A,B]$ is short-hand for a function that has the form
\[
\a(x,y) \= \xi  \frac{x y + A x + B y}{1 + \bar B x + \bar A y}.
\]
We will write $\a^{\leftarrow} = [\xi; B,A]$.

For later use, let us record the following. The proof is from repeated applications of Lemma \ref{leme2}.
\begin{cor}
\label{core1}
Let $\LL = 
 \{ (x_1, x_2, \a_3(x_1, x_2), \dots, \a_N(x_1, x_2)) : 
 (x_1, x_2) \in D_{\a_3} \cap 
\dots \cap D_{\a_N} \}
$ be a \wrs.
Assume $\a_n = [ \tau_n; \eta - \o B_n, B_n]$. 
Define 
\[\sigma_n = \begin{cases}
\bar \omega \eta & n = 1 \\
\eta & n = 2 \\
\bar \o \eta^2 \tau_n & n \ge 3 
\end{cases}
\]

Then for $n \ge 3, n \neq j,k$, we have
\be
\label{Eq64.1}
\notag
\g_n = [ \bar \omega^2 \eta^2  \sigma_n \bar \sigma_j \bar \sigma_k \frac{\bar B_j - \bar B_k}{B_j - B_k}; 
\o^2 \bar \eta^2  \sigma_k \frac{B_n - B_k}{\overline{ B_j} -  \overline{B_k}},  \o^2 \bar \eta^2 \sigma_j\frac{B_j- B_n}{\overline{ B_j} - \overline{ B_k}} ].
\ee
If $j$ or $k$  is $ 1$ or $2$, we interpret $B_1 = 0,$ and $B_2 =\bar \omega \eta $.
\end{cor}

We can now prove Theorem \ref{thmb3}.

{\bf Theorem \ref{thmb3}.}
{\em
Let $\a_3, \dots, \a_N \in \overline{\M}$, and let 
\be
\label{Eq56}
\notag
\LL = \{ (x_1, x_2, \a_3(x_1, x_2), \dots, \a_N(x_1, x_2)) : 
 (x_1, x_2) \in D_{\a_3} \cap 
\dots \cap D_{\a_N} \}.
\ee
Assume \ref{ass1}.
Then the following are equivalent:

(i) For any pair $j,k \in \{ 1, 2, \dots, N \}$, the projection of $\LL^\sharp$ onto the $x_j-x_k$ plane is injective.

(ii) For any pair  $j,k \in \{ 1, 2, \dots, N \}$, there are functions $\g_n$ in $\overline{\M}$ for $n \in \{ 1, \dots , N \} \setminus \{ j,k\}$
so that $\LL$ can be written as 
\[
\LL = \{ (x_1, \dots, x_N) \in \D^N : x_n = \g_n(x_j, x_k), (x_j, x_k) \in D_{\g_1} \cap \dots \cap D_{\g_N} \} .
\]
(iii) For any pair $j,k \in \{ 3, \dots, N \}$, the functions $\a_j$ and $\a_k$ form an edge pair.

Moreover, if these conditions hold, then  either $\LL$ is a \wrs, in which case it is either a retract or a graph over a \roms, or $\LL$ is not a \wrs, in which case it is a graph over a \roms\ of height one.
}

\bp
The equivalence of (i) and (ii) is Proposition \ref{prope2}.

To see that (ii) and (iii) are equivalent, let
 us change variables from writing $\LL$ as functions over $(x_1, x_2)$ to writing it as functions over $(x_n, x_k)$.
Let $\a_j = [ \tau_j; A_j, B_j]$.
When   can we write $\LL$ as a surface over any other pair?
By Lemma \ref{leme2}, we can do this with $n=1$ and any $k$ if and only if for every $j,k$ we have $\a_j, \a_k$ form an edge pair. 
So  condition (iii) is necessary, i.e.   that  for all $j,k$ in $\{ 3,4, \dots, N \}$  there are unimodular constants $\om_{jk}, \eta_{jk}$
satisfying
\be
\label{Eq57}
A_j + B_j \om_{jk} = \eta_{jk} = A_k + B_k \om_{jk} .
\ee
To prove (iii) is sufficient, we must investigate when we can write $\LL$ as a surface over a pair of coordinates containing neither $x_1$ nor $x_2$, which 
for notational simplicity we shall assume is the pair
$x_3, x_4$.

First, we rewrite $\LL$ as functions of $x_1,x_3$, so $\LL$ is
\[
\{ (x_1, \beta_2(x_1, x_3), x_3, \beta_4(x_1, x_3), \beta_5(x_1, x_3) , \dots) : (x_1,x_3) \in U \}
\]
for some set $U$.
We can do this because each pair $\a_3, \a_k$ is an edge pair.
Moreover, we have from \eqref{Eq55}
\begin{eqnarray}
\label{Eq58}
\notag
\beta_2 &=& [\bar \tau_3 \frac{\bar B_3}{B_3} ; \ - \tau_3 \frac{A_3}{\bar B_3}, \ \frac{1}{\bar B_3}]  =: [\xi_2; C_2, D_2]\\
\beta_k &=& [ \bar \tau_3 \tau_k \bar \om_{3k} \eta_{3k} \frac{\bar B_3}{B_3} ;\  \tau_3 \frac{A_k - A_3}{\bar B_3} ,\
\om_{3k} \bar \eta_{3k} \frac{B_k}{\bar B_3} ] =: [ \xi_k; C_k, D_k] \quad k \ge 4 .
\label{Eq59}
\notag
\end{eqnarray}

Now, we want to write $\LL = \{ \g_1 (x_3,x_4), \g_2(x_3, x_4), x_3, x_4, \g_5(x_3, x_4) , \dots \}$.
We will have $\gamma_1(x_3, x_4) = \beta_4^{-1} ( \underline{x_4}, x_3)$, and for other $n$ we have
$\gamma_n(x_3, x_4) = \beta_n(  \beta_4^{-1} (\underline{x_4}, x_3), x_3)$. By Lemma \ref{leme2} (with the order of the variables reversed)
$\g_n$ is in $\overline{\M}$ if and only if  the pair $\beta_4^\leftarrow, \beta_n^{\leftarrow}$ is an edge pair, which is the same as saying $\beta_4, \beta_n$ form an edge pair.

It is easy to see that $\beta_n, \beta_2$ form an edge pair, because, writing $\o,\eta$ for $\o_{2n}, \eta_{2n}$, 
\begin{eqnarray}
\notag
C_2 + D_2 \tau_3 \eta &=& -\tau_3 \frac{A_3}{\bar B_3} + \frac{1}{\bar B_3} \tau_3 \eta \\
\notag
&=&  \frac{\tau_3}{\bar B_3} \left[-  A_3 + A_3 + B_3 \o \right] = \tau_3 \o \frac{B_3}{\bar B_3}\\
\notag
C_n + D_n \tau_3 \eta  &=&  \tau_3 \frac{A_n - A_3}{\bar B_3} + \om \bar \eta \frac{B_n}{\bar B_3}\tau_3 \eta \\
&=& \frac{\tau_3}{\bar B_3} \left[ A_n + \o B_n -A_3\right] = \tau_3 \o \frac{B_3}{\bar B_3}.
\label{eqe107}
\end{eqnarray}

To see that $\b_4$ and $\b_n$ form an edge pair for $n \neq 2$, note that instead of rewriting in terms of $x_1$ and $x_3$, we could equally have written $\LL$ as
\[
\{ (x_1, \gamma_2(x_1, x_4),  \gamma_3(x_1, x_4), x_4, \gamma_5(x_1, x_4) , \dots)  \}
\]
This means 
\[ \b_n(x_1 , \b_4^{-1} (x_1, \underline{x_4})) =
\g_n(x_1, x_4) .
\]
Since $\g_n \in \overline{\M}$, we conclude from Lemma \ref{leme2} that $\b_4$ and $\b_n$ form an edge pair.

This means we can write $\LL$ as a graph over $x_3, x_4$. The same argument applies to any $x_j, x_k$,
so condition (iii) is also sufficient.

For the last assertion of the theorem, we either have $\LL$ is a \wrs \ or it is not.

Case (a): The functions  $\{\a_3, \dots, \a_N\}$ form a weak edge family: all the $\o_{jk}$ in \eqref{Eq57} are the same, and all the $\eta_{jk}$ are the same.
The argument in \eqref{eqe107} shows that being a \wrs\ is invariant under change of coordinates.
By Lemma \ref{lemd1}, we can assume that we have a good parametrization. Now divide 
$\{ \a_j \}$ into those that are in $\M$, say for $3 \le j \le k$, and those that are in $b\M$, say
$k+1 \le j \le N$. 
Then 
\[
\{ (x_1 , \dots, x_k) \in \D^k : x_j = \a_j(x_1, x_2), 3\le j \le k \}
\]
is a \roms, and $\LL$ is a graph over this.
If none of the $\a_j$ are in $\M$, then they all extend to be holomorphic on $\D^2$, and we have a retract.

Case (b): There are exactly 3 distinct values of $\o_{jk}$ in \eqref{Eq57}. In this case $\LL$ is a graph over a \roms\ of height $1$.

We must have 3 different possible $\o$'s in \eqref{Eq57}. Call them $\o_a, \o_b, \o_c$,
with corresponding $\eta_a, \eta_b, \eta_c$.
Choose a good parametrization. Then we can relabel to assume $\hatt(\a) = (\o_a, \o_b) $.
Notice that $\o_c$ cannot be in the arc $(\o_a, \o_b)$, since it must lie outside of $\hatt(\a)$.
So for every $j$ we have $\hatt(\a_j)$ is one of the three arcs $(\o_a,\o_b), (\o_a, \o_c), (\o_c,\o_b)$.
Reorder the coordinates so that
$\hatt(\a_3) = (\o_a,\o_b)$. Any other $\a_j$ with $\hatt(\a_j) = \hatt(\a_3)$ must have $A_j = A_3$ and $B_j = B_3$,
since we have $A_j + B_j \o_a = \eta_a = A_3 + B_3 \o_a$ and
$A_j + B_j \o_b = \eta_b = A_3 + B_3 \o_b$. Therefore $\a_j$ is just a unimodular multiple of $\a_3$.
Similarly there is essentially only one function, $\a_4$ say, with $\hatt(\a_4) = (\o_a,\o_c)$ and one $\a_5$
with $\hatt(\a_5) = (\o_c,\o_b)$.

Let
\[
\LL' = \{ (x_1, x_2, \a_3(x_1, x_2) : 
 (x_1, x_2) \in D_{\a_3} \cap \dots \cap D_{\a_k} \}. 
\]
Then $\LL'$ is a \roms\ of height $1$.
We claim that for any $j > 3$, $|\a_j|$ is bounded by $1$ on $D_{\a_3}$.
Indeed, if $\a_j$ is a unimodular multiple of $\a_3$, this is immediate. Otherwise, 
$\hatt(\a_3) \subsetneq \hatt(\a_j)$.
Now move in a continuous path of M\"obius automorphisms, moving the origin to $(a_1, a_2)$, giving functions $\tilde \alpha_n$.
The new surface will still satisfy condition (i), hence also (ii) and (iii), and the points $\tilde \o_a, \tilde \o_b,
\tilde \o_c$ will vary continuously. They cannot coalesce, since if we reach Case (a) after an automorphism, we must have been in Case (a) all along.
Therefore we have $\hatt ( \tilde \a_3) \subsetneq \hatt (\tilde \a_j)$ for all M\"obius automorphisms.
Let $(a_1, a_2)$ tend to a point in $\D^2$ that lies in the boundary  of $D_{\a_j}$. By \eqref{eq21}, this means $\hatt (\tilde \a_j)$ shrinks;
so $\hatt (\tilde \a_3)$ shrinks; and so by \eqref{eq21} again, this means $|\a_3|$ must tend to $1$.
Therefore $D_{\a_3} \subseteq D_{\a_j}$, and hence $\LL$ is a graph over $\LL'$.
\ep

\begin{lemma}
\label{leme4}
Let $\LL$ be a weak \roms, and suppose it contains a point $z'$ with $|z'_1| > \max_{j \ge 2} |z'_j|$.
Then for some $j \ge 2$ it has a good $\{ x_1, x_j \}$ parametrization.
\end{lemma}
\bp
Suppose $\LL$ is as in \eqref{Eq42}. By Theorem \ref{thmb3}, for each $j \ge 2$ there are functions $\a_{jk} \in \overline{\M}$ so that $\LL$ can be written as a graph over $\{ x_1, x_j \}$ of the form
$\{ x \in \D^N : x_k = \a_{jk}(x_1, x_j), \ k \neq 1, j \} $.

 Let \[
\Omega_j =  \cup_{ k \neq 1,j} b \hatt(\a_{jk}) .
\]
Let 
\[
E = \cup_{j=2}^N \cup_{\o \in \O_j}  \{ z \in \D^N : z_j = \omega z_1 \} .
\]
Since $E$ is $(2N-2)$-dimensional, it does not disconnect $\D^N$.
Start at $z'$ and follow a path in $\LL \setminus E$ until reaching a point $z''$ where
$|z''_1| = |z''_j| \ge \max_{k \neq 1,j} |z''_k|$.

To see that such a path exists, start at $z'$ and head radially towards $0$. If this path ever reaches the boundary of $\LL$, which is where some $|z_j| = 1$, then there must have been some point on the path where $|z_j| = |z_1|$. Otherwise $z'$ is in the same component of $\LL$ as $0$, and near $0$ one can write $\LL$ as a graph over 
$\{x_1, x_2 \}$ which range over an open set containing a neighborhood of $0$ in $\D^2$.
Therefore there are points where $|z_2| > |z_1|$, so by continuity we have some point $z''$ as above.

Let $\om = z''_j/z''_1$. We claim that at $z''$ we have $\max_{k \neq 1,j} |z''_k| < |z''_1 | = |z''_j |$.
Indeed, this follows from the fact that for any $\a \in \overline{\M}$, and any $\l \in \D \setminus \{ 0 \}$, we have the trichotomy
\be
\label{eq65.1}
\notag
|\a(\l, \o \l ) | \
\begin{cases}
< |\l| & \o \in \hatt(\a)\\
= |\l| & \o \in b \hatt (\a) \\
> |\l| & \o \in \T \setminus \overline{\hatt(\a)}.
\end{cases}
\ee
We have ruled out the second possibility by avoiding $E$, so for each $k$ we have
$|\a_{jk} (z_1'', \o z_1'')| < |z_1''|$, and hence $\o \in \hatt(\a_{jk})$.
Therefore $\cap_{k\neq 1,j} \hatt (\a_{jk})$   is non-empty, and we have a good $\{x_1, x_j \}$ parametrization.
\ep
\begin{cor}
\label{core5}
Every \roms\ $\LL$ is connected.
\end{cor}
\bp
Suppose for some $i,j$ we have $|z_i | = |z_j|$ for all $z \in \LL$. Then there is some unimodular $\o$ so that $z_j = \o z_i$, and it is sufficient to prove the result for the 
$(N-1)$-dimensional \car set we get from dropping the $j^{\rm th}$ component. So let us assume that we do not have two coordinates that always have the same modulus.

Let $z' \in \LL$. We will show there is a path connecting it to $0$. Indeed, we can connect it to a nearby point with the property that one of the coordinates is strictly larger in modulus than all the others. By the argument in the proof of Lemma \ref{leme4}, we can then connect this to a point $z''$ that has exactly two equal largest components, which for notational convenience we will assume are the first two: $|z_1''| = |z_2''|  > \max_{k> 2} |z_k''|$.
Moreover, from the proof of the lemma again we get that $\om = z_2''/z_1''$ is in $\cap_{j=3}^N \hatt (\a_j)$, where
we write $\LL = \{ (z_1, z_2, \a_j(z_1, z_2) ) : (z_1, z_2) \in \cap_{j=3}^N D_{\a_j} \}$.
Therefore for every $0 \le r \le 1$, we also have $(rz_1, rz_2) \in  \cap_{j=3}^N D_{\a_j}$, and hence $z'$ is connected to $0$.
\ep

{\bf Theorem
\ref{thme3}.}
{\em
{\roms}s are totally geodesic, and are \car sets.}

{\sc Proof of Theorem \ref{thme3}:} We shall take advantage of the fact that being a \wrs\ is invariant under
all automorphisms of the polydisk that map some point of the surface to the origin.
Let $\LL$ be as in \eqref{Eq42}, with $\{ \a_j \}$ an edge family. Applying a rotation automorphism,
we can assume each $\a_j = [1; A_j, B_j]$, and since they form an edge family we have unimodular $\om_0, \eta_0$ so that
\be
\notag
\label{Eq66.0}
A_j + B_j \om_0 = \eta_0 , \qquad 1 \le j \le m .
\ee

Step 1. We shall show that $\LL$ is totally geodesic in the infinitesimal sense. 

Let $a \in \LL$ and $X \in T_a \LL$, the complex tangent space. We will show that there is a complex geodesic in $\D^{m+2}$ with respect to $(a;X)$ that lies totally in $\LL$. By applying an automorphism, we can assume that $a=0$.
 Permute the coordinates so that $|X_1| \geq |X_2| \geq \dots$, and rescale so $X_1 = 1$ and $X_2 = c \in \overline{\D}$; then $X = (1, c, A_1 + B_1 c, \dots, A_m + B_m c )$.

Let \[
\Gamma \= \{ c \in \overline{\D} : |A_j + B_j c | \le 1,\ 1\le j \le m \} .
\]
For $c$ in $\Gamma$, we are looking for a geodesic of the form
\be
\label{Eq66}
f_c: \l \ \mapsto \ ( \l, \l \vp_{-c}(\om \l), \a_1(\l, \l \vp_{-c}(\om \l)), \dots,  \a_m(\l, \l \vp_{-c}(\om \l)) ),
\ee
where $\om = \om(c)$ is an appropriately chosen point in $\T$.
 It is obvious that the first coordinate function is a left inverse
for $f_c$, so it is a geodesic provided $f_c(\D) \subseteq \LL$.

Case: $|c| < 1$.
Write
\[
\a_j(\l, \l \vp_{-c}(\om \l)) \= \l \frac{ A_j + B_j \vp_{-c}(\om \l) + \l \vp_{-c}(\om \l)}{1+ \bar B_j \l + \bar A_j  \l \vp_{-c}(\om \l)}
\= \l \frac{p_j(\l)}{q_j(\l)} ,
\]
where $p_j$ and $q_j$ are the polynomials obtained by multiplying the numerator and denominator by $(1+ \bar c \om \l)$.
\beq
p_j(\l) &=& A_j(1+ \bar c \om \l) + B_j(c +\om \l) + \l(c + \om \l) \\
q_j(\l) &=& (1+ \bar c \om \l) + \bar B_j \l (1 + \bar c \om \l) + \bar A_j \l(c + \om \l) .
\eeq
Define unimodular $\l_0$ and $\o$ by
$\l_0 = - \bar \om_0 \eta_0$, and $\om = \bar \l_0 \vp_{c}(\om_0)$. Then we have
\beq
\vp_{-c}(\om \l_0) &\=& \om_0 \\
\l_0 \vp_{-c}(\om \l_0) &=& - \eta_0.
\eeq
It follows from these equations that $\l_0$ is a root of both $p_j$ and $q_j$.
Moreover, the product of the two roots of $p_j$ is $\bar\om(A_j + B_j c)$, and 
$\om \l^2 \overline{p_j(1/\bar \l)} = q_j(\l)$. So if $ |A_j + B_j c| < 1$,
then $\l p_j(\l)/q_j(\l)$ is a rational inner function of degree  $2$ with one  root at $0$ and the other inside the disk; if $ |A_j + B_j c| = 1$ then $ \l p_j(\l)/q_j(\l)  = (A_j + B_j c) \l$ is a unimodular multiple of $\l$.
Either way, $|\a_j(\l, \l  \vp_{-c}(\o \l))| < 1 $ for all $\l \in \D$.

Case: $|c|=1$. We set $\vp_{-c}(\om \l) = c$ for all $\l$, and $f_c$ becomes
\[
f_c(\l) = ( \l, c \l, \a_1(\l, c\l), \dots, \a_m (\l, c\l)).
\]
 Again, $\a_j(\l, c\l)$ is a unimodular multiple either of $\l$ or 
a Blaschke product of degree 2,
 so stays in $\D$ for
 all $\l$ in $\D$.

Step 2: We prove that $\LL$ is totally geodesic.

We wish to prove that for any two distinct points in $\LL$, there is a geodesic going through both.
Without loss of generality, we can assume that one point is $0$, and the second point $z$ satisfies
$|z_1| \ge |z_j|$ for every $2 \le j \le m+2$.
Indeed, we can make the stronger assumption that $|z_1| > |z_j|$ for every $ j \ge 2$, since if we
prove the assertion for a dense set, a normal families argument gives it for the general case;
and for any $i,j$ the set $\{ x \in \LL : |x_i| \neq |x_j| \}$ is an open dense set.
By Lemma \ref{leme4} there is a good $\{x_1, x_j \}$ parametrization for some $j$, which without
loss of generality we will assume is $j=2$.

We would like to find $c \in \Gamma$ so that $z_1 \vp_{-c}(\om(c) z_1) = z_2$; then $f_c$ from \eqref{Eq66} is the desired geodesic.

Let 
\[
Z \= \{ \xi \in \overline{\D} : |\a_j(z_1, \xi z_1)| \le |z_1|, 1 \le j \le m \}.
\]
Define $\Psi: \Gamma \to Z$ by
\[
\Psi (c) \= ((f_c(z_1))_{2} /z_1 ) = \vp_{-c} (\o(c) z_1 ),
\]
where $\o(c) = - \o_0 \bar \eta_0 \vp_c(\o_0) $ as in Step 1.
The map $\Psi$ is continuous, indeed it is a rational function of $c$ and $\bar c$,
and by Step 1 we know it maps into $Z$.
We wish to prove that $\Psi$ is surjective. 

$\Gamma$ is given by the intersection of $m+1$ disks, so is convex, and hence simply connected.
Moreover, $\Gamma$ contains $\cap_{j=1}^m T(\a_j)$ which is a non-empty open subset of $\T$, since we have a good $\{x_1, x_2 \}$ parametrization, so $\Gamma$ has interior.

$Z$ is given by the intersection of $m+1$ sets,
\[
Z = \overline \D \cap_{j=1}^m \{ \xi : |A_j +(B_j + z_1)\xi|^2  \le | (1 + \bar B_j z_1) + \bar A_j z_1 \xi|^2 \} .
\]
These are either disks (if $|B_j + z_1|^2 > |A_j z_1|^2$) 
 or the complement of a disk (if the opposite inequality holds), or a half-plane (if we have equality).
 None of them reduce to a point, since we have $\xi = z_2/z_1$ is in the interior of each set by the
 hypothesis that $|z_1| > |z_j|$ for all $j \ge 2$.

We make three claims:

(i) $\Psi(c) $ is in the boundary of $Z$ if and only if $c$ is in the boundary of $\Gamma$.

(ii) $\Psi$ is one-to-one on $b \Gamma$.

(iii) The interior of  $Z$ is connected.

Once these three claims are established, by the Jordan curve theorem $\Psi$ is a homeomorphism from
$b \Gamma$ to $bZ$ and maps the interior of $\Gamma$ onto the interior of $Z$.

To see (i), note that if $|c| = 1$, then $c \in b\Gamma$ if and only if $|A_j + B_j c | \le 1$ for all $j$.
Then $\Psi(c) = c$, and $|\a_j(z_1, c z_1)| \le |z_1|$.

If $|c| < 1$, then  $c \in b\Gamma$ if and only if $|A_j + B_j c | = 1$ for some $j$,
and $|A_k + B_k c | \le 1$ for all $k$.
The first equality  happens if and only if $\a_j(\l, \l \vp_{-c}(\o \l))$, which maps $\D$ to $\D$ by Step 1,
satisfies $\a_j(\l, \l \vp_{-c}(\o \l)) = (A_j + B_j c) \l$ by the Schwarz lemma.
So
$| \a_j( z_1, \Psi(c) z_1) | = |z_1|$ if and only if $|A_j + B_j c | = 1$.
Hence $\Psi(c) \in b Z$ if and only if $c \in b\Gamma$.

To see (ii), it is clear if $|c|=1$. If $|c|<1$, suppose $|A_j +B_j c| = 1$ and $\Psi(c) = \xi$.
Then $\a_j (z_1, \xi z_1) = (A_j + B_j c) z_1$. If any other $c'$ satisfies $\Psi(c') = \xi$, then 
we would have $(A_j + B_j c')z_1 = (A_j + B_j c) z_1$. so $c' = c$.

To see (iii), define
 $U_0 := \D$ and $U_j := \{ \xi : |\a_j(z_1, \xi z_1)| < |z_1|\} $.
We have  ${\rm int}(Z) = \cap_{j=0}^m U_j$, and the  arc $\hatt(\a)$ is in $bU_j$ for each $j$.
Now apply the map $\xi \mapsto \frac{1}{\xi-\o_0}$. Then each circle through $\o_0$ is mapped to a line,
and the sets $U_j$ get mapped to half-planes. The intersection of these half-planes
  contains the point to which $z_2/z_1$ gets mapped, so is a nonempty open convex set. Its inverse image
  is a connected open set with $\hatt(\a)$ in its boundary.

Step 3: To prove $\LL$ is also a \car set, we use Proposition \ref{propc3}.
\ep

{\bf Theorem
\ref{thme5}.}
{\em
Let $\LL$ be a \roms\ as in \eqref{Eq7}. Then $\LL$ is a Lempert set, and 
 $ D_{\a_3} \cap 
\dots \cap D_{\a_{m+2}}$ is a two-dimensional Lempert domain.
}

{\sc Proof of Theorem \ref{thme5}:} This follows from Theorem \ref{thme3}.
To prove that $\LL$ is a Lempert set, we must show that $l_\LL(\mu,\nu)  = d_\LL(\mu,\nu)$
for any points $\mu, \nu \in \LL$.
Since $\LL$ is totally geodesic, there is some holomorphic $f : \D \to \LL$  that maps
$a$ to $\mu$ and $b$ to $\nu$ for some $a,b \in \D$.
Then
\beq
d_\LL(\mu,\nu)  &\ = \ & \sup \{ \rho(\phi (\mu), \phi(\nu)) : \phi \in O(\LL, \D) \} \\
&=&
 \sup \{ \rho(\phi (f(a)), \phi(f(b))) : \phi \in O(\LL, \D) \}
\\
&\le& \rho(a,b).
\eeq
Taking the infimum over all possible choices of $f$ we get $d_\LL(\mu, \nu) \le l_\LL (\mu, \nu)$.
But we know there is some $f$ that has a left inverse $\psi$, so we get for this $f$
\beq
l_\LL (\mu, \nu) &\ \le \ & \rho(a,b) \\
&=& \rho( \psi(f(a)), \psi(f(b))) \\
&=& \rho(\psi(\mu), \psi(\nu)) \\
&\le& d_\LL(\mu, \nu) .
\eeq
Therefore $\LL$ is Lempert.

Since $\LL$ is a graph over $U := \cap_{j=1}^m D_{\a_j}$ of a function into $\D^m$, it is easy to check that 
$l_U(x,y) = l_\LL ( (x,\a(x)), (y,\a(y)))$ and
$d_U(x,y) = d_\LL ( (x,\a(x)), (y,\a(y)))$. So $U$ is also Lempert.
\ep

\section{Proof of Lemma \ref{lemd2}}
\label{secg}

It is well-known that the \car and Kobayashi metrics agree on the polydisk; this can be seen by direct computation,
or invoking Lempert's theorem.
They are given by
\be
\label{eq39}
d_{\D^N} (\l, \nu) \= \max_{1 \leq i \leq N} \rho(\l_i, \nu_i) .
\ee
Unlike in the case of strictly convex domains, complex geodesics that go through a pair of points $\l, \nu$ in $\D^N$ are rarely unique; but they are unique on the coordinates where the maximum in \eqref{eq39} is attained.
\begin{defin}
\label{DEFg1}
Let $\mu, \nu$ be distinct points in $\D^N$. We say they are $k$-balanced with respect to $\{ i_1, \dots, i_k \}$
for some set of $k$ numbers in $\{1, \dots , N\}$ if 
$d_{\D^N} (\mu,\nu) = \rho(\mu_{i_j}, \nu_{i_j})$ for $j \in \{ 1, \dots, k \}$.
\end{defin}
Suppose $(\mu,\nu)$ is $k$-balanced with respect to $\{1, \dots, k\}$.
Any complex geodesic $f : \D \to \D^N$ that maps $(a,b)$ to  $(\mu,\nu)$  has $f_i$ uniquely determined
for $1 \leq i \leq k$, namely it is the unique automorphism of $\D$ that sends $(a,b)$ to $(\mu_i, \nu_i)$.

For a proof of the following lemma see e.g. \cite[Lemma 9.3]{aly19}.
\begin{lemma}
\label{lemg1}
Let $D = \D \setminus \D(w, r)$ for some $w \in \D$ and positive $r < 1-|w|$.
Let $a,b$ be points in $D$. Then there is a map $g: D \to \D$ so that
$\rho(a,b) < \rho(g(a),g(b))$.
\end{lemma}
As a consequence, we have the following lemma, which dates back to an idea of P. Thomas (see \cite[Appendix by P.J. Thomas]{agmcvn}).
We say $\phi$ is a \car extremal if it attains the supremum in \eqref{Eq1}.
\begin{lemma}
\label{lemg2}
Let $V$ be a \car set, and $\phi: V \to \D$ a \car extremal.
Then $\phi(V)$ is dense in $\D$.
\end{lemma}
\bp
Suppose $\phi$ is an extremal for the pair $(\mu,\nu)$.
If the conclusion were false,  $\phi(V)$ would be contained in some domain $D$ as in Lemma \ref{lemg1}.
But then $g \circ \phi$ would also map $V$ to $\D$, and would separate $\mu$ and $\nu$ by more than $\phi$ did, contradicting its extremality.
\ep

\begin{lemma}
\label{lemg9}
Suppose $0$ and $x'$ are  in $V$.
Assume that $x'$ is a regular point of $V$ and that the germ of $V$ at $x'$ can be written
as
\be
\label{eq44}
\notag
(x_1, \dots, x_\dimv, \beta_{\dimv+1}(x), \dots, \beta_{\dimpol}(x) ). 
\ee
If $\tau \in \T^{\dimv+1}$ satisfies:
there exists $\l \neq 0$ so that
$\beta_{\dimv+1}(\tau_1 \l, \dots, \tau_\dimv \l) = \tau_{\dimv+1}\l$ and $|\beta_{j}(\tau_1 \l, \dots, \tau_\dimv \l)| \le |\l|$ for
all $j \geq \dimv+2$,
then $\tau \in \overline{\pi_{\dimv+1}(V)}$.
\end{lemma}
\bp
Let \[
\phi_\tau(x_1, \dots, x_N) \ := \ \frac{1}{\dimv+1} ( \overline{\tau_1} x_1 + \dots + \overline{\tau_{\dimv+1}} x_{\dimv+1}) .
\]
Then $\phi_\tau(0) = 0$ and $\phi_\tau(\tau_1 \l, \dots, \tau_{\dimv +1} \l)  = \l$, so $\phi_\tau$ is a \car extremal for $V$.
By Lemma \ref{lemg2}, this means $\phi_\tau(V)$ is dense in $\D$,
and hence $(\tau_1, \dots, \tau_{\dimv+1})$ is in $\overline{\pi_{\dimv+1}(V)}$.
\ep

\begin{lemma}
\label{lemg3}
Suppose $0$ and $x'$ are in $V$, and that
$x'$ satisfies
\be
\label{eq40}
|x_1'| = |x_2' |= \dots = |x_{d+1}'| >| x_j'| \geq 0 \qquad \forall\ d+2 \leq j \leq N .
\ee
If $d+1 = N$ we interpret \eqref{eq40} as saying all the coordinates of $x'$ having equal modulus.
Assume that $x'$ is a regular point of $V$ and that the germ of $V$ at $x'$ can be written
as
\be
\label{eq41}
(x_1, \dots, x_d, \beta_{d+1}(x), \dots, \beta_{\dimpol}(x) ). 
\ee
Then $\beta_{d+1}(\frac{x_1'}{|x_1'|} \l , \dots,\frac{x_d'}{|x_d'|} \l)$ is a constant multiple of $\l$ in some neighborhood of $|x_1'|$.
\end{lemma}
\bp
Composing with a rotation automorphism, we can assume that all the coordinates $x_j'$ are non-negative real numbers.
Let $g(\l) = \beta_{d+1}(\l, \dots, \l)/\l$. Hypothesis \eqref{eq40} says $g(x_1') = 1$.
Assume $g$ is non-constant; then it is open.
For $\tau \in \T^d$ and close to $(1, \dots, 1)$, define
\[
h_\tau(\l) = \beta_{d+1}(\tau_1 \l, \dots, \tau_d \l)/\l .
\]
From Hurwitz's theorem, there is a neighborhood $J$ of $1$ in $\T$ and a neighborhood $U$ of $x_1'$ in $\D$ 
so that for all $\tau$ sufficiently close to $(1, \dots, 1)$ and all $\tau_{d+1} $ in $J$, there is some
point $\l \in U$ so that $h_{\tau}(\l) = \tau_{\dimv+1}$.
Using the representation \eqref{eq41}, this means
if $(\tau_1, \dots, \tau_\dimv, \tau_{\dimv+1})$ is sufficiently close to $(1, \dots , 1)$ in $\T^{\dimv+1}$, then $V$ contains points of
the form
\be
\label{eq42}
\notag
\gamma_\tau(\l) := (\tau_1 \l, \dots, \tau_\dimv \l, \tau_{\dimv+1} \l, \beta_{\dimv+2} (\tau_1 \l, \dots \tau_{\dimv} \l), \dots, \beta_N(\tau_1 \l, \dots \tau_{\dimv} \l) )
\ee
where for each $j \ge \dimv+2$ we have $|\beta_j( \tau_1 \l, \dots \tau_{\dimv} \l)| < |\l|$.
(Note that the choice of $\l$ depends on $\tau$.)
Therefore $V$ contains a $(\dimv+1)$ balanced pair $(0, \gamma_\tau(\l))$.

By Lemma \ref{lemg9}, 
 $(\tau_1, \dots, \tau_{\dimv+1})$ is in $\overline{\pi_{\dimv+1}(V)}$,
for all $\tau$ in a neighborhood of $(1, \dots , 1)$
in $\T^{\dimv+1}$. By  Chevalley's Theorem \ref{thmc7} this is impossible since $V$ is $\dimv$-dimensional.
\ep

In Section \ref{secm} we want a variant of Lemma \ref{lemg3} where we drop the assumption of irreducibility.
\begin{lemma}
\label{lemg7}
Let $V$ be a \roms\ in $\D^N$ parametrized as $\{ (x, \alpha(x))  : x \in U \subseteq \D^\dimv \}$.
Assume $P \in V$ and $b \notin V$ are $(\dimv+1)$-balanced with respect to the first $\dimv+1$ coordinates.
Define $\tau_j$ to be the phase of $\vp_{b_j}(P_j)$. Let $W$ be a \car set that contains $V \cup \{ b\}$.
Either $W$ is more than $\dimv$-dimensional, or 
\be
\label{eq43}
\a_{\dimv+1} (\vp_{b_i}(\tau_i \l) ) = \vp_{b_{\dimv+1}}(\tau_{\dimv+1} \l) 
\ee
for all $\l \in \D$ with $ (\vp_{b_i}(\tau_i \l) ) \in U$.
\end{lemma}
\bp Note first that Lemma \ref{lemg3} can be rewritten as Lemma \ref{lemg6}, where the proof is identical except the last sentence becomes a dichotomy. Then \eqref{eq43} holds in a neighborhood of $|m_{b_1}(P_1)|$. By Corollary \ref{core5}, 
$U$ is connected, so the equality holds on all of $U$.
\ep

\begin{lemma}
\label{lemg6}
Let $W$ be a \car set in $\D^N$ containing $0$. Assume that there is a regular point $x' \neq 0 \in W$ such that locally
the germ at $x'$ can be written as  $(x_1, \dots, x_\dimv, \beta_{\dimv+1}(x), \dots, \beta_{\dimpol}(x) )$.
Suppose 
\[
|x_1'| = |x_2' |= \dots = |x_{\dimv+1}'| >| x_j'| \geq 0 \qquad \forall\ \dimv+2 \leq j \leq N .
\]
Then either  $\beta_{\dimv+1}(\frac{x_1'}{|x_1'|} \l , \dots,\frac{x_\dimv'}{|x_\dimv'|} \l)$ is a constant multiple of $\l$ in some neighborhood of $|x_1'|$ or
${\rm dim}(W) > \dimv$.
\end{lemma}

We can now prove Lemma \ref{lemd2}.

{\bf Lemma
\ref{lemd2}.}
{\em
Suppose $0$ is a regular point for $V$.
If $\om$ is in $\overline{T(\alpha)}$, then  for each $j$ the function $\l \mapsto
\alpha_j(\omegadotlambda) $ extends to a holomorphic map from $\D$ to $\D$ that is continuous from
$\overline{\D}$ to $\overline{\D}$.
Moreover, this extension is continuous on $W = \{ \omegadotlambda : \lambda \in \overline{\D}, \ \omega \in \overline{T(\alpha)} \}$ except possibly on a set of dimension $d-1$.
}

\bp
 We will show (i)   $\alpha_j(\omegadotlambda)$ is holomorphic on $\D$ 
if $\om$ is in a dense subset of $T(\a)$;
(ii) this remains true for all  $\om \in \overline{T(\a)}$;
(iii)  the function $\lambda \mapsto \alpha_j(\omegadotlambda)$
 is in the disk algebra (i.e.  it extends to be continuous on $\overline{\D}$).

(i) First assume that $\om \in T(\a)$. Define a one-dimensional set by
\be
\label{eq45}
\notag
W_\o \ := \left[ \{ (\o_1 \l, \o_2 \l, \dots , \o_{\dimv} \l) : \l \in \D \} \times \C^{N-\dimv} \right] \cap V .
\ee
The function $\alpha_j(\omegadotlambda)$ is holomorphic if we can show that
\be
\label{eq46}
W_\o \setminus \{ 0\} \subset \{ x: \ \forall\ j > \dimv, \ |x_j | < |x_1|  \} .
\ee
Indeed, if \eqref{eq46} holds, then the projection from $W_\o$ onto the first component is proper. 
This means that $W_\o$ must be single-sheeted over $x_1$,
and so $\alpha_j(\omegadotlambda) $
extends to be holomorphic on all of $\D$, and takes values in $\D$.

Suppose \eqref{eq46} does not hold. By Theorem \ref{thmc8} each of the functions $\alpha_j(\omega \cdot \l)$ is Nash,
so satisfies some polynomial equation $p_j(\l, \alpha_j(\omega \cdot \l)) =0$.
Let $D_j$ be the finite set of points where the discriminant of $p_j(x,y)$ (with respect to $y$) vanishes.
Then $\alpha_j(\omega \cdot \l)$ can be analytically continued to any simply connected subset of $\C \setminus D_j$.
We can find $\l \in \D$ so that $(\o\cdot \l, \a(\o \cdot \l))$ is a regular point of $V$ and
\be
\label{eq47}
\max_{j > \dimv} |\a_j(\o\cdot \l) | = |\l|.
\ee
If there is a unique $j > \dimv$ that attains the maximum in \eqref{eq47}, 
then Lemma \ref{lemg3} says $\a_j(\o \cdot \l) = \tau \l$ for some $\tau \in \T$.
Differentiating this equation at $0$ we get $| \sum_{i=1}^{\dimv} \al^i_j \, \o_i | =1$,
violating the condition that $\o \in T(\al_j)$.

So for the remainder of the proof of (i) we can assume there is an open set $U$ in $T(\a)$ such
that for all $\o$ in $U$ \eqref{eq46} fails and 
 for all points $\l$ where \eqref{eq47} holds, 
the maximum is never attained for a unique $j$.

Let us fix some $\o' \in U$ and a point $\l' \in \D\setminus\{0 \}$ so that $(\o' \cdot \l', \a(\o' \cdot \l'))$ is 
regular for $V$ and $\max_{j > \dimv} |\a_j(\o'\cdot \l') | = |\l'|.$
Without loss of generality, we can assume this maximum occurs for $j = \dimv+1$.
Let
\[
\tau_{\dimv+1}' \= \frac{\a_{\dimv+1}(\o'\cdot \l')}{\l'} .
\]
Since $\a_{\dimv+1}(\o'\cdot \l')$ cannot be a unimodular multiple of a coordinate function (as this would violate $\o' \in T(\a)$), by Hurwitz's theorem we have that for all $\tau_{\dimv+1}$ close
to $\tau'_{\dimv+1}$ there exists $\l$ close to $\l'$ such that $\a_{\dimv+1} (\o' \cdot \l) =
\tau_{\dimv+1} \l$.
Fix such a $\tau_{\dimv+1}$ and $\l$. Then there exist sequences $\l_{(i)} \to \l$ and
$\o_{(i)} \to \o$ and a set $S \subseteq \{ \dimv+1, \dots , \dimpol \}$ so that
\be
\label{eq48}
\max_{j > \dimv} |\a_j(\o_{(i)} \cdot \l_{(i)}) | = |\l_{(i)}|
\ee
and so that the maximum in \eqref{eq48} occurs precisely when $j \in S$.
By Lemma \ref{lemg9} we conclude that for all $\tau_{\dimv+1}$ close to $\tau'_{\dimv+1}$ 
the point $(\o', \tau_{\dimv+1})$ is in $\overline{\pi_{\dimv+1}(V)}$. 
The same argument works if we perturb $\o'$.
But this would mean that $\overline{\pi_{\dimv+1}(V)} \cap \T^{\dimv+1}$ is not contained in a $\dimv$ (real) dimensional subset of the $(\dimv+1)$ torus, contradicting Chevalley's theorem.

(ii)
We have established  the lemma for a dense set of  points ${\mathcal T}$ in $T(\a)$. 
Now assume that $\om \in \overline{T(\a)}$. Let $\o_N \in {\mathcal T}$ converge to $\o$.
Then the functions $\a_j(\o_N \cdot \l)$ will converge to a holomorphic function that
maps $\D$ to $\D$, and agrees with $\a_j(\o \cdot \l)$ on a neighborhood of the origin and hence
on the whole disk.

(iii) Let $ f(\l) := \alpha_j(\omegadotlambda)$. To see that $f$ is in the disk algebra, we use the fact that
$p_j(\l, f(\l)) = 0$. 
Let $\l_0 \in \T$. Let $\l_k \in \D$ be a sequence converging to $\l_0$, and suppose $f(\l_k)$ converges to some limit $w$.
Then $p_j(\l_k, f(\l_k)) =0$ converges to $p_j(\l_0, w)$; so $w$ must be one of the finitely many roots of $p_j(\l_0, \cdot)$.
Let these roots be $w_1, \dots, w_\dimv$, and let $\varepsilon > 0$ be small enough that the closed disks of center $w_i$ and radius
$\varepsilon$ are disjoint. There exists $\delta > 0$ so that if $\l \in \D$ and $|\l - \l_0| < \delta$, then $f(\l) \in
\cup_{i=1}^\dimv \D(w_i, \varepsilon)$, since otherwise one could find a sequence $\l_k \to \l_0$ so that $f(\l_k)$ converges to some limit that is not one of the $w_i$. But since $f$ is continuous on $\D$, it must be that there exists a single $w_i$ so that
$f( \D \cap \D(\l_0, \delta)) \subseteq \D(w_i , \varepsilon)$.
As $\varepsilon$ can be made arbitrarily small, we conclude that $f$ is continuous at $\l_0$.

(iv) Similarly, 
to see that $\alpha_j$ is continuous at $\l_0\cdot \omega_0 \in W$, let $\l_k \cdot \omega_k \in W$ converge to $\l_0\cdot \omega_0$,
and suppose $\al_j (\l_k \cdot \omega_Nk$ converges to $w$.
As $\al_j$ is a Nash function, 
there is some polynomial $Q$ so that $Q(x, \al_j(x))= 0$ in a neighborhood of $0$ where $\al_j$ is holomorphic, and hence wherever $\al_j$ extends holomorphically. 
Therefore $Q(\l \cdot \om , \al_j(\l \cdot \om)) = 0$ for $\l \in \D$ and $\om \in \overline{T(\al)}$.
(Note that $Q$ does not depend on $\om$).
Therefore $Q(\l_k \cdot \omega_k, \al_j(\l_k \cdot \omega_k )= 0$.

As $Q$ is continuous,  $Q( \l_0\cdot \omega_0, w) = 0$, so we have that $w$ is one of the zeroes of $Q (\l_0\cdot \omega_0, \cdot)$.
It is possible that $Q (\l_0\cdot \omega_0, \cdot)$ is identically $0$; but this can only occur for $\l_0 \cdot \omega_0$ lying in a $d-1$ dimensional algebraic set. 
Otherwise, we conclude as before that $\al_j$ is continuous at $\l_0\cdot \omega_0$.
\ep

\section{Proof of Lemma \ref{lemd3}}
\label{sech}

\begin{lemma}
\label{lemh1}
Let $\om \in b T(\a) \cap b T(\a_j)$. Then there exists $\eta \in \T$
so that for all $\lambda \in \D$,  $\a_j(\omegadotlambda) = \eta \l$.
\end{lemma}
\bp
By Lemma \ref{lemd2}, we know that $\a_j(\omegadotlambda) \in \hol(\D,\D)$, and maps $0$ to $0$.
Since
$\om  \in b T(\a_j)$, the derivative of the function is unimodular at $0$, and so the result follows
from the Schwarz lemma.
\ep

\begin{lemma}
\label{lemh2}
Let $c_1, \dots, c_\dimv \in \C$. The set
$ \{ | \sum_{i=1}^\dimv c_i \om_i | : \om \in \T^\dimv \}$ is the interval $[s,t]$,
where $t = \sum_{i=1}^\dimv |c_i|$ and $s = \max(0, 2 \max_{1 \leq i \leq k} (|c_i|) -  \sum_{j=1}^\dimv |c_j|)$.
\end{lemma}
\bp
Without loss of generality, assume that all the $c_i$ are positive and decreasing.
The set $\{ \sum_{i=1}^2 c_i \om_i \}$ is the closed annulus $\{ z : c_1 - c_2 \leq |z| \leq c_1 + c_2 \}$.
(If  $c_1 = c_2$ we get a closed disk).
Adding each extra $\{ c_j \om_j: \om_j \in \T \}$ in succession, one gets fatter annuli. The outer radius will increase by $c_j$,
the inner one will decrease by $c_j$ (unless it reaches $0$).
\ep
We won't explicitly use Lemma \ref{lemh2}, but it helps conceptually.
Note in particular that any complex number with modulus in the interval $(s,t)$ can be written as a sum $\sum c_i \om_i$ with
$k-1$ degrees of freedom in choosing the numbers $\om_i$.

\begin{lemma}
\label{lemh3} 
Let $P \in \C[x_1, \dots,  x_\dimv][ z]$ be irreducible, and assume that the degree of $P$ with respect to $z$ is $m$.
Let $W$ be an open subset of $\T^\dimv$.
Let $f_1, \dots, f_m$ be the branches over $W$, counting multiplicity (so $P(x,f_j(x)) = 0$ for
$x \in W$). Suppose that $|f_1| = 1$ on $W \cap \T^\dimv$.
Then $f_1 \cdots f_m$ is meromorphic and has modulus $1$ a.e. on $\T^\dimv$.
\end{lemma}
\bp
Let $\gamma(z) = i \frac{1+z}{1-z}$ be the Cayley transform, and $\delta(z) = \frac{z-i}{z+i}$ be its inverse.
The function $R(x,z) = P(\delta(x_1), \dots, \delta(x_\dimv), \delta(z))$ is a rational function;
let $Q(x,z)$ be its numerator in simplest form. Notice that $Q$ must also be irreducible, since if $Q$ 
factored as $Q_1 Q_2$ then the zero set of $P$ would be the union of the zero sets of $Q_1 \circ \gamma$
and $Q_2 \circ \gamma$.

Let $g_j = \gamma \circ f_j \circ \delta$, and let $W' = \gamma(W) \subseteq \R^\dimv$.
Then the hypotheses of the Lemma say that $g_1$ is real valued on $W'$.
Write $Q = q_1 + i q_2$, where $q_1, q_2 \in  \R[x_1, \dots,  x_\dimv][ z]$.
Since $g_1$ is real-valued on $W'$, we have $q_1(g_1)$ and $q_2(g_1)$ vanish on $W'$.
Therefore the zero set of $Q$ is contained in both the zero set of $q_1$ and the zero set of $q_2$.
Since the degrees of both $q_1$ and $q_2$ are at most the degree of $Q$, and this is irreducible,
we get that $q_1$ and $q_2$ must be constant multiples of $Q$ (perhaps zero). 
This means that $Q$ is a constant multiple of a polynomial in $\R[x][z]$.
Therefore each $g_j$ is  either real or occurs in a complex-conjugate pair.
The corresponding $f_j$ is then either unimodular, or occurs in a pair whose product is unimodular.
Therefore $f_1 \dots f_m$ is unimodular on $W$.
Since $f_1 \dots f_m$ is the product of all the roots of $P$, it is a rational function of the coefficients, hence rational in $x$.
\ep

\begin{lemma}
\label{lemh3.5}
Suppose there exists $\vare > 0$ so that, for all $a \in V \cap B(0,\vare)$, if we write as in \eqref{eq20}
\[
\tilde \a_j(y_1, \dots, y_\dimv) \ = \vp_{a_j} \circ \a_j(\vp_{-a_1}(y_1), \dots, \vp_{-a_\dimv}(y_\dimv) ) ,
\]
then there is an open set $U$ in $\T^d$ such that $U \cap b\overline{T( \tilde \al_{\dimv +1}) } $ is non-empty and 
 for all $\o \in  b\overline{T( \tilde \al_{\dimv +1}) } \cap U $ there exists $\eta \in \T$ so that $\tilde \a_{d+1} (\omegadotlambda) = \eta \l$ for all $\l \in \D$.

 Then $\al_{\dimv +1}$ is a rational function, and has modulus one a.e. on $\T^\dimv$.
\end{lemma}
\bp
Since $  b\overline{T(  \al_{\dimv +1}) } $ is non-empty, the interval $(s,t)$ from Lemma \ref{lemh2}, with $c_i = \frac{\partial \a_{\dimv+1}}{\partial x_i} (0)$,
has $s < 1 < t $.
Write $\beta$ for $\al_{\dimv +1}$.  Let us write $\Phi_\mu(x) = (\vp_{\mu_1}(x_1), \dots, \vp_{\mu_\dimv}(x_\dimv))$ and 
\be
\label{eq49}
\notag
\b_\mu(x) \= \vp_{\b(\mu)} \circ \b \circ \Phi_{-\mu}(x) .
\ee
There is an open set $G$ in $\C^\dimv$ so that 
 for all $\mu$ in $G$ and for all  $\tau$ in $bT(\b_\mu) \cap U  $  there exists $\eta \in \T$  (depending on $\mu$ and $\tau$) so that
\be
\label{eq50}
\b_\mu(\taudotlambda) \= \eta \l .
\ee
We can assume that $U$ is invariant under multiplication by unimodular scalars (since every  $T(\b_\mu)$ is).

Since $U = - U$ and
 $\lim_{\mu \to 0} \Psi_{\mu}(x) = x$, we can shrink $G$ to $G'$, a smaller neighborhood of $0$,  so that
\[
U' \= \cap_{\mu \in G'} \Phi_{-\mu}(U) 
\]
is non-empty.
Define
\begin{eqnarray}
\label{EQh1} \notag
X &\ :=\ &
\bigcup_{\mu \in G' }  bT(\b_\mu) 
\cap U'
 \subseteq\ \T^\dimv\\
\notag
Y &\ :=\ & \{ \Phi_{-\mu}(\taudotlambda ): \l \in \D,\ \mu \in G', \tau \in  bT(\b_\mu) 
\cap U' 
 \} \ \subseteq \ \D^\dimv 
\\ 
\notag
Z &\ :=\ &\bigcup_{\mu \in G' }  T(\b_\mu)
 \cap U'
 \subseteq\ \T^\dimv
\end{eqnarray}
Claim:  $X$ has non-empty interior in $\T^\dimv$, and $\beta$ is unimodular 
almost everywhere 
on $X$.

Proof of Claim:
Let $\om \in  bT(\b)$.
Let $F: \D^\dimv \times \T^\dimv \to \R$ be defined by
\[
F(\mu,\xi) \ := \  \left| \sum_{i=1}^\dimv \frac{1-|\mu_i|^2}{1 - |\beta(\mu)|^2} \frac{\partial \beta}{\partial x_i}\Big |_\mu
\xi_i
\right|^2 .
\]
Then 
\[
X \=
\cup_{\mu \in G'} 
 \{ \xi \in U': F(\mu,\xi) = 1 \} .
\]
The function $F(\mu,\xi)$ is real analytic in $\mu$ and $\xi$, and $F(0,\o) = 1$ for all $\o$ in $U'$. If its
derivative at $(0,\om)$ with respect to $\mu$ is non-zero,
the implicit function theorem says for any choice of $\xi$ in a neighborhood of $\om$
one can find $\mu$ in $G'$ so that $F(\mu,\xi) = 1$. This would imply that $X$ has interior.

Calculating, we get 
the following for $\frac{\partial F}{\partial \mu_k}$, where we use superscripts to denote partial derivatives of $\beta$ at $0$:
\be
\label{eq51}
\notag
\frac{\partial F}{\partial \mu_k} |_{(0,\o)} =
( \sum_{i=1}^\dimv \beta^{ki} \o_i) ( \sum_{j=1}^\dimv \overline{ \beta^j \o_j} ).
\ee
Since $| \sum_{j=1}^\dimv \overline{ \beta^j \o_j}| = 1$, this derivative is zero only if
the first factor is, so $\frac{\partial F}{\partial \mu_k} |_{(0,\o)} = 0$ implies 
\be
\label{eq511}
\sum_{i=1}^\dimv \beta^{ki} \o_i = 0.
\ee
To know $X$ has interior, we  just need to know that there is some point $(0,\o)$  at 
which the derivative of $F$ with respect to $\mu$ is non-zero; so if this is not true, we get that for
  all $\o$ in $  bT(\b) \cap U'$, \eqref{eq511} still holds.

 There are two possibilities:
 
 (i) the span of $\{ \o \in  U' \cap bT(\b) \}$ is all of $\C^d$, 
 in which case each $\b^{ki}$ must be $0$.

 (ii) $\M = \vee \{ \o \in  U' \cap bT(\b) \}$ has a non-zero vector $\eta$ in $\M^\perp$.
 In this case, as $\o$ ranges over $U' \cap bT(\b)$, a $(d-1)$-dimensional subset of $\T^d$,
 we must have that $\sum_{i=1}^d \o_i \bar \eta_i = 0$. 
Suppose three (or more) $\eta_i$ are non-zero, including say $\eta_1$.
Then in $\M$, $x_1$ is a linear function of the other variables; but as
we vary $\o_2, \dots, \o_d$, the sum $\sum_{i>1} \frac{\bar \eta_i}{\bar \eta_1}\o_i$ has interior,
and therefore cannot always have modulus $1$.
So Case (ii)  can happen if and only if exactly two of the components of $\eta$ are non-zero, and their
 ratio $- \bar \gamma$ is unimodular.  For notational simplicity, let us assume these are $\eta_1 $ and $\eta_2$.
Then $bT(\b) = \{ \o \in \T^d: \o_2 = \g \o_1 \}$,
$| \b^1 + \gamma \b^2| = 1$, and $\b^{k1} + \gamma \b^{k2} = 0$.

Now apply the same analysis to $\b_\nu$ for $\nu \in G'$. We have
\beq
\b_\nu(x) &=& \vp_{\b(\nu)} \circ \beta \circ \Phi_{-\nu}(x) \\
\b_\nu^i(x) &=&
\vp'_{\b(\nu)} (\b(\Phi_{-\nu}(x))) \b^i(\Phi_{-\nu}(x))\vp'_{-\nu_i}(x)\\
&=& \frac{1-|\b(\nu)|^2}{(1-\overline{\b(\nu)} \b(\Phi_{-\nu}(x))^2} \b^i(\Phi_{-\nu}(x)) \frac{1-|\nu_i|^2}{(1+\bar \nu_i x_i)^2}
\\
\b_\nu^{ki}(x) &=& \frac{2(1-|\b(\nu)|^2)}{(1-\overline{\b(\nu)} \b(\Phi_{-\nu}(x))^3}\b^k(\Phi_{-\nu}) \frac{1-|\nu_k|^2}{(1+\bar \nu_k x_k)^2}\b^i(\Phi_{-\nu}(x))  \frac{1-|\nu_i|^2}{(1+\bar \nu_i x_i)^2} \\
&&-  \frac{1-|\b(\nu)|^2}{(1-\overline{\b(\nu)} \b(\Phi_{-\nu}(x))^2} \b^{ki}(\Phi_{-\nu}(x))
\frac{1-|\nu_k|^2}{(1+\bar \nu_k x_k)^2} \frac{1-|\nu_i|^2}{(1+\bar \nu_i x_i)^2}\\
&& -2 \delta_{ki} \frac{1-|\b(\nu)|^2}{(1-\overline{\b(\nu)} \b(\Phi_{-\nu}(x))^2} \b^i(\Phi_{-\nu}(x))
 \frac{1-|\nu_i|^2}{(1+\bar \nu_i x_i)^3} \bar \nu_i\\
\b_\nu^{ki}(0) &=& \frac{2}{(1-|\b(\nu)|^2)^2} \b^k(\nu) (1-|\nu_k|^2)(1-|\nu_i|^2)\b^i(\nu) \\
&&- \frac{1}{1-|\b(\nu)|^2} \b^{ki}(\nu) (1-|\nu_k|^2)(1-|\nu_i|^2) \\
&&-2 \delta_{ki} 
\frac{1}{1-|\b(\nu)|^2} \b^{i}(\nu) (1-|\nu_i|^2) \bar \nu_i .
\eeq
If we are always in Case (i), then 
every $\b_\nu^{ki}(0) $ vanishes,
 so for every $i \neq k$ we get
\be
\label{eq51.1}
\b^{ki}(\nu) = 2 \frac{\b^i(\nu) \b^k(\nu)}{1-|\b(\nu)|^2} .
\ee
The left-hand side of \eqref{eq51.1} is holomorphic in $\nu$, the right-hand side is not.
Therefore both sides must be constant, and equal to $\b^{ki}(0) =0$. This means $\b^i(\nu) $ is identically zero for all except possibly one value of $i$. But if $\b$ depends on only one variable, $x_1$ say, then 
either $\b^1(0)$ has modulus less than $1$, in which case $T(\b) = \T^\dimv$ and $b T(\b)$ is empty, or the modulus is greater than or equal to $1$, in which case $T(\b)$ is empty.

Suppose instead we are in Case (ii). Then all the partial derivatives of $\b$ except the first two vanish, so
$\b$ is just a function of $x_1, x_2$. Moreover, we have
\[
\gamma(\nu) = - \frac{\b^1 (\nu)}{\b^2(\nu)} \frac{1-|\nu_1|^2}
{1-|\nu_2|^2}
\]
has modulus $1$. But if we let $|\nu_1| = |\nu_2| = r$ be small, this says that 
the holomorphic function $\frac{\b^1 (\nu)}{\b^2(\nu)}$ has modulus $1$ on $r\T^2$.
By the maximum principle, it must be constant on some neighborhood of $0$.
But now varying the modulus of $\nu^1$ and $\nu^2$ leads to a contradiction.
So Case (ii) cannot occur for generic perturbations.

Therefore we conclude that $F$ is non-constant, and hence $X$ has interior.

By Lemma \ref{lemd2}, 
a.e point of $X$ is a point of continuity of some $\beta_\mu$,
and by \eqref{eq50} the values  are unimodular, so $\beta$ is unimodular a.e. on $X$.
\oec

We cannot apply Theorem \ref{thmc2} directly to $\beta$ yet, because we do not know that $\beta$ is holomorphic
on $Y$. Lemma \ref{lemh3} will allow us to circumvent the difficulty.
From Tworzewski's Theorem \ref{thmc8}, we know that $\b$ is a Nash function, so there is some irreducible  polynomial $P \in \C[x_1, \dots, x_\dimv][x_{\dimv +1}] $ such that $P(x, \b(x)) = 0$. 
Let $W$ be the interior of $X $.
 The function $\beta$ is one of the branches, which without loss of generality we can assume is $f_1$.
The product $\beta f_2 \cdots f_m$ extends to a rational function $g$,
which is $(-1)^m P_0/P_m$ if we write $P(x,z)  = \sum_{j=0}^m  P_j(x) z^j$.
By Lemma \ref{lemh3}, this function has modulus one on a non-empty
open subset of $\T^\dimv$, so since it is rational, we conclude that $g$ has modulus one a.e. on $\T^\dimv$.

Now we repeat the argument for $\Phi_b (V)$, where $b = (\mu, \alpha(\mu))$ and  $\mu$ ranges over $G'$.  We conclude that for all $\mu$ in $G'$, 
\[
\vp_{\beta(\mu)} (\beta \circ \Phi_\mu)  \vp_{\beta(\mu)} (f_2 \circ \Phi_\mu) \cdots \vp_{\beta(\mu)} (f_m \circ \Phi_\mu)
\]
is rational and unimodular a.e. on $\T^\dimv$, and, since $\Phi_\mu$ is a rational diffeomorphism of  $\T^\dimv$, 
the function 
\be
\label{eq52}
\vp_{\beta(\mu)} (\beta ) \vp_{\beta(\mu)} (f_2 )
\cdots \vp_{\beta(\mu)} (f_m)
\ee
is also rational and unimodular a.e. on $\T^\dimv$.
Let $c = \beta(\mu)$; then \eqref{eq52} becomes
\be
\label{eq53}
g_c \ := \ \vp_c (\beta) \vp_c(f_2) \cdots \vp_c (f_m) 
\ee
is rational and unimodular a.e. on $\T^\dimv$ for $c$ in a non-empty open subset of $\D$, and hence by the identity principle
for all $c$ in $\D$.

Let $E = \{ \tau \in \overline{T(\a)} : \b(\taudotlambda) {\rm \ is\ a\ finite\ Blaschke\ product} \}$.
Write $E = \cup_{k=1}^\i E_k$ where $E_k$ is the set of $\tau$ for
which $\b(\taudotlambda)$ is of degree $k$. If
the closure of some $E_k$ has interior, then by 
Rudin's theorem \ref{thmc2}, we conclude
that $\b$ is a rational function of degree $k$.
Suppose then that this does not happen, so by 
 the Baire category theorem, the complement of $E$ is dense.

We can find $\tau$ arbitrarily close to $\om$, lying in
$T(\a)$, and such that the function $\beta(\taudotlambda)$, which we know from Lemma \ref{lemd2} is
in the disk algebra, is not a finite Blaschke product. Therefore there is some point $\zeta \in \T$
such  that $\beta( \tau \cdot \zeta) = c \in \D$. Since $g_c$ is a rational function that is unimodular on 
$\T^\dimv$ whenever its numerator and denominator do not both vanish, we get from
\eqref{eq53} that there must be some $j$ between $2$ and $m$ so that $\beta( \tau \cdot \zeta) \overline{f_j( \tau \cdot \zeta)} = 1$. We can repeat this argument for points that lie in a neighborhood
of $\tau \cdot \zeta$ in $\T^\dimv$. 
We conclude that the function
\be
\label{eq54}
(1 - \beta \overline{f_2}) \cdots (1 - \beta \overline{f_m}) \= 0
\ee
on an open neighborhood of $\tau \cdot \zeta$ in $\T^\dimv$. Since the left-hand side of
\eqref{eq54} is real analytic, it must be zero on any connected open set to which all the functions 
extend analytically. In particular, \eqref{eq54} holds on $Z$,
and therefore by continuity it holds a.e. on $X$.
Since $\beta$ is unimodular a.e. on $X$, we therefore have that
\be
\label{eq55}
\notag
( \beta -f_2) \cdots (\beta - f_m) \= 0
\ee
a.e. on $X$.
Therefore the discriminant of $P$ vanishes on $X$, and therefore is identically $0$, which contradicts the fact that $P$ is irreducible. 
\ep

\begin{lemma}
\label{lemh4}
Suppose there exists $\vare > 0$ so that, for all $a \in V \cap B(0,\vare)$, if we write as in \eqref{eq20}
\[
\tilde \a_j(y_1, \dots, y_\dimv) \ = \vp_{a_j} \circ \a_j(\vp_{-a_1}(y_1), \dots, \vp_{-a_\dimv}(y_\dimv) ) ,
\]
then 
there is an open set $U$ in $\T^d$ such that $U \cap b\overline{T( \tilde \al_{\dimv +1}) } $ is non-empty and contained in $\overline{T(\tilde \alpha)}$.

 Then $\al_{\dimv +1}$ is a rational function, and has modulus one a.e. on $\T^\dimv$.
\end{lemma}

\bp 
This follows from Lemmas \ref{lemh3.5} and \ref{lemh1}.
\ep

The set $ \{ \g \in \overline{\M}: \g(0) = 0 \}$ is the set of rational functions $\kappa_b$ as in \eqref{eq4.1} 
where $b_3 \neq 0$ and the maximum of $\{ |b_1|, |b_2|, |b_3| \}$ is less than or equal to the sum of the other two. 

\begin{lemma}
\label{lemh6}
Suppose $\g$ is a rational function of $\dimv$ variables, that is unimodular on $\T^\dimv$ a.e., and such that $\g(0) = 0$. Suppose there is an open set $U$ in $\T^\dimv$ such that $U \cap b\overline{ T(\gamma)} \neq \emptyset$ and so that for all $\o \in U \cap b \overline{T(\gamma)}$ there exists $\eta \in \T$ so that $\g(\o \cdot \l) = \eta \l$ for all $\l \in \D$. Suppose moreover that this is true for all small perturbations of $\gamma$. Then $\dimv = 2$ and $\gamma \in \overline{\M}$.
\end{lemma}
\bp
We must have $\dimv \geq 2$, since if $\dimv=1$, the assumption $\g(\om \l) = \eta \l$ violates Assumption \ref{ass1}.
Let  us write $\om = (\om', \om_\dimv)$ with $\om' \in \T^{\dimv-1}$, and $ \C^\dimv \ni x = (x', x_\dimv)$. For $c \in \D$, let
\[
\gc(x',x_\dimv)
\=
\vp_{\g(0,c)} \circ \g (x', \vp_{-c}(x_\dimv)) .
\]
Let us fix some $\o'$ so that $(\om', \om_\dimv) \in U \cap  b\overline{ T(\gamma)} $.
The hypotheses imply that for $c$ close to $0$, there are $\om_c, \eta_c \in \T$ so that, if $\om \in U \cap b T(\gamma)$, then
\be
\label{eq76}
\gc(\om' \cdot \l, \om_c \l) \= \eta_c \l \qquad \forall \ \l \in \D .
\ee

Assuming the claim holds, we have that 
\be
\label{eq77}
\g (\om' \cdot \l, \vp_{-c} (\om_c \l)) \=
\vp_{-\g(0,c)} (\eta_c \l), 
\ee
and the right-hand side of \eqref{eq77} is
rational of degree at most  $1$ in $\l$.
  Let $c$ tend to some point $\tau \in \T$.
The left-hand side of \eqref{eq77} converges to $\g(\om' \cdot \l, \tau)$, and the right-hand side converges to either a  linear fractional transformation  of $\l$ or a point, so is either degree $1$ or $0$ in $\l$.

If $\dimv =2$, this means $\g$ is of degree at most one in $x_1$. If $\dimv > 2$, perturbing $\om'$ we can conclude that
$\g$ is of total degree at most one with respect to $x_1, \dots, x_{\dimv-1}$. 
Repeating the argument with $x_\dimv$ and $x_1$ interchanged, we conclude that $\g$ is of degree at most one with respect to each variable separately, and of total degree at most $2$, with at most one quadratic term. After permuting coordinates, this means that the numerator of $\g$ is of the form
\[ A_1 x_1 + \dots + A_\dimv x_\dimv + x_1 x_2 .
\]
By Theorem \ref{thmc4} the denominator must be the reflection of this, and non-zero at $0$, which requires $\dimv = 2$.
Finally the requirement at $c=0$ that $\g(\om \cdot \l) = \eta \l$ means $|A_1 \om_1 + A_2 \om_2| = 1$,
so $\g \in \overline{\M}$.

To prove the claim, note that, by Theorem \ref{thmc4}, there is some polynomial $P_c$  so that
\[
\gc ( x) \= x^I \frac{P_c(x)}{\revP_c(x)} .
\]
The coefficients of $P_c$  can be chosen analytically in $c$ and $\bar c$, since $\gc$ is analytic in $c,\bar c$ and
$P_c$ is determined, up to a constant, by the zero set of $\gc$.
Then \eqref{eq76} becomes
\be
\label{eq78}
P_c(\om' \cdot \l, \om_c \l) \= \eta_c \l \revP_c(\om' \cdot \l, \om_c \l).
\ee
Comparing coefficients of different powers of $\l$, \eqref{eq78} gives a system of equations relating $\om_c$ and $\eta_c$.
The equation that comes from $\l$ to the power of $1$ is of the form
\be
\label{eq79}
\notag
E_c \om_c + F_c \= \eta_c ,
\ee
where $E_c$ and $F_c$ are analytic in $c,\bar c$. 
The fact that $\om_c, \eta_c$ are unimodular means that for $c$ close to zero, there is a unimodular solution $\om_c$ to the quadratic equation
\be
\label{eq80}
(E_c \bar F_c) \om_c^2 + (|E_c|^2 + |F_c|^2 - 1) \om_c + (\bar E_c F_c) \= 0 .
\ee
There are two possibilities:

(i) One of the equations from \eqref{eq78} is an independent condition on $\om_c$.

(ii) All the equations from comparing powers of $\l$ in \eqref{eq78} reduce to \eqref{eq80}.

If possibility (i) holds, then by the Euclidean algorithm we can get a linear equation for $\om_c$ with real analytic coefficients. Since it has a unimodular solution for $c$ near $0$, it will have a unimodular solution for all $c$.

If possibility (ii) holds, note that the product of the two roots of \eqref{eq80} is real analytic in $c$ and unimodular.
Therefore one of them must be in $\overline{\D}$, as required.
%
\ep

{\bf Lemma \ref{lemd3}.}
{\em
Suppose $\om \in b \overline{T(\alpha_j)} \cap b {T(\a)} $ and moreover $\om \in \cap_{\dimv< \ell \le N, \ \ell \neq j} T(\alpha_\ell)$.
Then $\alpha_j \in \M$.
}

\bp
Without loss of generality, $j = \dimv+1$. Perturbing to $\tilde \a$, we will perturb $b T(\a_{\dimv +1})$ and for small
perturbations this will stay in $\cap_{j=\dimv +2}^{\dimpol} T(\tilde \a_j)$. So by Lemma \ref{lemh4}, we conclude that
$\a_{\dimv +1}$ is rational and unimodular a.e. on $\T^\dimv$, satisfying the hypotheses in the first sentence of Lemma \ref{lemh6}. 
By Lemma \ref{lemh1} the next two sentences are also satisfied. So $\dimv = 2$ and $\a_{\dimv +1}$ is in $\overline{\M}$.
This means $\a_{\dimv +1} = [ \tau; A, B ]$ where $1, |A|, |B|$ satisfy the triangle inequality, but perhaps not strictly.

If $|A| + |B| = 1$, then $\overline{T(\a_{\dimv +1})} = \T^\dimv$, so it has no boundary.
If $|A| = 1 + |B|$ (or vice versa), then $T(\a_{\dimv +1})$ is empty.
So we must have that $(1, |A|, |B|)$ are triangular, and $\a_{\dimv +1}$ is in $\M$.
\ep


\section{Proof of Lemma \ref{lemd4} }
\label{seci}

\begin{lemma}
\label{lemh5}
Suppose $\a, \b$ are in $\M$, and for all small perturbations we have
$T(\tilde \a) = T(\tilde \b)$.
Then $\a$ is a unimodular multiple of $\b$.
\end{lemma}

\bp
Assume $\a = [\tau_1;A,B]$ and $\b = [\tau_2;C,D].$  The end points of $T(\a)$ are the solutions of $|A \om + B|=1$,
which are the  roots of the equation
\be
\label{eq73}
\notag
A \bar B \om^2 + [|A|^2 + |B|^2 - 1] \om + \bar A B \= 0.
\ee
If $T(\a)$ and $T(\b)$ have the same end-points, then we must have
\begin{eqnarray}
\label{eq74}
\frac{|A|^2 + |B|^2 - 1}{A \bar B} &\=& 
\frac{|C|^2 + |D|^2 - 1}{C \bar D} \\
\label{eq75}
\frac{ \bar A B}{A \bar B} &=& 
\frac{ \bar C D}{C \bar D}
\end{eqnarray}
We can re-express \eqref{eq75} as saying that $A/B$ is a real multiple of $C/D$.

Now let $a = (a_1, a_2, \a(a_1,a_2), \b(a_1,a_2))$, and let 
$\tilde \a = \vp_{a_3} \circ \a(\vp_{-a_1}, \vp_{-a_2})$ and 
$\tilde \b = \vp_{a_4} \circ \b(\vp_{-a_1}, \vp_{-a_2})$.
Then by \eqref{eq21}, $A/B$ becomes $\frac{1-|a_1|^2}{1-|a_2|^2} \frac{\a^1}{\a^2}$,
where $\a^i$ is the $i^{\rm th}$ partial derivative at $a$.
After this change of variables, condition \eqref{eq75} becomes
$\a^1/\a^2$ is a real multiple of $\b^1/\b^2$.

Now let us specialize to the case that $\b = \kappa = [1;-1,-1]$, which we can do by Theorem \ref{thmb2}.
Then \eqref{eq75} says $A/B$ is real, and \eqref{eq74} implies
$|A|^2 + |B|^2 = 1 \pm |AB|$.
We have
\[
\frac{\a^1(x,y)}{\a^2(x,y)} \= \frac{A+ ( |A|^2 -|B|^2 +1) y + \bar A y^2}{B + (|B|^2 - |A|^2 +1) x + \bar B x^2}\]
For all small $x,y$, this has to be a real multiple $r$ of the corresponding expression for $\kappa$, which is
\[
\frac{\kappa^1(x,y)}{\kappa^2(x,y)} \= \frac{1 -y + y^2}{1 - x + x^2}.
\]
Cross-multiplying, we get
\[
[A+ ( |A|^2 -|B|^2 +1) y + \bar A y^2][1 - x + x^2] \= r \, [B + (|B|^2 - |A|^2 +1) x + \bar B x^2] [ 1 -y + y^2] .
\]

Comparing  coefficients, we see that $A = rB = \bar A$, so $A$ and $B$ are real, and
$A^2 - B^2 + 1 = -rB = r(B^2 - A^2+1)$. If $r =1$, we get $A = -1 = B$.
Otherwise, we get from the equality $(r^2 -1) B^2 + 1 = r((1-r^2 ) B^2 +1)$ that $B^2 = \frac{1}{(1+r)^2}$, 
and using this in the equality $(r^2 -1) B^2 + 1 = -rB$ we get $B = \frac{-2}{1+r}$, a contradiction.
\ep

{\bf Lemma \ref {lemd4}.}
{\em
Let $\dimv=2$.
Let $\ell \ge \dimv+2$.
Suppose there exists $\vare > 0$ so that, for all $a \in V \cap B(0,\vare)$, if we write as in \eqref{eq20}
\[
\tilde \a_j(y_1, \dots, y_\dimv) \ = \vp_{a_j} \circ \a_j(\vp_{-a_1}(y_1), \dots, \vp_{-a_\dimv}(y_\dimv) ) ,
\]
then $ b \overline{T( \tilde \al_{\dimv+1}) }\cap \dots \cap b \overline{T(\tilde \al_{\ell}) }\cap T(\tilde \al_{\ell+1}) \cap \dots \cap T(\tilde \al_N)$ is non-empty.
 Then $\al_{\dimv+1}, \dots , \al_{\ell}$ are in $\M$,
 and they form an edge family. 
}

\bp
For each $j$ between $\dimv +1$ and $\ell$ we know $\hatt(\a_j)$ is a proper open arc in $\T$.
We can apply
Lemma \ref{lemh4} to conclude that each $\a_j$ is rational and unimodular on $\T^2$.
By Lemmas \ref{lemh1} and  \ref{lemh6}, they are all in $\overline{\M}$.
 Since $b \overline{T(\a_j) }\neq \emptyset$,
they must all lie in $\M$.

Let $\o \in b\overline{T(\a_{\dimv +1}) }\cap \dots \cap b\overline{T(\a_\ell) }\cap T(\a_{\ell+1}) \cap \dots \cap T(\a_N)$.
There are unimodular $\eta_j$ so that
\be
\label{eq81}
\notag
\a_j(\om_1 \l, \om_2 \l) \= \eta_j \l,\qquad \dimv +1 \le j \le \ell .
\ee
For any $ \dimv +1 \le i \neq j \le \ell$,  define $\g_i$ by 
$\gamma_{i} (x_1, x_j) = \a_{i}(x_1, \a_{j}^{-1} (x_1,\underline{ x_j}))$.
Then
\beq
\g_i(\om_1 \l, \eta_j \l) &\=& \a_i(\om_1 \l, \a_j^{-1}(\om_1 \l, \eta_j \l)) \\
&=& \a_i(\om_1 \l, \om_2 \l) \\
&=&
\eta_i \l .
\eeq
By hypothesis, this also holds for small perturbations, since $\tilde \g_i(x_1, x_j)
= \tilde \a_i( x_1, \tilde  \a_{j}^{-1} (x_1,\underline{ x_j}))$.
As $\g$ is unimodular on $\T^2$ a.e.,
 by Lemma \ref{lemh6}, each $\g_i \in \overline{\M}$.
Similarly, $\a_i( \a_{j}^{-1} (\underline{ x_j}, x_2),x_2) \in \overline{\M}$,
 so by 
Theorem \ref{thmb3}, either the set $\{ \a_{\dimv +1}, \dots, \a_\ell \}$ form an edge family,
or there is some function, $  \a_{\dimv +1}$ say, so that
for every $ \dimv +2 \le n \le \ell$, we have
$
 D_{\a_{\dimv+1}}  \subseteq D_{\a_n} 
$.
But the latter case would violate Assumption \ref{ass2}, as we could just throw those functions out
and get a lower dimensional \car set.
\ep

\section{Proof of  Lemma \ref{lemd5}}
\label{secj}

{\bf Lemma \ref{lemd5}.} {\em
Let $\dimv \geq 3$.
Let $\ell \ge \dimv+2$.
Suppose there exists $\vare > 0$ so that, for all $a \in V \cap B(0,\vare)$, if we write as in \eqref{eq20}
\[
\tilde \a_j(y_1, \dots, y_\dimv) \ = \vp_{a_j} \circ \a_j(\vp_{-a_1}(y_1), \dots, \vp_{-a_\dimv}(y_\dimv) ) ,
\]
then $ b {T( \tilde \al_{\dimv+1}) }, \dots , b {T(\tilde \al_{\ell}) }$ all share
$b_{\dimv -1}T(\a)$ at some point that lies in $ T(\tilde \al_{\ell+1} )\cap \dots \cap T(\tilde \al_N)$.
 Then $\al_{\dimv+1}, \dots , \al_{\ell}$ are in $\M$, and they all depend on the same pair of variables.
Moreover, they form an edge family.
 }
 
 \bp
Choose two functions $\a_j, \a_{j'}$ for
$j,j'$ between $\dimv +1$ and $\ell$. By Lemma \ref{lemh4} they are both rational, and unimodular a.e. on
$\T^\dimv$, and by Lemmas \ref{lemh1} and  \ref{lemh6} we can conclude that they are both actually in $\M$,
just as we argued in the proof of Lemma \ref{lemd4}.

If $\a_j$ and $\a_{j'}$ do not depend on the same pair of variables, then $bT(\a_j) \cap bT(\a_{j'})$ would be $(\dimv-2)$-dimensional, so it would be impossible for both $bT(\a_j)$ and $bT(\a_{j'})$ to share $b_{\dimv-1}T(\a)$ at a point
$\om \in bT(\a_j) \cap bT(\a_{j'})$. 

So we know all the functions $\a_j$ depend on the same pair of variables.
Now we are in the same situation as in Section \ref{seci}.
The same arguments that proved Lemma \ref{lemd4} now show that
$\{ \al_{\dimv +1}, \dots, \a_\ell \}$ form an edge family.
\ep


\section{Proof of Lemma \ref{lemd7}}
\label{secka}

\begin{lemma}
\label{lemka2}
Let $\a_j = [1; A, B]$, let $\o \in bT(\a_j (x_1, x_2))$, and let $\o_j = A \o_1 + B \o_2$.
Then there is a direction one can perturb $a_1$, keeping $a_2 = 0 = \dots = a_\dimv$, so that, with $\o_2$ staying fixed,
$\o_1$ and  $\o_j$ move.
\end{lemma}
 \bp
Since $T(\a_j)$ is $\T$ invariant, we can assume $\o_2 = 1$.
We have $ A \o_1 + B = \o_j$.
Now apply a small perturbation $\Phi_a$, where only $a_1 \neq 0$.
Then $bT(\tilde \a_j)$ will contain the point $(\tilde \o_1, 1)$ where
$| \tilde A \tilde \o_1 + \tilde B | = 1$, and  $\o_j$ will become $\tilde \o_j = \tilde A \tilde \o_1 + \tilde B \in \T$.
We want to show that $\tilde \o_1$ must move.

Suppose it does not. Then $\tilde \a_j( \o_1 \l, \l) = \tilde \o_j \l$ for all $\l$ in $\D$.
This means
\be
\label{eqk1}
\a_j(\vp_{-a_1}(\o_1 \l), \l) = \vp_{-\a_j(a_1,0)}(\tilde \o_j \l) .
\ee
Write $\tilde a_j$ for $\a_j(a_1,0)$.
Expanding, \eqref{eqk1} becomes
\be
\label{eqk2}
\frac{(A + \l)( \o_1 \l  + a_1 ) +B\l(1 + \bar a_1  \o_1 \l)}{(\bar A \l + 1)(1+ \bar a_1  \o_1 \l)+ \bar B ( \o_1 \l + a_1 )}
\=
\frac{ \tilde \o_j \l + \tilde a_j}{1 + \overline{\tilde a_j} \tilde \o_j \l } .
\ee
The numerator $p(\l)$ and denominator $q(\l)$ on the left-hand side of \eqref{eqk2} must therefore share a root for each $a_1$.
It is well-known that two quadratics $E_1 \l^2 + F_1 \l + G_1$ and  $E_2 \l^2 + F_2 \l + G_2$
share a root if and only if their resultant vanishes, which means
\be
\label{eqk3}
(G_1E_2-G_2E_1)^2 = (E_1 F_2 - E_2F_1)(F_1 G_2 - F_2 G_1) 
\ee
In our case
\[
p(\l) = (1+\bar a_1 B)\o_1 \l^2 + (A \o_1 + B + a_1) \l + A a_1,
\]
and
$ E_2 = \o_1 \bar G_1, F_2 = \o_1 \bar F_1, G_2 = \o_1 \bar E_1$, so to check \eqref{eqk3} we need to check
\beq
\lefteqn{
[ |A|^2 |a_1|^2 - |1+\bar a_1 B|^2 ]^2 } \\
&& \stackrel{?}{=} 
\left| (1+ \bar a_1 B) \o_1 ( \bar A \bar \o_1 + \bar B + \bar a_1) + \bar A \bar a_1 (A \o_1 + B  +a_1) \right|^2
\eeq
To see this cannot hold identically, look at the fourth order terms in $a_1$. The left-hand side has no term in $a_1^2 |a_1|^2$, and the coefficient on the  right-hand side  
  of $a_1^2 |a_1|^2$ is $ \bar A \bar B \bar \o_1$, which cannot be zero.

Swapping the roles of $j$ and $1$, we get that $\o_j$ must also vary as $a_1$ is perturbed with $a_2$ constant.
\ep 

{\bf Lemma \ref{lemd7}.}
{\em
Let $\a_\ell $ and $\a_{\ell'}$ be in $\BB$.
Then either $\a_\ell$ and $\a_{\ell'}$ depend on the same pair of variables, or on disjoint pairs of variables.
}

 Suppose  there are functions $\a_j(x_1, x_2)$ and $\a_{j'}(x_1, x_3)$ that are both in $\BB$.
Choose $m$ as in Lemma \ref{lemk1} with $\a_m$ in $\M$ and a function of $x_1$ and $x_2$,
and so that $\{x_2, x_3, \dots, x_\dimv, x_m \}$ is a good parametrization, with corresponding functions $\b$.

Consider \[
\b_{j'}(x_2, \dots, x_m) \= \a_{j'} (\b_1(x_2, x_m), x_3) .
\]
This is a rational function of 3 variables that is unimodular when $|x_2|=|x_3| =|x_m| = 1$.
Since we have a good $\{x_2, x_3, \dots, x_\dimv, x_m \}$ parametrization, we know $T(\beta)$ is non-empty.
We claim $\overline{T(\b_{j'})} \neq \T^\dimv$.
Indeed, since $\a_{j'}$ and $\b_1$ are in $\M$, we know $|\a_{j'}^1| + |\a_{j'}^3| > 1$ and $|\b_1^2| + |\b_1^m| > 1$. Therefore
\[
|\b_{j'}^2| + |\b_{j'}^3| + |\b_{j'}^m | = |\a_{j'}^1|(|\b_1^2| + |\b_1^m|) +  |\a_{j'}^3| > 1 ,
\]
and so $\overline{T(\b_{j'})} \neq \T^\dimv$.

Suppose we know that $b T(\b_{j'})$ shares $b_{\dimv -1} T(\b)$.
Then we can use Lemmas \ref{lemd2} and \ref{lemh1} to show that
 the hypotheses of Lemma \ref{lemh6} are satisfied for $\b_{j'}$ and small perturbations.  The conclusion of the lemma would say that $\b_{j'} $ is a function of two variables, a contradiction. 

We can't quite do this, but we can after automorphisms.
Our argument will just depend on $T(\a)$, so we can assume for notational convenience 
  that for each pair $(x_i, x_j)$ 
there is at most one function 
 $\a_{ij}(x_i, x_j) = \a_{ji}(x_i,x_j)$,
where $\a_{ij} = [1; A_{ij}, B_{ij}]$.
Let \[
\Delta := T_0(V) \cap \overline{\D}^\dimpol  \= \{ y \in \overline{\D}^\dimpol : A_{ij} y_i + B_{ij} y_j = y_{ij} \} .
\]
(Recall \eqref{eqtt}.) 
For the rest of the proof of Lemma \ref{lemd7},  any number denoted $\o_i$ is unimodular,
any number denoted $z_i$ is in $\D$, and any number $\z_i \in \overline{\D}$, for any subscript.
We shall use $*$ to mean other entries are from $\D$, and $+$ to mean they are from $\overline{\D}$.
Any point in the set $\Delta$, the intersection of the $\dimv$-dimensional tangent space at $0$ with ${\overline \D}^\dimpol$, is uniquely determined by its first $\dimv$ coordinates.
For a point $(\o_1, \dots, \o_\dimv, \z_{ij})$ in $ \Delta \cap (\T^{\dimv} \times {\overline \D}^{\dimpol - \dimv})$, the coordinate
 $\z_{ij}$ depends only on $\o_i$ and $\o_j$.

Since $\a_{12}$ is bounding\footnote{Note that $\a_{12}$ and $\o_{12}$ mean $\a_{ij}$ and $\o_{ij}$ for $i=1$ and $j=2$.}, we have a point
$\o^{(0)}= (\o_1^{(0)}, \dots, \o_\dimv^{(0)}, \o_{12}^{(0)}, *) $ in $\De$.
Move $\o_3$ in $\T$, keeping the other $\dimv-1$ coordinates between $1$ and $\dimv$ constant, until some other coordinate becomes unimodular.

Case: the first such coordinate to become unimodular is $\o_{13}$ or $\o_{23}$.
We can assume that only one new coordinate becomes unimodular.

Indeed, by Lemma \ref{lemka2}, by applying a small M\"obius perturbation to $V$, we can assume that we don't have 
both $\z_{13}$ and $\z_{23}$ simultaneously becoming unimodular, since a small perturbation in  $a_1$ will 
affect $\o_{13}$ but not $\o_{23}$. 
Likewise, if 
  some  other point $\z_{3j}$ became unimodular at the same time, then we could perturb $\o_j$ to make the new modulus less than $1$. Since all the other terms $z_{kj}$ for $k > 3$ are in $\D$, they will stay there for small perturbations of $\o_j$.

So we can assume we have a point $\o' = (\o_1', \dots, \o_d', \o_{12}', \o_{13}', *) $ in $\De$.
Moreover, by construction we have $\pi_{1, \dots, d} \o' \in \overline{T(\a)}$.
Change variables to $\beta(x_2, \dots, x_d, x_{12})$. 

We would like to show $\pi_{2,\dots,d, 12}\o' \in \overline{T(\b)}$.
From  Lemma \ref{lemk2} (with $m = 12$) the point $\o^{(0)}$ is in $bT(\beta)$.
We have
\[
\b_{13}(x_2, x_3, x_{12}) = \a_{13}(\b_1(x_2, x_{12}), x_3 ),
\]
and $\b_{13}^3 = \a_{13}^3 \neq 0$.
Since $\z_1 = C_2 \o_2 + D_{12} \o_{12}$, where $\b_1 = [1; C_2, D_{12}]$,
there is one direction around the circle such that perturbing $\o_2$ in that direction causes  $|\z_1|$ to decrease.
So if we rotate $\o_{12}$ a small amount in that direction, and then move $\o_3$ to stay in $bT(\b_{13})$, 
we get a point $\o'' = (z_1, \o_2, \dots , \o_d, \o_{12}, \o_{13}, *) \in bT(\b) \cap b\overline{T(\b_{13})}$.
We now have two degrees of freedom: if we perturb $\o_2$ and $\o_{12}$ in $\o''$ independently, we can find an appropriate perturbation of $\o_3$ so that $\o_{13}$ remains unimodular, and all the other terms that were in $\D$ remain there.
Now we can apply Lemma \ref{lemh6} to conclude that $\b_{13}(x_2, x_3, x_{12})$ is actually a function of two variables, and this gives us our desired contradiction.

Case: The first coordinate to become unimodular when we start with $\o^{(0)}$ and rotate $\o_3$ is neither $\o_{13}$ nor $\o_{23}$. Then it is some other coordinate, $\o_{34}$ say.
Again, by Lemma \ref{lemka2} we can assume that we don't hit two unimodular numbers simultaneously.
Continue rotating $\o_3$, but now also rotate $\o_4$ to keep $\o_{34}$ (which only depends on $\o_3$ and $\o_4$) unimodular. We continue until we get another unimodular hit,
which we can assume is either $\o_{35}$ or $\o_{45}$. Either way, we now have two unimodular numbers with one coordinate in common, say $\o_{34}$ and $\o_{35}$.
Now repeat the earlier argument, but replace $x_1$ with $x_3$ as the variable that gets removed, to be replaced by $x_{34}$.
\ep

\section{Proof of Lemma \ref{lemd6}}
\label{seck}

\begin{lemma}
\label{lemk0}
Let $\a, \b \in \M$, and let $W = \{ (x,y,\a(x,y),\b(x,y) ) \} \cap \D^4$.
Suppose $T(\a)$ and $ T(\b)$ overlap, neither is contained in the other, and the four end-points of  $\hatt(\a)$ and $\hatt(\b)$ are distinct. 
Then $W$ is not a \car set at $0$, i.e. for all $\vare > 0$ there exists $P \in W \cap B(0,\vare)$
so that $d_W(P,0) > d_{\D^4} (P,0)$.
\end{lemma}
\bp
We can assume $\a = [1; A,B]$ and $\b = [1;C,D]$, and that $\hatt(\a) = (\o_1, \o_3)$
and $\hatt(\b) = (\o_2, \o_4)$ with $\arg(\o_1) < \arg(\o_2) < \arg(\o_3) < \arg(\o_4)$. 
The circle $C_1 = \{ \zeta : |A + B \zeta = 1| \}$ intersects the unit circle at $\o_1$ and $\o_3$.
The circle $C_2 = \{ \zeta : |C + D \zeta = 1| \}$ intersects the unit circle at $\o_2$ and $\o_4$.
So there must be two points where $C_1$ and $C_2$ intersect, one in $\D$ and one in $\C \setminus \overline \D$. Let $\zeta$ be the point in $\D$ where
\[
|A + B \z| = 1 = |C+D\z| .
\]
Let $A + B \z = \o_\a$ and $C + D \z = \o_\b$. Define
\[
f(x_1, x_2, x_3, x_4) = \frac{1}{3} (x_1 + \bar \o_\a x_3 + \bar \o_\b x_4 ) .
\]
For $t$ very close to $0$, notice that
\beq
f(t, \zeta t, \a(t,\zeta t), \b(t, \zeta t) )
&=& f(t, \zeta t, \o_\a t, \o_\b t) + O(t^2)
\\
&=& t + O(t^2) .
\eeq
So if we can show that $F(W) \subsetneq \D$, then by Lemma \ref{lemg1} we are done.

Suppose $\overline{F(W)} \supseteq \T$. Then for every $x_1 \in \T$ there is some $x_2$ so that
\[
\a(x_1, x_2) = \o_\a x_1, \quad \b(x_1, x_2) = \o_\b x_1 .
\]
So
\[
\a^{-1} (x_1, \underline{\o_\a x_1}) = \b^{-1}(x_1, \underline{\o_\b x_1}) .
\]
We have
\[
\a^{-1} (x_1, \underline{\o_\a x_1}) = \left(\frac{\bar B \o_\a}{B}\right) \,  x_1 \, m_{\frac{B \bar \o_\a}{\bar B} \zeta} (x_1) := g(x_1) ,
\]
a Blaschke product of degree 2.

Now we swap the roles of $x_1$ and $x_2$, and find a point $\zeta' \in \D$ where $|A \zeta' + B| = 1 = |C \zeta' + D |$. We get unimodular constants $\eta_\a, \eta_b$ and a Blaschke product $h$ of degree $2$ that vanishes at $0$ so that $\partial W$ contains
\[
\Gamma_2 := \ \{ (h(x_2), x_2, \eta_\a x_2, \eta_\b x_2) : x_2 \in \T \} 
\]
as well as
\[
\Gamma_1 := \ \{ (x_1, g(x_1), \o_\a x_1, \o_\b x_1) x_1 \in \T \} .
\]
For convenience, write $g(x) = x b(x)$ and $h(x) = x c(x)$ where $b,c$ are of degree $1$.
The intersection $\Gamma_1 \cap \Gamma_2$ occurs when $x = h(g(x))$, which is where
\be
\label{eqko1}
1 = b(x) c(g(x) ).
\ee
The right-hand side of \eqref{eqko1} is a Blaschke product of degree $3$, so there are three distinct solutions, $\xi_1, \xi_2, \xi_3$. Let $ \upsilon_j = g(\xi_j) = \xi_j b(\xi_j)$.
Since $b$ is degree $1$, we get the three values
\be
\label{eqko2}
b(\xi_j) = \frac{\upsilon_j}{\xi_j} 
\ee
are distinct.
We have
\begin{eqnarray}
\label{eqko3}
\a(\xi_j, \upsilon_j) &= \o_\a \xi_j &= \eta_\a \upsilon_j \\
\b(\xi_j, \upsilon_j) &= \o_\b \xi_j &= \eta_\b \upsilon_j .
\label{eqko4}
\end{eqnarray}
Equation \eqref{eqko3} says
\be
\label{eqko5}
 \frac{\upsilon_j}{\xi_j}  = \omega_\a \bar \eta_\a,
 \ee
  and this 
 holds unless $(\xi_j, \upsilon_j)$ is a singularity of $\alpha$, which happens
if $-1 = \bar B \xi_j + \bar A \upsilon_j$, i.e. it fails only
when
\be
\label{eqko6}
A + B \upsilon_j \bar \xi_j = - \bar \upsilon_j .
\ee
Since the three values in \eqref{eqko2} are distinct, only one can equal $\omega_\a \bar \eta_\a$,
so two of them must be in $b \hatt(\a)$ and satisfy \eqref{eqko6}.

Similarly, two of them must be in $b \hatt(\b)$, which means
\be
\label{eqko7}
C + D \upsilon_j \bar \xi_j = - \bar \upsilon_j .
\ee
So we conclude that for some $j \in \{ 1,2,3\}$ both  equations \eqref{eqko6} and \eqref{eqko7} hold, and hence
$b \hatt(\a) \cap b \hatt (\b)$ is non-empty. 
This is a contradiction, and hence either $F$ or the analogous function when we swapped $x_1$ and $x_2$ must map $W$ into a proper subset of $\D$. 
\ep

\begin{lemma}
\label{lemk02}
Let $\a, \b \in \M$, and let $W = \{ (x,y,\a(x,y),\b(x,y) ) \} \cap \D^4$.
Suppose  $\hatt(\a)$ and $\hatt(\b)$ intersect at one point, but are not an edge pair.
Then $W$ is not a \car set at $0$, i.e. for all $\vare > 0$ there exists $P \in W \cap B(0,\vare)$
so that $d_W(P,0) > d_{\D^4} (P,0)$.
\end{lemma}
\bp
We proceed as in the proof of Lemma \ref{lemk0}, but
now we can find a point $\z \in \T$ so that $|A+ B \zeta| = 1 = |C + D \zeta|$.
As above, this leads to finding a $j$ so that both
\eqref{eqko6} and \eqref{eqko7} hold. This would say that $\a$ and $\beta$ form an edge pair at 
$ \upsilon_j \bar \xi_j $, a contradiction.
\ep

{\bf Lemma \ref{lemd6}.} {\em
Suppose for some pair $i,j \in \{ 1, \dots, d \}$ the set ${\EE} = \{ \a_k \in \BB : \a_k \ {\rm is\ a\ function\ of\ } x_i,x_j \}$ is non-empty.
Then
the functions in $\EE$ form an edge family.
}

\bp
Let us fix some pair of variables, say $x_1$ and $x_2$ for definiteness, and let $\BB'$ be the set of all functions in $\BB$ that
are functions only of the pair $(x_1, x_2)$. Let $I'$ be the index set for $\BB'$, so $\a_\ell \in \BB'$ if and only if $\ell \in I'$.

For each $\ell \in I'$, since it shares $b_{d-1} T(\a)$, there is an open set $U$ in $\T^d$ so that $U \cap b\overline{T(\a_\ell)} = U \cap b T(\a)$ is non-empty. We can assume that $U$ is invariant under the action of $\T$.
Let us unpack these statements. We know $\a_\ell \in \M$, so $b\overline{T(\a_\ell)} = b T(\a_\ell)$ and being in this set depends only on the ratio $\om_1 : \om_2$.
 So if $\om = (\om_1: \o_2 : \dots \o_\dimv) \in U \cap b T(\a_\ell)$ and $\a_k \in \BB \setminus \BB'$, then $\om \in T(\a_k)$.

For $\a_\ell \in \BB'$, $\a_\ell$ just depends on $x_1$ and $x_2$, so let us define $\hatt(\a_\ell(x_1, x_2)) = \{ \om_2: | \a_\ell^1(0) + \om_2 \a_\ell^2 (0) | < 1 \}.$
The intersection 
 $  \cap_{\ell \in I'} \hatt(\a_\ell(x_1,x_2) )$ is a non-empty open arc $(\om_2^\bullet, \om_2^\dagger) \subsetneq \T$.

Without loss of generality we can assume that there is some minimal bounding function $\alpha_j$ in $\BB'$ that shares
$b_{\dimv-1} T(\a)$ at a point with first coordinate $1$ whose second coordinate is $\om_2^\bullet$ (it must be $\om_2^\bullet$ or $\om_2^\dagger$ and the  case with $\om_2^\dagger$ is handled with the obvious changes).
Relabel the indices so that $I' = \{ \dimv +1, \dots, p, p+1, \dots q \}$,
where if $\dimv +1 \le \ell \le p$ then $bT(\a_\ell)$ shares $b_{\dimv -1} T(\a)$ at $(1:\om_2^\bullet: \dots) $, and if $p+1 \le \ell \le q$
then $bT(\a_\ell)$ does not share $b_{\dimv -1} T(\a)$ at $(1:\om_2^\bullet: \dots) $, but does share it at $(1: \om_2^\dagger: \dots)$.
We want to prove that either $p=q$, or all the functions in $\BB'$ have another edge point in common
where they form an edge family (as in Figure \ref{figk1}).

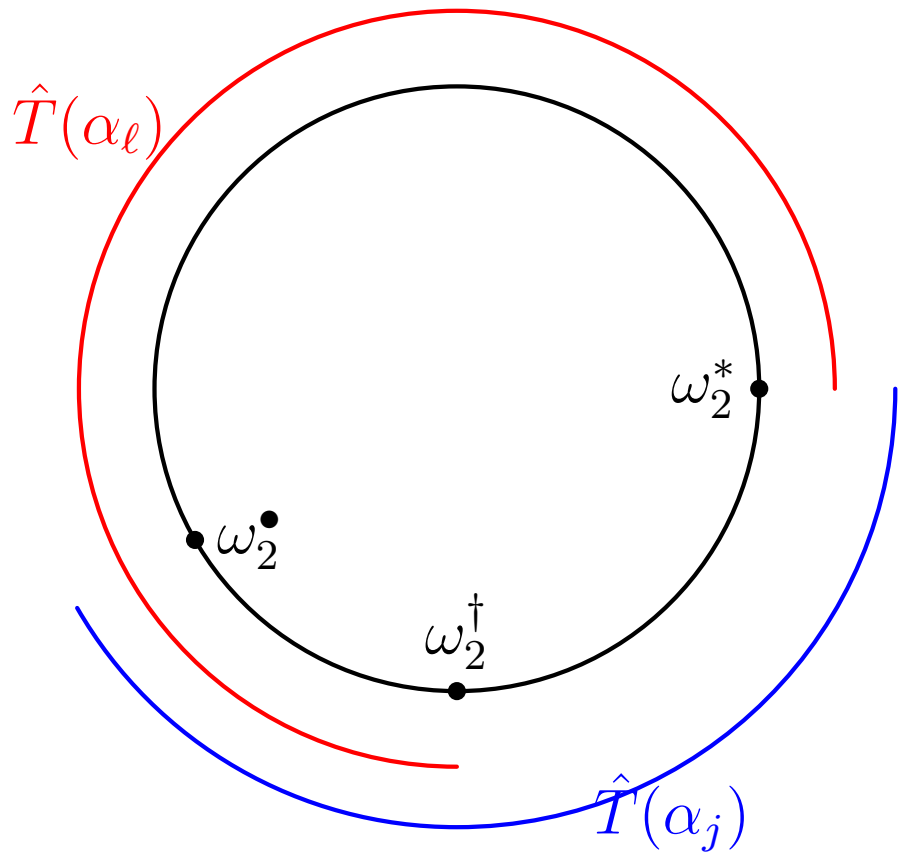
\begin{figure}[htbp]
  \centering

\begin{tikzpicture}[scale=2, line cap=round, line join=round]

  \def\R{1.0}        
  \def\RarcA{1.25}   
  \def\RarcB{1.45}   

  \def\angStar{0}    
  \def\angLeft{210}   
  \def\angRight{270}  

  \draw[thick] (0,0) circle (\R);

  \coordinate (wstar) at (\angStar:\R);
  \coordinate (wleft) at (\angLeft:\R);
  \coordinate (wright) at (\angRight:\R);

  \fill (wstar)  circle (0.03);
  \fill (wleft)  circle (0.03);
  \fill (wright) circle (0.03);

  \node[left] at (wstar) {$\omega_2^*$};
  \node[right]  at (wleft) {$\omega_2^\bullet$};
  \node[above] at (wright) {$\omega_2^\dagger$};

  \draw[thick,red]
    (0,0) ++(\angStar:\RarcA)
      arc[start angle=\angStar, end angle=\angRight, radius=\RarcA]
      node[midway, left] {$\hat T(\alpha_{\ell})$};

  \draw[thick,blue]
    (0,0) ++(\angStar:\RarcB)
      arc[start angle=\angStar, end angle=-150, radius=\RarcB]
      node[midway, right] {$\hat T(\alpha_{j})$};

\end{tikzpicture}
\caption{Edge pair with opposite orientations }
\label{figk1} 
\end{figure}

Case 1: Suppose $\a_j$ is a minimal bounding function, and $\hatt(\a_j)$ is the arc $(\o_2^\bullet, \o_2^\star) \supsetneq \hatt(\a)$. By Lemmas \ref{lemk0} and \ref{lemk02}, 
if  $p+1 \le \ell \le q$ then $\a_\ell$ and $\a_j$ must form an edge pair at $\o_2^*$.
But for any $j'$ between $d+1$ and $p$, $\a_{j'}$ must also form an edge pair with $\a_\ell$, and this must also be at $\o_2^*$. So $\BB'$ are an edge family at $\o_2^*$.

Case 2: Suppose $\a_m$ is a minimal bounding function, and $\hatt(\a_m)$ is the arc $(\o_2^\bullet, \o_2^\bullet) = \hatt(\a)$. Then for every other $\ell$ in $\{ d+1, \dots, q \}$ we must have that $\a_\ell$ and $\a_m$ form an edge pair. It is then possible  that there is some $\a_j$ that forms an edge pair with $\a_m$ at $\o_2^\bullet$, and another $\a_\ell$ that forms an edge pair with it at $\o_2^\dagger$; and since $\a_\ell$ and $\a_j$ form an edge pair themselves, that must be at $\o_2^*$.
Then by Theorem \ref{thmb3}, $\a_\ell$ and $\a_j$ are both bounded by $1$ on $D_{\a_m}$, which means that in the $\ell$ and $j$ coordinates we have graphs over the others. This is ruled out by Assumption \ref{ass2}.
\ep

\section{Proof of Lemma \ref{lemd9}}
\label{secl}

Recall that we have chosen coordinates so that the bounding functions $\BB$ are given by
\[
\BB \= \bigcup_{j=\dimv+1}^{\dimv+m_2} \EE_j \= \{ \a_\ell : \dimv+1 \le \ell \le m_3 \} ,
\]
$\BB'$ is the perhaps larger set $\cup \EE'_j$ that includes the full edge family for each $\EE_j$,
and
\[
G \= \bigcap_{\a_\ell \in \BB'} \{ x \in \D^\dimv : |\a_\ell(x) | < 1 \} .
\]

For $\a = [\tau; A,B] \in \M$, let us define $\pp(\a)$ to be the unique  point $\tau$ in $b \hatt(\a)$ with the property that for $\vare$ positive and sufficiently small, $\tau e^{i\vare} \in \hatt(\a)$;
and let $\ppm(\a)$ be the other point in $b \hatt (\a)$. If one orients the arc $\hatt(\a)$ counter-clockwise, $\pp(\a)$  is the initial point of the arc, and $\ppm(\a)$ is the final point.
Suppose $\LL$ is a \roms, given by
\be
\label{eq88.1}
\LL = \{ (x_1, x_2, \a_k(x_1, x_2)) : 3 \le k \le m \} .
\ee
Let each $\a_k = [ \tau_k; A_k, B_k]$. Since they form an edge family, there are unimodular $\o, \eta$
so that $A_k + B_k \o = \eta,\  \forall k$. Suppose $\o \in \pp(\a_k)$ for every $k$. Then for every M\"obius automorphism, each $\hatt(\tilde \a_k)$ will be an arc with the same initial point. Therefore their intersection  $\cap_{k=3}^m \hatt(\tilde \a_k)$ will just be the minimal such arc, which we know is never empty. So $(x_1, x_2)$ will be a good parametrization at every point of $V$.

\begin{lemma}
\label{leml01}
Let $\a = [ \tau; A, B]$ be in $\M$. Let $\o \in b\hatt(\a)$. Then $\o = \pp(\a)$ if and only if
$\arg( \bar A B \o) \in [0,\pi) $.
\end{lemma}
\bp
This follows from 
\[
|A + B \o e^{i\vare}|^2 = |A|^2 + |B|^2 + 2 \Re (\bar A B \o)(\cos \vare + i \sin \vare) .
\]
\ep
\begin{lemma}
\label{leml02}
Let $c_1, \dots , c_N$ be distinct complex numbers. Then there exists some choice $j,k$ so that for every $i \neq j, k$ we have 
\[
\arg \frac{c_j - c_i}{c_k - c_i} \ \in \ [-\pi,0) .
\]
\end{lemma}
\bp
Let $P$ be the convex hull of the set $\{ c_i \}$. Choose $c_j$ to be a vertex of this polygon, 
let $E$ be an edge of $P$ starting at $c_j$, and choose $c_k$ to be the other end of $E$.
If $c_i \in E$, then $\frac{c_j - c_i}{c_k - c_i} = -1$.
For every other $c_i$, they all lie in the half-plane on one side of the line that extends $E$.
Swapping $c_j$ and $c_k$ if necessary, we can assume that as 
 $t$ moves from $c_j$ to $c_k$, the argument of $t - c_i$ 
decreases for every $c_i$ not in $E$, and 
 by a total amount less than $\pi$. Therefore the difference in  argument between  $c_j - c_i$ and $c_k - c_i$ is negative, but greater than $-\pi$.
\ep

To simplify our calculations, notice that if we start with a \roms\ $\LL$ as in \eqref{eq88.1}, then applying a rotation automorphism to $\D^m$, we can assume that $\o = 1 = \eta$ and that each $\tau_k = 1$.
So we can write $\a_k = [ 1; 1-B_k, B_k]$.
\begin{lemma}
\label{leml03}
Let $\LL$ be the \roms\ given by \eqref{eq88.1}, where for each $\ell \ge 3$ we have $\a_\ell = [1; 1-B_\ell, B_\ell]$. Let us rewrite $\LL$ as a graph over $x_j, x_k$: 
\be
\label{eq88.3}
\notag
\LL = \{ x \in \D^m :  x_i = \g_i(x_j, x_k) \}
\ee
Then for every $i \neq j,k$ we have \be
\label{eq88.2}
\g_i = [ \frac{\bar B_j - \bar B_k}{B_j - B_k} ; \frac{B_i - B_k}{\bar B_j - \bar B_k} , \frac{B_j - B_i}{\bar B_j - \bar B_k} ].
\ee
If $i =1$, we take $B_1 = 0$, and if $i=2$, we take $B_2 = 1$ in \eqref{eq88.2}.
\end{lemma}
\bp
This is a special case of Corollary \ref{core1}.
\ep

\begin{cor}
\label{corl1}
Suppose $B_j$ and $B_k$ are adjacent vertices in the convex hull of $\{ B_i : 1 \le i \le m \}$.
Then if we change variables to $\g_i(x_j, x_k)$, we have that all the arcs $\hatt(\g_i)$ share a common initial point (or share a common final point) and $(x_j, x_k)$ is a good parametrization at every point of $\LL$.
\end{cor}
\bp
The first assertion follows from formula \eqref{eq88.2} and Lemma \ref{leml01}: we have either 
$(B_j-B_i)/(B_i - B_k)$ is always in the upper half-plane, or always in the lower half-plane.

The second assertion follows by continuity. As we move $a$ to get $\tilde \g_i$, the direction of the arcs from the anchor point doesn't change. If they all go in the same direction, their intersection equals the smallest such arc, which we know never goes to zero since $\g_i \in \M$. 
\ep

\begin{lemma}
\label{leml04}
Let $\LL$ be a \roms, normalized as in Lemma \ref{leml03} by
$\a_\ell = [1; 1-B_\ell, B_\ell]$.
There is some choice of $j,k$ in $\{1, \dots, m \}$ so that
if one parametrizes $\LL$ by $\g_i(x_j, x_k)$, then 
for all the functions $\g_i$ that lie in $\M$, the arcs $\hatt(\g_i)$ share  a common initial point.
\end{lemma}
\bp
If all the arcs $\hatt(\a_j)$ share $1$ as a common initial point, we take $j=1, k=2$.
If they share it as a final point, take $j=2,k=1$. Otherwise, we are in a situation like Figure \ref{figk1}.
 Choose some $j$ so that $\partial_{-} \hatt (\a_j) = \partial_{-} (\cap_{i=3}^m \hatt(\a_i))$.  Call this point $\o$. Then we know $| (1-B_j) + B_j \o| = 1$, and for all $i$ we have $| (1-B_i) + B_i \o| \le 1$. Let $\eta=  (1-B_j) + B_j \o$. Then the half-plane
$\{ z \in \C : \Re [ \bar \eta( 1- z + z \o) ] \le 1 \}$ has $B_j$ in its boundary, and no other $B_i$, since any $B_i$ in the boundary would have to satisfy $  (1-B_i) + B_i \o = \eta$. Therefore $B_j$ is a vertex of the convex hull of $\{ B_i \}$,
and so by Corollary \ref{corl1} we can choose $B_k$, an adjacent vertex, to make the 
lemma hold.
\ep

\begin{lemma}
\label{leml05}
There is some choice of coordinates so that $V$ has a good parametrization and so that each edge family $\{ \g_i \}$  has the property that  the anchor point is the initial point for  each arc $\hatt(\g_i)$.
\end{lemma}

\bp
Let us start with a good parametrization $(x_1, \dots, x_\dimv)$.
We want to show that for each edge family, we can make a choice as in Lemma \ref{leml04} with the additional property that we still have a good parametrization for all the other functions.
We will do it for each edge family in turn, so 
it is sufficient to show we can do this for the edge family $ {\mathcal E} = \{ \a_n (x_1, x_2) : n \in \I \}$.
We can normalize this family as in Lemma \ref{leml03} to have $\a_n = [1; 1-B_n, B_n]$.
We want to find $j,k \in \I \cup \{ 1,2\}$ so that if we replace $x_1$ and $x_2$ by $x_j$ and $x_k$,
(i) the new functions $\g_n$ have the property that every $\hatt(\g_n)$ has $1$ as its initial point,
and (ii) we still have a good parametrization.

Note first that if the anchor point $1$ is the initial point (or final point by symmetry) of every $\hatt(\a_n)$, we don't have to do anything. 
Otherwise, we choose $(x_j,x_k)$ so that the arcs $\hatt(\a_j) = (1, \o_2^\dagger)$ and $\hatt(\a_k) = (\o_2^\star, 1) $ are minimal,
 the first among all arcs that have $1$ as an initial point, the second among those that have it as a final point.

By Corollary \ref{cord3}, we understand the structure of $T(\a)$. Replacing $x_1, x_2$ by  coordinates $x_j, x_k$ from $\I$ will not affect $T(\a)$ outside the first two coordinates, so we only need to worry about those  $\a_p$ that 
are not in  any of the edge families, and for which at least one of $\a_p^1$ and $\a_p^2$ is non-zero.
Choose some $\a_p $ that is not in  any of the edge families;
for notational convenience, let us call this function $\psi$. We know $T(\psi) \supset \overline{T(\a)}$,
since $\psi$ does not share the boundary anywhere.

We know that if $(\o_1 : \o_2) \in \cap \{  T(\a_n(x_1, x_2)) : \a_n \in {\mathcal E} \}$
and $(\o_3 : \dots : \o_\dimv)$ is such that $(\o_1: \o_2 : \o_3 : \dots : \o_\dimv) \in T(\a)$, then
$| \sum_{i=1}^\dimv \a^i_\ell \o_i| < 1$ for all $\ell$ between $\dimv + 1$ and $\dimpol$.
Let 
\[
R = \sup \{   | \sum_{i=3}^\dimv \psi^i \o_i | :  (\o_3 : \dots : \o_d) \in \pi_{3,\dots,d} T(\a) \}.
\] 
As $\{ \psi^1 \o_1 + \psi^2 \o_2 \}$ is circled, we must have 
\be
\label{eq88.4}
\notag
|\psi^1 \o_1 + \psi^2 \o_2| + R< 1.
\ee

Let us define $D(\a_k) = \{ \o_2 \in \D: |\a_k^1 + \a_k^2 \o_2| < 1 \}$, let $S(\a_k) = bD(\a_k)$,
and let  $D_R(\psi) = \{  \o_2: |\psi^1 + \psi^2 \o_2| < 1-R \}$, and 
$S(\psi) = bD_R(\psi)$.
$D$ stands for disk and $T$ for torus, and we have $\hatt(\a_k) = D(\a_k) \cap \T$.

Now, let us swap coordinates $x_2$ and $x_j$, keeping $x_1$.
We will have new functions $\b_i(x_1, x_3, \dots, x_\dimv, x_j) = \alpha_i(x_1, \beta_2(x_1, x_j), x_3, \dots , x_\dimv)$ where
$\b_2(x_1, x_j) = \a_j^{-1}(x_1, \underline{x_j})$.
Recall that we chose $j$ so that $\hatt(\a_j)$ is an arc with initial point the anchor point $1$, and final point 
 $\o_2^\dagger$, and with the property that $\om_2^\dagger \in T(\a_n)$ for every other $n \in \I$.

Note that if $g(x_1, x_j) = f(x_1, \beta_2(x_1, x_j))$, then
for any $\tau_2$, if we define $\tau_j = \a_j^1  + \a_j^2 \tau_2 = (1-B_j) + B_j \tau_2$, then $\tau_2 = \b_2^1 + \b_2^j \tau_j$, so
\beq
g^1 1 + g^j \tau_j &\=& f^1 + f^2 \b_2^1 + f^2 \b_2^j \tau_j \\
&=& f^1 + f^2( \b_2^1 + \b_2^j \tau_j) \\
&=& f^1 + f^2 \tau_2.
\eeq
Therefore $\tau_2 \in D(f)$ if and only if $\tau_j$ is in $D(g)= \{ \o_j : | g^1 + g^j \o_j | < 1 \}$.

Define $\o_j^\dagger = (1-B_j) + B_j \om_2^\dagger$,
and $\o^\dagger = (1 : \o_3 : \dots : \o_\dimv: \o_j^\dagger)$,
where $ (\o_3 : \dots : \o_d) $ is some fixed point in $ \pi_{3,\dots,d} T(\a)$.
Then $\o_j^\dagger$ is in $\cap_{i \in \I, i \neq j} \hatt(\b_i) \cap b \hatt(\b_2)$, and 
$\o^\dagger$ is in $T(\a_\ell(x_1, x_3, \dots, x_\dimv, x_j))$ for all other $\a_\ell$.
If we perturb $\o_j^\dagger$ slightly so that $\o^\dagger$ falls into $T(\b_2)$, then we get a point
in $\cap T(\b)$ for all $\b$ in the new coordinates, and hence $\{ x_1, x_3, \dots, x_\dimv, x_j \}$ is a good parametrization.

There is a more powerful, geometric way to see this. Let us fix $\o_1 = 1$, and $\o_3, \dots ,\o_\dimv$ as above. Define $\hatd(\a) = \cap_{i=\dimv+1}^\dimpol \{ z_2 \in \D: |\a^1_i + \a^2_i z_2 | < 1 \}$. The boundary of the convex set  $\hatd(\a)$ in Figure \ref{figl1} consists of the four arcs
going through $\o_2^*, \o_2^\dagger, U_2, V_2$. Changing variables as in Figure \ref{figl2},
the new coordinates will be good if and only if the boundary of $\hatd(\beta)$ contains an arc of the unit circle
(in the picture, the boundary of  $\hatd(\beta)$ is the four arcs 
going through $\o_j^*, \o_j^\dagger, U_j, V_j$).

Claim:
$1$ is in  $S(\psi)$.

If $1$ is not in $S(\psi)$,  we have a situation like Figure \ref{figl1}.

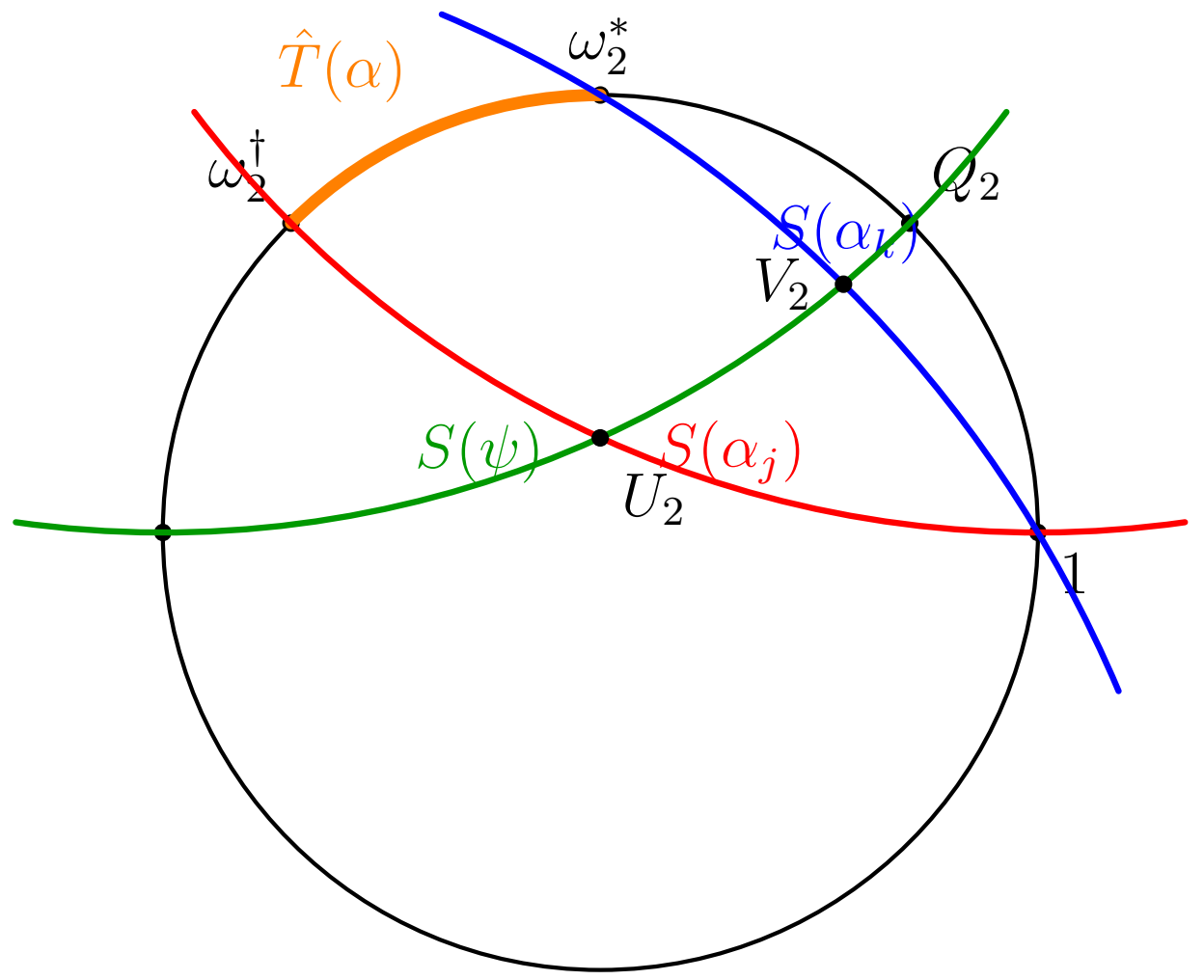
\begin{figure}[!ht]
  \centering

\begin{tikzpicture}[scale=3, line cap=round, line join=round]
\useasboundingbox (-1.25,-1.25) rectangle (1.25,1.25);

\pgfmathsetmacro{\rtwo}{sqrt(34)/2}
\pgfmathsetmacro{\rone}{1+sqrt(2)}
\pgfmathsetmacro{\rthree}{1+sqrt(2)}
\pgfmathsetmacro{\eps}{8}

\coordinate (P) at (1,0);
\coordinate (Q) at ({1/sqrt(2)},{1/sqrt(2)});
\coordinate (R) at (0,1);
\coordinate (S) at ({-1/sqrt(2)},{1/sqrt(2)});
\coordinate (Tleft) at (-1,0); 

\coordinate (C1) at (1,{1+sqrt(2)});     
\coordinate (C2) at (-1.5,-1.5);         
\coordinate (C3) at (-1,{1+sqrt(2)});    

\draw[thick] (0,0) circle (1);

\foreach \X/\pos/\lab in {
  P/below right/{$1$},
  Q/above right/{$Q_2$},
  R/above/{$\omega_2^*$},
  S/above left/{$\omega_2^\dagger$},
  Tleft/left/{}
} {
  \fill (\X) circle (0.02) node[\pos] {\lab};
}

\draw[line width=2.4pt, orange]
  (0,0) ++(90:1) arc[start angle=90, end angle=135, radius=1]
  node[midway, above left, orange] {$\hatt(\alpha)$};

\draw[very thick, red]
  (C1) ++(-90+\eps:\rone)
  arc[start angle=-90+\eps, end angle=-135-\eps, radius=\rone]
  node[midway, right] {$S(\a_j)$};

\pgfmathsetmacro{\aPtwo}{atan2(1.5,2.5)}
\pgfmathsetmacro{\aRtwo}{atan2(2.5,1.5)}
\draw[very thick, blue]
  (C2) ++(\aPtwo-\eps:\rtwo)
  arc[start angle=\aPtwo-\eps, end angle=\aRtwo+\eps, radius=\rtwo]
  node[midway, above] {$S(\a_k)$};

\draw[very thick, green!60!black]
  (C3) ++(-45+\eps:\rthree)
  arc[start angle=-45+\eps, end angle=-90-\eps, radius=\rthree]
  node[midway, left] {$S(\psi)$};

\path[name path=redcircle]   (C1) circle (\rone);
\path[name path=greencircle] (C3) circle (\rthree);

\path[name intersections={of=redcircle and greencircle, by={RGone,RGtwo}}];

\path (RGone); \pgfgetlastxy{\xA}{\yA};
\path (RGtwo); \pgfgetlastxy{\xB}{\yB};

\ifdim\yA<\yB
  \coordinate (Tint) at (RGone);
\else
  \coordinate (Tint) at (RGtwo);
\fi
\path[name path=bluecircle]  (C2) circle (\rtwo);
\path[name path=greencircle] (C3) circle (\rthree);

\path[name intersections={of=bluecircle and greencircle, by={BGone,BGtwo}}];

\path (BGone); \pgfgetlastxy{\xC}{\yC};
\path (BGtwo); \pgfgetlastxy{\xD}{\yD};

\ifdim\yC<\yD
  \coordinate (Vint) at (BGone);
\else
  \coordinate (Vint) at (BGtwo);
\fi

\fill (Vint) circle (0.02);
\node[anchor=east, xshift=-2pt] at (Vint) {$V_2$};

\fill (Tint) circle (0.02);
\node[anchor=west, yshift=-12pt] at (Tint) {$U_2$}; 

\end{tikzpicture}

\caption{In $\omega_2$ plane, $\omega_1 = 1$ }
\label{figl1} 
\end{figure}

\begin{figure}[!ht]
  \centering

\begin{tikzpicture}[scale=3, line cap=round, line join=round]
\useasboundingbox (-1.25,-1.25) rectangle (1.25,1.25);

\pgfmathsetmacro{\rone}{1+sqrt(2)}
\pgfmathsetmacro{\eps}{8}

\coordinate (P) at (1,0);
\coordinate (R) at (0,1);
\coordinate (S) at ({-1/sqrt(2)},{1/sqrt(2)});

\pgfmathsetmacro{\thetaU}{140}
\coordinate (U) at ({cos(\thetaU)},{sin(\thetaU)});

\coordinate (C1) at (1,{1+sqrt(2)});
\coordinate (C2) at (2.23750267, 6.14748805);

\pgfmathsetmacro{\rtwo}{veclen(1-2.23750267, 0-6.14748805)}

\draw[thick,red] (0,0) circle (1);

\foreach \X/\pos/\lab in {
  P/below right/{$1$},
  R/above/{$U_j$},
  S/above/{$\omega_j^\dagger$}
} {
  \fill (\X) circle (0.02) node[\pos] {\lab};
}

\draw[very thick]
  (C1) ++(-90+\eps:\rone)
  arc[start angle=-90+\eps, end angle=-135-\eps, radius=\rone]
  node[midway, below] {$S(\b_2)$};

\path[name path=blackcircle] (C1) circle (\rone);
\path[name path=bluecircle]  (C2) circle (\rtwo);

\path[name intersections={of=blackcircle and bluecircle, by={X,P}}];

\fill (X) circle (0.02) node[shift={(-7pt,-7pt)}] {$\omega_j^*$};

\pgfmathsetmacro{\aPtwo}{atan2(0-6.14748805, 1-2.23750267)}
\pgfmathsetmacro{\aUtwo}{atan2(sin(\thetaU)-6.14748805, cos(\thetaU)-2.23750267)}

\draw[very thick, blue]
  (C2) ++(\aPtwo+\eps:\rtwo)
  arc[start angle=\aPtwo+\eps, end angle=\aUtwo-\eps, radius=\rtwo]
  node[midway, above] {$S(\b_k)$};

\coordinate (Cg) at (-2,0);
\pgfmathsetmacro{\rg}{sqrt(5)}
\path[name path=greencircle] (Cg) circle (\rg);

\coordinate (Gtop) at (0,1);   
\coordinate (Gbot) at (0,-1);

\pgfmathsetmacro{\angG}{atan2(1,2)} 

\draw[very thick, green!60!black]
  (Cg) ++(\angG:\rg)
  arc[start angle=\angG, end angle=-\angG, radius=\rg];

\path (X); \pgfgetlastxy{\xX}{\yX};
\path (S); \pgfgetlastxy{\xS}{\yS};
\path (C1); \pgfgetlastxy{\xC}{\yC};

\pgfmathsetmacro{\angXblack}{atan2(\yX-\yC,\xX-\xC)}
\pgfmathsetmacro{\angSblack}{atan2(\yS-\yC,\xS-\xC)}

\draw[very thick, orange]
  (C1) ++(\angXblack:\rone)
  arc[start angle=\angXblack, end angle=\angSblack, radius=\rone];

\path[name intersections={of=greencircle and blackcircle, by={GbK1,GbK2}}];
\path (GbK1); \pgfgetlastxy{\xA}{\yA};
\path (GbK2); \pgfgetlastxy{\xB}{\yB};
\ifdim\xA>\xB
  \coordinate (Gblack) at (GbK1);
\else
  \coordinate (Gblack) at (GbK2);
\fi

\path[name intersections={of=greencircle and bluecircle, by={GbB1,GbB2}}];
\path (GbB1); \pgfgetlastxy{\xC}{\yC};
\path (GbB2); \pgfgetlastxy{\xD}{\yD};
\ifdim\xC>\xD
  \coordinate (Gblue) at (GbB1);
\else
  \coordinate (Gblue) at (GbB2);
\fi

\fill (Gblack) circle (0.02) node[shift={(8pt,10pt)}] {$V_j$};   
\fill (Gblue)  circle (0.02) node[shift={(10pt,-16pt)}] {$Q_j$}; 

\end{tikzpicture}

\caption{In $\omega_j$ plane, $\omega_1 = 1$ }
\label{figl2} 
\end{figure}
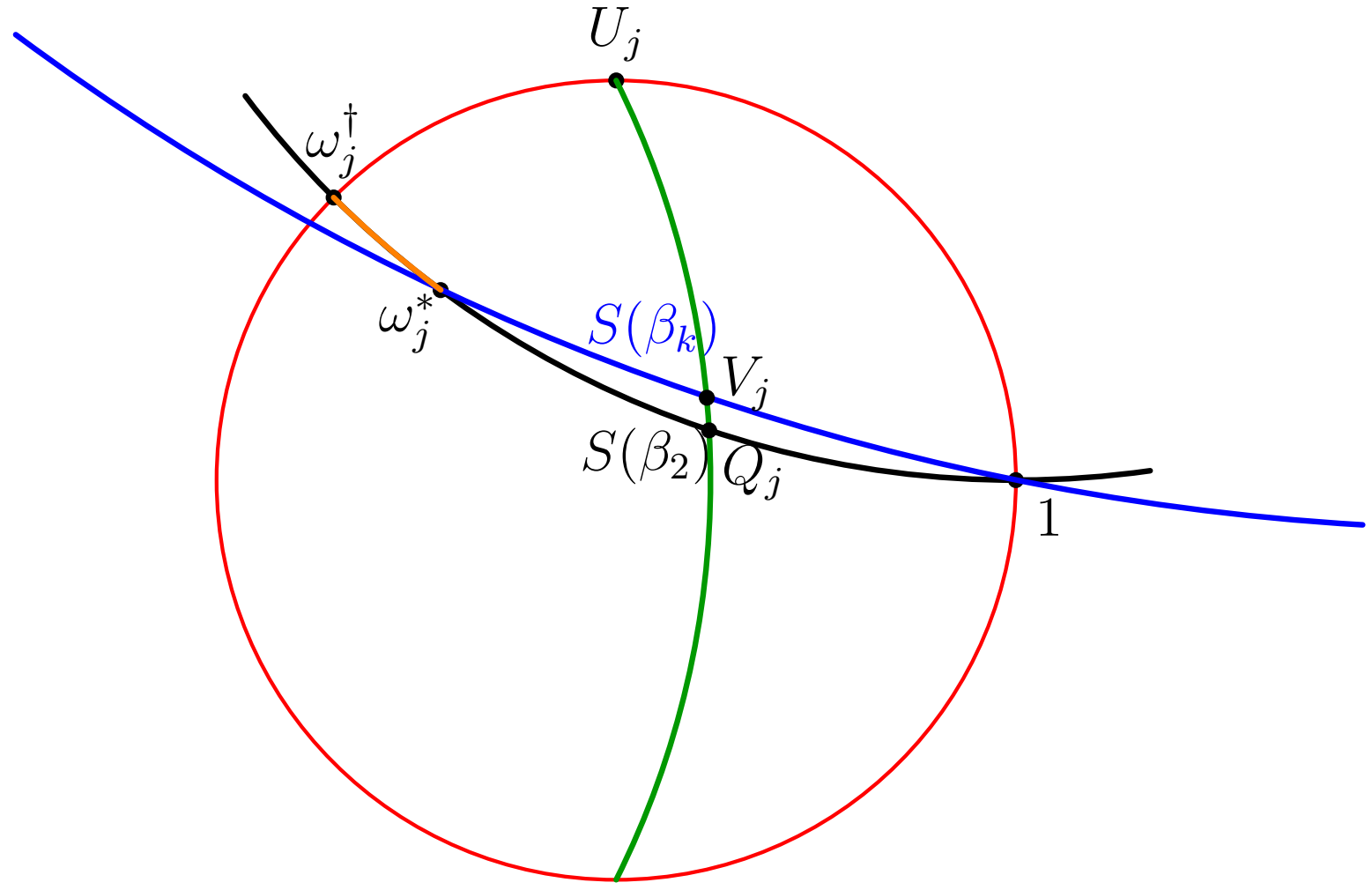

Suppose first that the only functions that provide part of the boundary of $\hatd(\a)$ are 
$x_2, \a_j, \a_k, \psi$, as in the picture Figure \ref{figl1}.
Then when we change coordinates to $(x_1, x_j)$ as in Figure \ref{figl2}, we see
that $\phi = \psi(x_1, \b_2(x_1, x_j), x_3, \dots, x_\dimv)$ is bounding at $U$. This means that $\phi(x_1, x_j) \in \M$, and 
$\phi, \beta_2$ and $\beta_k$ are all in an edge family by Lemma \ref{lemd6}.
Since $b\hatt(\beta_2)$ and $b\hatt(\beta_k)$ only have the point $1$ in common,
this must also be in $b\hatt(\phi)$, and hence $1$ was in $b\hatt(\psi)$ all along.
This is a contradiction.

If there are more functions, a similar argument holds. In Figure \ref{figl1}, draw straight lines from $1$ to $\o_2^\dagger$ and $\o_2^*$. If there is some $\psi$ with $1$ not in $D_R(\psi)$, choose one with $S(\psi) \cap \T$ closest to $\hatt(\a)$. There will be an arc of $b\hatd(\a)$ that comes from
$b(D_R(\psi))$ (the green arc in the pictures). This will be connected to arcs in $S(\a_{j'})$ and $S(\a_{k'})$ that are also in
$b\hatd(\a)$. Now we repeat the argument with $j'$ and $k'$ to conclude that $1$ had to be in $S(\psi)$ all along.
\oec

So we can conclude that $(1, \dots, 1)$ is in $T(\a_\ell)$ for every $\ell$ not in $\I$.
Now we can use Corollary \ref{corl1} to choose a pair $x_j, x_k$ 
so that all the arcs $\{ \hatt(\g_i) : i \in \I\}$ have a common initial point $1$, and this will be a good parametrization since $(1, \dots, 1) \in T(\g_\ell)$ for all other $\g$.
Now we repeat this argument with the other sets of coordinates associated with the other edge families.
\ep

The importance of all the arcs going in the same direction is the following.
\begin{cor}
\label{corl06}
There is some choice of coordinates $\{ x_{i_1}, \dots, x_{i_\dimv}  \}$ that is good for every M\"obius automorphism applied to $V$.
\end{cor}
\bp
Choose a good parametrization at $0$ as in Lemma \ref{leml05}.
We can assume this is $\{ x_1, \dots, x_\dimv \}$. By continuity, at every $a \in V$,
all the edge families that depend on the pair $x_{2i-1}, x_{2i}$ have $\hatt(\a_j(x_{2i-1}, x_{2i}))$ go in the same direction from the anchor point, so their intersection is just the minimal one. This is never empty, since every function is in $\M$. So $T(\tilde{\a})$ will be a product of bands, and non-empty.
\ep

\vskip 20pt
{\bf Lemma \ref{lemd9}.}
{\em
Let $j' > m_3$. Then $\a_{j'} \in O(G, \D)$.
}

\bp
Let us write $\psi$ for $\a_{j'}$.

Let $V_1$ be the set of points in $V_{reg}$ where locally $V$ can be parametrized by 
the first $\dimv$ coordinates. Note that $V \setminus V_1$ is an algebraic set of dimension at most 
$\dimv-1$, so $V_1$ is connected.

We know $\psi$ is not bounding, and
it is  a Nash function that is holomorphic near $0$.
Suppose it does not extend to be holomorphic from $G$ to $\D$.
This can happen because either 

(i) there is a point $a \in G$ that is a singular point for
$\psi$, or 

(ii)  there is a point $a \in G$  where $\psi$ has modulus $1$.

In the first instance, $P= (a, \a(a))$ is
in $V \setminus V_1$.
In the second, we can assume $P \in b V$
(by taking the first point on a path where this happens).

Claim: (i) cannot happen.

Move $P$ to $0$ by a M\"obius automorphism $\Phi$. 
Then $0$ is not a regular point for 
every $\tilde{\alpha}_\ell$.
But it is a regular point for the lower dimensional set $\{ (x_1, \dots, x_\dimv, \tilde{\alpha}_\ell (x)) : \alpha_\ell \in \BB' \}$.
By the argument in Cor. \ref{corl06},
 we can assume $\cap \{ T(\tilde{\a}_\ell) : \a_\ell \in \BB' \} \neq \emptyset$.

 Choose a point $Q = (\om_0 \cdot \l_0, \tilde{\a}(\om_0 \cdot \l_0))$ that is regular for $\tilde V = \Phi(V)$, with $\l_0\in \D$ and $\om_0 \in \cap \{ T(\tilde{\a}_\ell) : \a_\ell \in \BB' \}$.
Now apply another  automorphism interchanging $0$ and $Q$.
In our new coordinate system (which we shall  denote $\tilde{\tilde{\a}}$), the point $P$ has now moved to $Q = (\om_0 \cdot \l_0, \tilde{\tilde{\a}}(\om_0 \cdot \l_0))$ and 
$0$ is a regular point.
If we choose $|\l_0|$ small enough, we can also ensure that $\{ x_1, \dots, x_d \}$ is a good parametrization at $Q$
 and $\om_0 \in \cap \{ T(\tilde{\tilde{\a}}_\ell) : \a_\ell \in \BB' \}$.

Moreover, we must have $ \cap \{ T(\tilde{\tilde{\a}}_\ell) : \a_\ell \in \BB' \} \subseteq T(\tilde{\tilde{\psi}})$,
since otherwise $\tilde{\tilde{\psi}}$ would be bounding, and hence in $\M$. This would mean $\psi \in \M$, and so has no singular points in the closure of the set where it is bounded by $1$.

Since $0$ is regular, we can now apply Lemma \ref{lemd2}, to conclude that 
$\tilde{\tilde{ \psi}}(\om \cdot \l)$ is holomorphic 
for $\l \in \D$, and for a non-empty  open set of $\om$.

Now let $\vare_i$ be close to $0$.  Let the coordinates of $\o\in \cap \{ T(\tilde{\tilde{\a}}_\ell) : \a_\ell \in \BB' \}$ be denoted by $\o_i$. 
Then $\l \mapsto \tilde{\tilde{\psi}}(\vp_{\vare_i} (\om_i \l))$ is
holomorphic for $\l \in \D$, and with $\l = \l_0$ the set $\{  (\vp_{\vare_i} (\om_i \l_0)) \}$ is a neighborhood of $\om_0 \cdot \l_0$.
This means that $\psi$ is holomorphic 
on one-dimensional slices in many directions through $a$.

By Lemma \ref{leml08}, $\psi$ is a holomorphic function on all of $G$ except possibly on $Y \cap G$, a set of dimension at most $d-2$. Since $G \setminus Y$ is simply connected, the chosen branch has trivial  monodromy (i.e. analytic continuation along closed loops is single-valued), whence  $\psi$ is a well-defined bounded holomorphic function on $G \setminus Y$. So by Riemann's extension theorem, $\psi$ extends analytically to all of $G$.

Claim: (ii) cannot happen.

By Corollary \ref{corl06},   we  can assume  $\{x_1, \dots, x_\dimv \}$
are a good 
parametrization for every M\"obius automorphism of $V$.
Therefore 
 $\cap_{\a_\ell \in \BB'} T(\tilde \a_\ell)$ is non-empty  at $P$.
By  \eqref{eq21}, we know that $T(\tilde{\psi})$ will shrink to nothing at $P$. So as we move to $P$ along a path inside $V_1$, at some stage a new function must become bounding. 
It must be in $\M$, and it either depends on a new pair of variables, or an old pair.
Neither of these is possible.
In the former case, there would be a constraint in $T(\tilde \alpha)$ involving the new pair of variables, and hence this would have been a constraint in $T(\alpha)$ always.
In the latter case, the function would have to lie in one of the edge families ${\mathcal E'}_j$ by
Lemmas \ref{lemd6} and \ref{lemd7}, so would have been in $\BB'$ all along.
\ep

\begin{lemma}\label{leml08}
Let $\psi$ be a Nash function of $x= (x_1, \dots, x_\dimv)$, and $F$ be an irreducible polynomial so that
$F(x,\psi(x)) = 0$. Let $X$ be the locus set of $F\in \mathbb C[x][z]$, i.e.  the projection of $W$ onto the first $d$ coordinates, where  \[
W = \{( x,w)\in \C^{d+1}: F(x,w) = 0, \frac{\partial F}{\partial w} (x,w) = 0 \}.
\] Let $Y$ be the set of non-smooth points of $X$.
Let $\l_0 \in \D, \ \o_0 \in \T^\dimv$, and let $a = \o_0 \cdot \l_0$. 
 Assume for all $\vare_i$ sufficiently close to 0 and $\o$ close to $\o_0$, the function $\psi(\phi_{\vare_i}(\o_i \l))$ is holomorphic in
 $\l$ for all $\l \in \D$.
 
 If $a$ is not in $Y$, then $\psi$ is holomorphic at $a$.
\end{lemma}

Note that $X$ is where the discriminant of $F$, thought of as an element of $\C[x][w]$, vanishes, so $X$ is  an algebraic set of dimension $d-1$, and $Y$ is of dimension at most $d-2$ (as a set of singular points of $X$).

\bp
If $a $ is not in $X$, then $\psi$ is holomorphic by the implicit function theorem. 
Suppose that $a \in X \setminus Y$. We shall use a local argument to show that $\psi$ has a Puiseux series expansion near $a$, and then use the fact that $\psi$ is holomorphic on disks through $a$ to conclude that it is actually analytic.

There is a local biholomorphic change of variables that maps $a$ to $0$ and 
the locus set $X$ to $\{ x_1 = 0 \}$.
Under the change of variables $F$ becomes $H \in {\mathcal O}(x)[w]$, and we can assume that $H$ is monic in $w$; the function $\psi$ becomes $\chi$ after the change of variables.
We write $x = (x_1, x')$. The family of analytic discs $\l \mapsto (\vp_{\vare_i} (\om_i \l))$ is also transformed via this biholomorphism. We can choose a subfamily that intersect $\{x_1=0\}$ transversely and whose union covers $\{0\}\times U'$, where $U'$ is a neighborhood of $0$ in $\mathbb C^{d-1}$. In particular, among them we can distinguish a family of the form $\l\mapsto (\l, \sigma_{x'}(\l))$, where $\sigma_{x'}(0)=x'$.

We now use the Puiseux theorem with a parameter in the form that is stated in \cite[Theorem 4.24]{den08}  (see also \cite[p. 170]{loj91}, or rather its proof, since the theorem as stated assumes that $d=2$). Since the locus set is $x_1=0$, the discriminant of $H$ does not vanish if $x_1\neq 0$ is close to $0$, so Theorem 4.24 applies.
Then there are germs $h_\mu$ of holomorphic functions near $(0,0)$ and $p\in \mathbb N$ such that
$$H(x_1^p, x', w) = \prod_{\mu = 0}^{p-1} (w - h_\mu( x_1, x')).$$
Consequently, $\chi(x_1^p, x') = h_{\mu_0}(x_1, x')$ for some $\mu_0$. Let $h$ denote $h_{\mu_0}$.

Let us expand $h(x_1,x) = \sum_n x_1^n \mathfrak h_n(x')$. Our hypothesis is that $\chi$ is holomorphic on one-dimensional curves, which implies that $\sum \l^{n/p} \mathfrak h_n(\sigma_{x'}(\l))$ is holomorphic near $0$. This is possible only if $\mathfrak h_n$ vanish whenever $n\not\equiv 0 \pmod p$. This proves that $\chi$ is holomorphic at $0$.
\ep

\section{Proof of Lemma \ref{lemd10}}
\label{secm}

\begin{lemma}
\label{LEMp3}
Let $g: \D \to \D$ be holomorphic, with $g(0) = 0$ and $|g'(0)| = c < 1$.
Let $0 < R < 1$.
 Let $\vare $ satisfy 
\be
\label{eq98}
0 <  \vare < \frac{R(1-c)}{8} .
\ee
If  $\eta,\zeta \in \C$ and $|\eta|, |\zeta| < \vare$, then for all $\l$ with 
 $1 > |\l| > R$ we have $\rho(g(\l), \eta) < \rho(\l,\zeta)$.
\end{lemma}
\bp
Let $g'(0) = \o c$ with $|\o| = 1$.
By Schwarz's lemma, we have
\be
\label{eq99}
\left| \vp_{-\o c}( \frac{g(\l)}{\l} ) \right| \ \le \ |\l|
\ee
for all $\l$ in $\D$.
Let 
\[
M \= \max \{ |\vp_{-\o c}(\l) : |\l| \le s \} \= \frac{c+s}{1 + sc}.
\]

Let $|\l| = s > R$.
Then \eqref{eq99} says that   $\left| \vp_{-\o c}( \frac{g(\l)}{\l} ) \right|  \le s$, so
$\left|  \frac{g(\l)}{\l}  \right|  \le M$, so
\[
| g (\l) | \ \leq \ s  \frac{c+s}{1 + sc}.
\]
Then
\be
\label{eq100}
 \rho(g(\l), \eta) \ \le \ \frac{ s \frac{c+s}{1 + sc} + \vare}{1 + \vare s \frac{c+s}{1 + sc}}.
\ee
We also have
\be
\label{eq101}
\rho(\l, \zeta) \ \ge \
s \frac{1 - \frac \vare s}{1 - \vare s } .
\ee
Comparing \eqref{eq100} and \eqref{eq101}, we want to show that
\be
\label{eq102}
s ( 1 - \frac \vare s) (1 + sc + \vare s (c+s) ) \ \stackrel{?}{>} \
(1 - \vare s) (s(c+s) + \vare (1 + sc) ) .
\ee
Dividing by $s$ and rewriting \eqref{eq102}, we get
\beq
(1-s)(1-c) & \ \stackrel{?}{>} &\ \vare \left[ \frac{2}{s} + 2 c -2s(c+s)\right] 
 - \vare^2 \left[ (1-s)(1-c) 
\right]
\\
\Leftrightarrow (1-c) & \ \stackrel{?}{>} &\ \frac{2\vare}{s} \left[ (1+ s + s^2) + cs \right] 
 - \vare^2 (1-c) ,
\eeq
and this follows from  hypothesis \eqref{eq98}. 
\ep

\begin{prop}
\label{PROm1}
Suppose $V$ is an irreducible  $\dimv$-dimensional \car set in $\D^\dimpol$ with $0 \in V$.
Let $X$ be a smooth variety in a neighborhood of $0$, with $0 \in X$ and $X \not \subseteq V$.
Any \car set $W$ containing $V \cup X$ has dimension larger than $ \dimv$.
\end{prop}
\bp
Assume that the dimension is $\dimv$. We will obtain a contradiction.
Let $V = \{ (x_1, \dots, x_\dimv, \a_j(x) ) \ : \dimv +1 \le j \le \dimpol  \}$.
We break the proof into 3 steps. In the first, we assume that no $\a_j$ is a coordinate function.
In the second, we drop that assumption, but take $d=2$. In the third, we consider the most general case.

{\bf Step 1:}
Assume $V$ satisfies the second part of  Assumption \ref{ass1}, i.e. 
none of the $\a_j$ are coordinate functions.

Let $\omega \in T(\alpha)$. Let $\vare >0$ be such that for all $b = (b', b_{\dimv +1}, \dots, b_{\dimpol}) $ 
with $\| b \| \le \vare$, if we set 
  $\tilde \a_j = \vp_{\a_j(b')} \circ \a_j ( \vp_{-b_1} (x_1), \dots , \vp_{-b_\dimv}(x_\dimv)) $, then we have $\omega \in T(\tilde \a)$
  and there exists $c_0 < 1$ so that 
  $|\sum_{i=1}^\dimv \tilde \a_j^i (0) \omega_i | \le c_0 $ for each $\dimv +1 \le j \le \dimpol$.
  
 Let $b \in X\setminus V$ have $ \| b \| < \vare$ and $|\vp_{\a_j(b')}(b_j) | < \vare$ for each $\dimv +1 \le j \le \dimpol$.
For $\l \in \D$, look at 
\begin{eqnarray}
\notag
\lefteqn{
d( ( \vp_{-b_i}(\om_i \l),  \a_j (  \vp_{-b_i}(\om_i \l))) , (b', b_{\dimv +1}, \dots, b_N))}
\\
&=&
\max [ \max_{1 \le i \le \dimv} \rho(\vp_{-b_i}(\om_i \l), b_i) , \max_{\dimv +1 \le j \le \dimpol} \rho( \a_j (  \vp_{-b_i}(\om_i \l)), b_j) ].
\label{eq89}
\end{eqnarray}
For $1 \le i \le \dimv$ the distance in \eqref{eq89} is $|\l|$.
When $\l=0$, the distance in \eqref{eq89} is attained for some $j \ge \dimv +1$.

Shrinking $\vare$ if necessary, choose $s > \frac 12$ and less than $1 - \frac{8\vare}{1 - c_0}$.
We can now apply Lemma \ref{LEMp3} to $\psi(\l) =\tilde \a_j (\om \cdot \l)$,  with $\eta = \vp_{\a_j (b')}(b_j)$ and $\zeta = 0$, to conclude that for all $\l$ with $| \l | = s$, we have 
\be
\label{eq90}
\rho(\a_j(\vp_{-b_i}(\o_i \l)), b_j) < |\l| .
\ee
Therefore by continuity there is some $\l = \l(b)$ where
\eqref{eq89} is attained in at least $\dimv +1$ slots, including the first $\dimv$.

Now repeat this for some path $b(t) \to 0$. 
We get a $(\dimv +1)$ (or more)  balanced pair for $b(t)$ and a point of
the form $(\vp_{-b'}(\omega \cdot \lambda_t), \alpha(\vp_{-b'}( (\omega \cdot \lambda_t)))$.

From \eqref{eq90} we know there is a closed curve around $0$ in the disk of radius $s$ such that for $\l$ on this curve
\be
\label{eq91}
\max_{\dimv +1 \le j \le \dimpol} \rho(\a_j(\vp_{-b_i}(\o_i \l)), b_j)  = |\l| .
\ee
There is an arc for which the maximum on the left-hand side of \eqref{eq91} is attained for some $j$, which without loss of generality we can assume is $\dimv +1$.
For each point $\l$  on this arc, there is some unimodular $\o_{\dimv +1}(\l)$ such that
\be
\label{eq92}
\a_{\dimv +1}( \vp_{-b_1} (\omega_1 \l), \dots, \vp_{-b_\dimv} (\omega_\dimv \l)) \= \vp_{b_{\dimv +1}} (\omega_{\dimv +1} \l) .
\ee
If $\o_{\dimv +1}(\l)$ is constant on the arc, then by analyticity \eqref{eq92} holds for all $\l$ in $\D$.
Differentiate \eqref{eq92} with respect to $\l$ and evaluate at $\l = 0$.
We get
\be
\label{eq93}
\notag
\sum_{i=1}^\dimv \a_{\dimv +1}^i (b) (1- |b_i|^2) \omega_i = (1- |b_{\dimv +1}|^2) \omega_{\dimv +1} .
\ee
Now take a limit as $b(t) \to 0$.
The left-hand side tends to $\sum \a_{\dimv +1}^i (0) \omega_i$, which has modulus less than $1$ since $\om \in T(\a)$,
and the right-hand side has modulus $1$. This is a  contradiction, so for $b$ sufficiently close to $0$ we must have
that $\om_{\dimv +1}$ is not constant.

If $\om_{\dimv +1}$ is not constant, we can use \eqref{eq92} to solve for $\om_{\dimv +1}$ in terms of $\l$. 
\[
\frac{1}{\l}
\vp_{-b_{\dimv +1}} \circ
\a_{\dimv +1}( \vp_{-b_1} (\omega_1 \l), \dots, \vp_{-b_\dimv} (\omega_\dimv \l)) \= \omega_{\dimv +1} 
\]
In particular, we get, for a fixed small $b$, an open set of $\om_{\dimv +1}$ satisfying \eqref{eq92} for some $\l$ that also satisfies \eqref{eq91}. Perturbing $(\om_1, \dots, \om_\dimv)$ we conclude that there is an open set of $(\om_1, \dots, \om_\dimv, \om_{\dimv +1})$ in $\T^{\dimv +1}$ for which both these equations hold.

Now apply an automorphism $\Phi = (\vp_{b_1}, \dots, \vp_{b_N})$ of $\D^\dimpol$ that moves $b$ to $0$. 
We have $\Phi(V) $ contains a point $(\om_1 \l, \dots, \om_{\dimv +1} \l, y_{\dimv+2} , \dots, y_N )$
where $|y_j| \le |\l|$ for each $j \ge \dimv+2$.

Now we argue as in the proof of Lemma \ref{lemg3} that $\overline{\pi_{\dimv +1}(W)}$ contains a non-empty open set
in $\T^{\dimv +1}$, so by Chevalley's Theorem \ref{thmc7} the dimension of $W$ must be at least $\dimv +1$.

{\bf Step 2:}  Suppose $\dimv=2$. We will prove the case  where we allow some $\a_j$ to be  unimodular multiples of coordinate functions by induction on $\dimpol$.

We can assume that $V$ has the form
\be
\label{eq94}
\{ (x_1, x_2,  \a_j(x) ) : x \in U \} 
\ee
where for $3 \le j \le p-1$ the functions $\a_j$ are unimodular multiples of coordinate functions, and for $p \le j \le \dimpol$
they are not. 
By a rotation automorphism we can assume each $\a_j$ is either $x_1$ or $x_2$ for $ j < p$.
Let $b \in X \setminus V$.

 Let $\omega \in \cap_{j=p}^\dimpol T(\a_j)$.
Choose $ \frac 12 < s < 1$, let $c = \max_{p \le j \le \dimpol } | \om_1 \a_j^1 (0) + \om_2 \a_j^2(0)|$, 
and assume $d(b,V)$ is smaller than the right-hand side of \eqref{eq98}.

Consider
\be
\notag
d(b, (\om \cdot \l, \alpha( \om \cdot \l)) .
\ee
Then for $|\l| = s$, this will be active only in the slots $\{1, \dots, p-1 \}$ by Lemma \ref{LEMp3}.
Rotating around the circle $|\l| = s$, we will reach a value of $\l$ where it is active with respect to at least 2 different indices.

Claim: we can assume these two different indices correspond to the two different coordinate functions.

Indeed, let 
\beq
J_1 &=& \{ 1 \} \cup \{ j: 3 \le j \le p-1,\ {\rm and\ } \a_j(x)= x_1 \} \\
J_2 &=&
 \{ 2 \} \cup \{ j: 3 \le j \le p-1,\ {\rm and\ } \a_j(x)= x_2\} .
 \eeq
For any given $\l$ on the circle of radius $s$, we have
\be
\label{eq95}
d(b, (\om \cdot \l, \alpha( \om \cdot \l)) =
\max_{1 \le j \le p-1} \rho( b_j, \om_j \l) .
\ee
Suppose the maximum always occurs for $j \in J_1$.

 Let $(1:\tau) \in \cap_{j=p}^\dimpol T(\a_j)$. Then the disk $\{ (\mu, \tau \mu) : | \mu | < 1 \} \subset U$ (where $U$ is the base of $V$, from \eqref{eq94}).
Fix some $\mu \in \D \setminus \{0\}$.
 Let us apply the M\"obius automorphism $\Psi = (\vp_{\mu}, \vp_{\tau \mu}, \vp_{\a_j(\mu, \tau \mu)})$ to $\D^\dimpol$.
 Then for $j \in J_1$, we have $b_j$ goes to $\vp_{\mu} (b_j)$, and for $j \in J_2$, $b_j$
 goes to $\vp_{\tau \mu}(b_j)$.
 The set  $\cap_{j=p}^\dimpol T(\tilde \a_j)$ will still be non-empty.
 So \eqref{eq95} will become
 \begin{multline}
\label{EQm2} \notag
d(\tilde b, (\om \cdot \l, \tilde \alpha( \om \cdot \l)) = \max_{1 \le j \le p-1} \rho( \tilde b_j, \om_j \l) 
\\
= \max ( \max_{j \in J_1} \rho( \bar \om_j  \vp_{\mu}(b_j), \l),
\max_{j \in J_2} \rho( \bar \om_j  \vp_{\tau \mu}(b_j), \l) )
\end{multline}
For fixed $\mu$ and $\tau$, choose $\omega$ so that $\{ \bar \om_j \mu : j \in J_1 \}$ and
$\{ \bar \om_j \tau \mu : j \in J_2 \}$ are disjoint. As $b \to 0$, we have $\vp_\mu(b_j) \to -\mu$ and
$\vp_{\tau \mu}(b_j) \to -\tau \mu$. This means the points from $J_1$ and the points from $J_2$ tend to a set that is on the circle of radius $|\mu|$.  So we can assume we can find a $\l$
on the circle of radius $s$ so that \eqref{eq95} is attained for at least 2 values of $j$, one in $J_1$ and one in $J_2$, since as we go around the circle of radius $s$ we cannot always be closest to the points from $J_1$. Relabelling the coordinates if necessary, we can assume these are the values $j=1$ and $j=2$.
\oec

Moreover, as we move along the arc $\{\l :  \rho(\bar \om_1 b_1,  \l) = \rho(\bar \om_2 b_2,  \l) \}$ there will be some third value $3 \le  j_3 \le p-1$  where
$\rho(\bar \om_1 b_1,  \l) = \rho(\bar \om_{j_3} b_{j_3},  \l)$ (see Figure \ref{figm1}).  In particular, as we move along this arc there has to be a point at
which the pair $b$ and $P = (\om \cdot \l, \a (\om \cdot \l))$ is at least 3-balanced, with respect to $x_1, x_2$ and $x_j$ for some $j > 2$ (which may or may not also be smaller than $p$).

\begin{figure}[htbp]
  \centering
\centerline{\includegraphics[scale = 0.11, angle = 0]{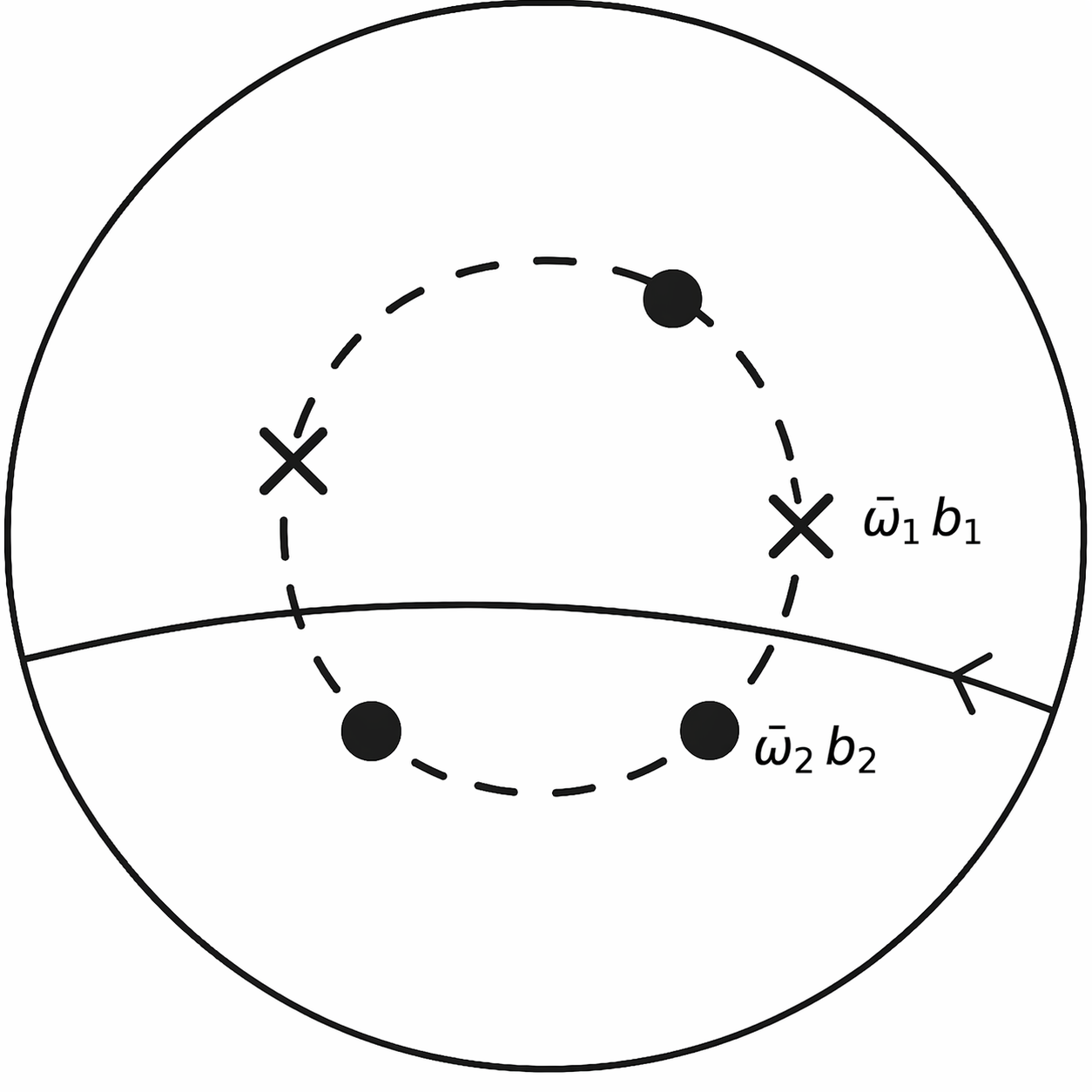}}
\caption{The arc $\{\l :  \rho(\bar \om_1 b_1,  \l) = \rho(\bar \om_2 b_2,  \l) \}$}
\label{figm1}
\end{figure}

If the pair $(b,P)$ is exactly 3-balanced, we can apply Lemma \ref{lemg7}, and conclude that either ${\rm dim}(W) > \dimv$ or
\be
\label{eq96}
\a_j(\vp_{b_1}(\om_1 \l), \vp_{b_2}(\om_2 \l)) \= \vp_{b_j} (\om_j \l) \qquad \forall \l \in \D .
\ee
There are now two subcases. 

If $j < p$, then $\a_j$ is a coordinate function, say $\a_j(x_1, x_2) = x_1$ for definiteness.
Then \eqref{eq96} says $b_j = b_1$ and $\om_j = \om_1$. If this happens for a non-discrete set of $b$'s, we can project onto the $\dimpol-1$ coordinates that do not include $j$, and proceed by induction.

If $j \ge p$, differentiate \eqref{eq96} with respect to $\l$ and evaluate at $\l = 0$.
We get
\be
\label{eq97}
\a_j^1 (b_1, b_2) \om_1 (1-|b_1|^2 ) + \a_j^2 (b_1, b_2) \om_2 (1- |b_2|^2 ) 
\= \om_j (1- |b_j|^2 ).
\ee
Now let $b \to 0$. The left-hand side will tend to a limit of modulus less than one, since $(\om_1:\om_2) \in T(\a_j)$, and the right-hand side will have a unimodular limit. 
This contradiction means that ${\rm dim}(W) > \dimv$.

 If $(b,P)$ is  $4$-balanced or more, perturb $\omega$ slightly.
Either we now get a 3-balanced pair, and apply Lemma \ref{lemg7}; or we get a continuum where Lemma \ref{lemg9} applies.
Now we proceed as in the 3-balanced case. If some index $2 < j < p$ is active in the balanced pair $(b,P)$ then we can remove this coordinate and reduce to $(\dimpol-1)$ dimensions. If the active indices are $\ge p$, then again we have
\eqref{eq96} and \eqref{eq97} and we get the same contradiction as $b \to 0$.

{\bf Step 3:}
The most general $V$ is of the following form.
There are {\roms}s ${\mathcal L}_r$, for $1 \le r \le m_1$, which have the form
\[
{\mathcal L}_r \= \{ (y,z, \alpha_j(y,z) ) : (y,z) \in U_r, 1 \le j \le h_r \} .
\]
The set $\{ \a_j : 1 \le j \le h_r \}$ is an edge family, and $U_r = \{ (y,z) \in \D^2 : |\a_j(y,z)| < 1, 1 \le j \le h_r \} $.
Then $V$ is a graph over ${\mathcal L}_1 \times \dots  \times {\mathcal L}_{m_1} \times \D^{m_2}$,
where $k = 2m_1 + m_2$.
Let us write
\be
\label{eq97.1}
V = \{ (x_1, \dots, x_{\dimv}, \a_j(x) ) : x \in U_1 \times \dots \times U_{m_1} \times \D^{m_2} \} .
\ee
We can assume that $\a_j$ is a coordinate function only in the range $\dimv +1 \le j \le p-1$,
and by Step 1 we can assume $p-1 > \dimv$.
For $1 \le r \le \dimv$, let 
\[
J_r = \{ r \} \cup \{ j: \dimv +1 \le j \le p-1,\ {\rm and\ } |\a_j(x)|= |x_r| \}
\]

Let $\om \in \cap_{j=p}^\dimpol T(\a_j)$. Let $c = \max_{p \le j \le \dimpol} |\sum_{i=1}^\dimv \a^i_j (0) \om_i|$.
Assume $d(b,V) < \frac{1-c}{16}$. Let $ \frac 12 < s_i < 1$.
Let $\l_1, \dots , \l_{m_1 + m_2}$ be complex numbers of modulus $s$.
Then
\be
\label{eqm30}
x' \= (\om_1 \l_1, \om_2 \l_1, \om_3 \l_2, \om_4 \l_2, \dots, \om_{2 m_{1} -1} \l_{m_1}, \om_{2m_1} \l_{m_1},
\l_{m_1 + 1}, \dots, \l_{m_1 + m_2})
\ee
is in $U_1 \times \dots \times U_{m_1} \times \D^{m_2}$, so the point $(x',\a(x')) \in V$.
By Lemma \ref{LEMp3}, the distance $d(b, (x', \a(x'))$ will only have active slots in the range $\{ 1, \dots, p-1 \}$.
Since we can rotate each $\l_i$ independently, by the argument in Step 2, we can choose $x'$ in the form \eqref{eqm30} so that
\[
\rho(\bar \o_{2i-1} \l_i, b_{2i-1}) = \rho (\bar\o_{2i} \l_i b_{2i}) = |s_i| + o(1)
\]
as $b \to 0$. So by increasing the smaller $|s_i|$ we 
can find a point of the form \eqref{eqm30} so that
 $d(b, (x', \a(x'))$ is at least $\dimv$ balanced and it is active with respect to at least one index in every $J_r$.
Relabelling the indices if necessary, we can assume that $1, \dots, \dimv$ are active, 
so
\begin{multline}
\notag
d(b, (x', \a(x')) = \rho(\bar \om_1 \l_1 , b_1) = \rho( \bar \om_2 \l_1 , b_2) =
\dots
= \rho(\bar  \om_{2m_1} \l_{m_1}, b_{2 m_1}) \\
\notag
= \dots = \rho(\l_{m_1+ m_2}, b_\dimv) .
\end{multline}
Now apply a M\"obius automorphism $\Phi$  that moves $(b_1, \dots, b_\dimv, \alpha(b_1, \dots, b_\dimv))$ to $0$.
After the automorphism, $V$ becomes $\{ (x, \tilde \a(x) ) : x \in \tilde U \}$ as in \eqref{eq97.1}.
The point $b$ goes to $\tilde b$, whose first $\dimv$  coordinates are $0$, and  $(x', \a(x'))$ goes to $(y, \tilde \a(y))$ where $y$ is of the form
\[
y = (\tau_1 s, \tau_2 s, \dots, \tau_\dimv s) ,
\]
for some $\tau \in \T^\dimv$, since $\rho(0, y_i) = s$ for $1 \le i \le \dimv$.
Moreover, since $T(\a)$ is open,  if $d(b,V)$ is small enough, we will have 
for each $1 \le r \le m_1$, the pair 
$(\tau_{2 r -1}, \tau_{2r})$ 
is in $T(\tilde \a_j)$ for those
$\a_j$ corresponding to the \roms\ $\Phi({\mathcal L}_r)$.

We have $ (y, \tilde \a(y))$ is at least $\dimv$-balanced with respect to $\tilde b$, and the first $\dimv$ indices are active.
 Now let $s \to 0$. The points $P_s = (\tau \cdot s, \tilde \a (\tau \cdot s))$ are all in $\tilde V$,
and eventually $d(P_s, \tilde b)$ will become active in a coordinate larger than $\dimv$. 
So there must be some value of $s$ at which $(P_s, \tilde b)$ is at least $(\dimv +1)$-balanced, and 
where the coordinates $1, \dots, \dimv$ are all active.
Now we can finish the argument as we did in Step 2. 
\ep

Lemma \ref{lemd10} follows immediately from Proposition \ref{PROm1}.

\section{Proof of Theorem \ref{thmb5}}
\label{secn}

An $n$ point extremal Pick problem on $\D^N$ is called degenerate if there is some $(n-1)$ point subproblem that is also extremal. It is called strictly $n$-dimensional if there is no interpolating function that depends on $n-1$ variables.

When $n=2$ the problem of characterizing the uniqueness set of a $2$ point Pick problem can be solved using Lempert's theorem. A similar phenomenon was proven to be true if $n=3$ \cite{kosthree}. In particular, the following theorem was proved.
\bt
\label{thmn1} 
Let $(\mu_1, \mu_2, \mu_3)$ and $(w_1, w_2, w_3)$ be the data for a strictly 3-dimensional 
nondegenerate
extremal Pick problem in $\D^N$.
Then there are an analytic disc $f: \D\to \D^N$ containing $\mu_1,$ $\mu_2$ and $\mu_3$ in its range and a Blaschke product $b$ of degree 2 such that 
every solution $\phi$ of the Pick problem satisfies 
\be
\label{eqn1}
\notag
\phi(f(\l)) = b(\l) .
\ee
\et

Let us make some simplifying  reductions  in this  theorem.
We can clearly assume that $\mu_1=0$, $w_1=0$ and that $b(\l) = \l m_{\g}(\l)$ for some $\g\in \D$.
Since the problem is   nondegenerate, none of the components  $f_j$ is a M\"obius map, and  since the problem is 3-dimensional,
at least three of them are Blaschke products of degree 2, say $f_{i_1},\ldots, f_{i_k}$. The generic case is when this is the case for exactly three components, i.e. $k=3$.  If $k>3$,
some additional book-keeping is required.

Since $f_j$ are not  Blaschke products of degree 2 for $j\not\in \{i_1,\ldots, i_k\}$, there is a solution to the problem that depends only on the $z_{i_1},\ldots, z_{i_k}$ variables.
Moreover for each $i_j$, there is some solution that depends only on three variables, one of which is $z_{i_j}$.
 Permuting, we can assume $i_j =j$, $j=1,\ldots, k$.
 Let $f_{j}(\l) = \om_{j}\l \vp_{a_{j}}(\l)$.

Let us write  $\conv(a_1, \ldots, a_k)$ to denote the interior of the convex hull of the set $\{ a_1, \dots , a_k \}$.
 It then follows
from \cite[Lemma 11]{kosthree} that
 $\gamma$ lies in ${\conv}(a_1, \ldots, a_k)$.

Calculations carried out in \cite{kosthree} and \cite{maci23} show the following
\begin{lemma}
\label{corn1} 
Let $a_{j_1}$, $a_{j_2}$, $a_{j_3}$ and $\g\in \overline{ \conv} (a_{j_1}, a_{j_2}, a_{j_3})$. 

(i) Suppose $\g$ lies on a side with vertices $a_{j_1}$, $a_{j_2}$, not equal to either endpoint, then there is $F_0$ in $b\mathcal M$ such that $F_0(\l \vp_{t a_{j_1}} (\l), \l \vp_{t a_{j_2}} (\l)) = \l \vp_{t\g}(\l)$, $t\in (0,1)$.

(ii) Suppose that $\g\in   \conv (a_{j_1}, a_{j_2}, a_{j_3})$. Then one can find solutions
$\phi_1(z) = F_1(F_2(z_{j_1}, z_{j_2}),z_{j_3})$ and $\phi_2(z)=F_3(z_{j_1}, F_4(z_{j_2}, z_{j_3}))$ where $F_r$ are in $b{\M}$ and of the form
$F_r = [\frac{\bar A_r}{|A_r|}; A_r, B_r]$ with $|A_r| + |B_r| = 1$.

The set $\{z:\ \phi_1(z_{j_1}, z_{j_2}, z_{j_3}) = \phi_2(z_{j_1}, z_{j_2}, z_{j_3})\}$ coincides with the variety
\[
\{ (z_{j_1}, z_{j_2}, z_{j_3}) \in \D^3 : z_{j_3} = \alpha(z_{j_1}, z_{j_2}) \}
\]
where $\a$ is the function in $\M$ given by
\[
\alpha = [- \o_1 \o_2 \bar \o_3 \frac{a_{j_2}-a_{j_1}}{\bar a_{j_2} - \bar a_{j_1}} ; \o_2 \frac{a_{j_2}- a_{j_3}}{\bar a_{j_2} - \bar a_{j_1}}, \o_1 \frac{a_{j_3} -a_{j_1}}{\bar a_{j_2} - \bar a_{j_1}}] .
\] 

\end{lemma}

With these tools at hand we can prove the theorem. Recall that a weak \roms\ is like a \roms\ but the functions are allowed to be in $\overline{\M}$.

{\bf Theorem \ref{thmb5}.}
{\em
Let $N \ge 3$.
A strictly 3-dimensional nondegenerate 3 point extremal Pick problem in $\D^\dimpol$ never has a  unique solution. 
 After
an automorphism of $\D^\dimpol$   the uniqueness variety is 
the Cartesian product of a   weak \roms\ of height $m$ and $\D^{\dimpol-m-2}$. Moreover, every such set  is the uniqueness variety for some 
strictly 3-dimensional nondegenerate 3 point extremal Pick problem.
}

\bp
Let $(\mu_1, \mu_2, \mu_3)$ and $(w_1, w_2, w_3)$ be the data for a strictly 3-dimensional nondegenerate extremal Pick problem in $\D^\dimpol$. Apply automorphisms and make reductions as in the comment after Theorem~\ref{thmn1}.

 By \cite[Sec. 3.1]{kosthree}, there are $a_1, \dots, a_k  \in \D$, 
 not all collinear, 
a point $\g$ in $\conv(a_1, \ldots, a_k)$
and
$\o_1, \dots , \o_k \in \T$, so that all solutions $\phi$ of the Pick problem satisfy
\be
\label{eqn3}
\phi(\o_1 \l \vp_{t a_1}(\l), \o_2 \l \vp_{t a_2} (\l), \dots, \o_k \l \vp_{t a_k} (\l),*) = \l \vp_{t \g }(\l)
\ee
for every $\l$ in $\D$ and $t \in (0,1]$, where $*$ means for any choice of the other variables.
From this we can conclude that
the uniqueness variety ${\mathcal U}$   contains the three real dimensional set swept out by varying $t$ and $\l$. Let $\mathcal U(a_1,\ldots, a_k)$ be its complexification to a complex variety in $\mathbb D^k$ of dimension 2. This complexification can be found explicitly.
If $\b_j = [\frac{a_{2}-a_{1}}{\bar a_{2} - \bar a_{1}} ; \frac{a_{2}- a_{j}}{\bar a_{2} - \bar a_{1}}, \frac{a_{j} -a_{1}}{\bar a_{2} - \bar a_{1}}]$, then $\l \vp_{t a_j}(\l)= \beta_j( \l \vp_{t a_1}(\l), \l \vp_{t a_2}(\l))$.
So if we let
\be
\label{eqn7}
\notag
\a_j = [-\bar \o_1 \bar \o_2 \o_j \frac{a_{2}-a_{1}}{\bar a_{2} - \bar a_{1}} ; \o_2 \frac{a_{2}- a_{j}}{\bar a_{2} - \bar a_{1}}, \o_1 \frac{a_{j} -a_{1}}{\bar a_{2} - \bar a_{1}}],
\ee
then  $\al_3,\ldots, \al_k\in \overline{\mathcal M}$ and \eqref{eqn3} says
 $\mathcal U(a_1,\ldots, a_k)\times \D^{N-k}$ is contained in $\mathcal U$,
 where
 \[
 \mathcal U(a_1,\ldots, a_k):=\{ x \in \D^k: x_j=\alpha_j(x_1, x_2), \ 3 \le j \le k \}.
 \]
Since
\[
 \o_2 \frac{a_{2}- a_{j}}{\bar a_{2} - \bar a_{1}} + \o_1 \frac{a_{j} -a_{1}}{\bar a_{2} - \bar a_{1}} \bar \o_1 \o_2
 =
 \o_2 \frac{a_{2}- a_{1}}{\bar a_{2} - \bar a_{1}},
 \]
 $\{ \alpha_3, \ldots, \alpha_k \}$ form a weak edge family.

So to prove the first part of the theorem, it is enough to show that $\mathcal U$ is actually equal to $\mathcal U(a_1, \ldots, a_k)$ and for this it suffices to show that $\mathcal U$ is contained in the graph of a holomorphic mapping of two  variables.

We shall use induction, with respect to the number of $a_j$'s 
forming the convex hull in which $\g$ is contained. 
Precisely, we shall show the following: if  $\g\in \conv(a_1, \ldots, a_\dimi)$,
for some $\dimi \le k$, where $a_1$, $a_2$, and $a_3$ are all distinct,
and
\be\label{eqcorn}
F(\l \vp_{a_1}(\l), \ldots, \l \vp_{a_\dimi}(\l)) = \l \vp_\g(\l)
\ee
 for some $F\in \mathcal O(\D^\dimi, \D)$, then there are  finitely many  $F_j\in \mathcal O(\D^\dimi, \D)$ that satisfy \eqref{eqcorn} and such that the intersection $\cap_j \{F_j =F\}$ is contained in
  the  graph of a holomorphic mapping of $x_1$ and $x_2$.
  For $\dimi=3$ the assertion follows from Lemma~\ref{corn1}, part (ii):
   the set where
  $\phi_1 = \phi_2$ is
an algebraic surface, and so agrees with 
 the set ${\mathcal U}(a_1, a_2, a_3)$ which is the graph of $\a_3$.
  Let $\dimi>3$. Up to a permutation of coordinates we can assume that $\g\in \overline \conv(a_1, \ldots,a_{\dimi-2}, a_{\dimi-1})$ and that $\g\in \conv(a_1, \ldots,a_{\dimi-2}, (1-s) a_{\dimi-1} +s a_\dimi)$ for $s>0$ small enough.

By Lemma \ref{corn1} (i), there is 
$F_s\in b\mathcal M$  such that $F_s(\l \vp_{a_{\dimi-1}}(\l), \l \vp_{a_\dimi} (\l)) =\l \vp_{(1-s) a_{\dimi-1} + s a_{\dimi}}(\l)$. 
By the inductive hypothesis applied to $\{ a_1,\ldots, a_{\dimi-2}, (1-s) a_{\dimi-1} + s a_\dimi \}$, we can find
functions $F_{j,s}$ 
satisfying
\[
F_{j,s} (\l \vp_{a_1}(\l), \ldots, \l \vp_{a_\dimi -2}(\l), \l \vp_{(1-s) a_{\dimi-1} + s a_\dimi}(\l) ) = \l \vp_\g(\l)
\]
and such that the intersection of the sets on which they all agree is contained in
a graph of $x_1$, $x_2$.
Do this for three distinct  small positive  $s$, say $s_1, s_2$ and $s_3$. 
Then the functions 
\[
F_{j_i;i}(z):=F_{j_i,s_i}(z_1, \ldots, z_{\dimi-2}, F_{s_i}(z_{\dimi-1}, z_\dimi)),
\]
for $i =1, 2,3$ 
and $1 \le j_i \le p_i$
form the set of  functions we are looking for.

Indeed,
\beq
\lefteqn{F_{j_i;i}(\l \vp_{a_1}(\l), \ldots, \l \vp_{a_\dimi}(\l))}\\
 &\ =&F_{j_i,s_i}(\l \vp_{a_1}(\l), \ldots, \l \vp_{a_{\dimi-2}}(\l), F_{s_i}(\l \vp_{a_{\dimi-1}}(\l), \l \vp_{a_{\dimi}}(\l)))\\
&\ =&F_{j_i,s_i}(\l \vp_{a_1}(\l), \ldots, \l \vp_{a_{\dimi-2}}(\l), \l \vp_{(1-s) a_{\dimi-1} + s a_{\dimi}}(\l))\\
&\ =&\l \vp_\g(\l);
\eeq
so \eqref{eqcorn} is satisfied.

 The set where all $F_{j_1, s_1}$ coincide,
 as $j_1$ ranges from $1$ to $p_1$,
 is a two-dimensional set in $\C^{\dimi -1}$, a graph over $x_1$ and $x_2$, 
 given by $x_\ell = \b_\ell(s_1, x_1, x_2)$ for $3 \le \ell \le n-1$. The set where all $F_{j_1;1}$ coincide is 3-dimensional,
since we want 
$F_{s_1}(x_{n-1}, x_n) = \b_{n-1}(s_1, x_1,x_2)$.

As $F_{s_1} \neq F_{s_2}$, the set where all $F_{j_2;2}$ also coincide will be one dimension lower,
and will have the  constraints
\begin{eqnarray}
\notag
F_{s_1}(x_{n-1}, x_n) &=& \b_{n-1}(s_1, x_1,x_2)\\
\label{eqn5}
 F_{s_2}(x_{n-1}, x_n) &=& \b_{n-1}(s_2, x_1,x_2).
\end{eqnarray}
Since $F_{s_i} \in b\M$, solving \eqref{eqn5} gives quadratic equations for $x_{n-1}$ and $x_n$,
with coefficients involving $\b_{n-1}(s_1, x_1,x_2)$, $\b_{n-1}(s_2, x_1, x_2)$, and the parameters of $F_{s_1}$ and $F_{s_2}$.
Requiring  also
\[
 F_{s_3}(x_{n-1}, x_n) = \b_{n-1}(s_3, x_1,x_2)
\]
reduces the quadratic equation to a linear equation, so $x_{n-1}$ and $x_n$ are both functions of $x_1$ and $x_2$.

Now we can drop the assumption that the $a_j$ be distinct. Assume $\{ a_1, \dots, a_m \}$ is a maximal set of distinct points, and by what we have just proved $\U (a_1, \dots, a_m)$ is the graph of a function $G$
of $x_1$ and $x_2$. For each $\ell$ between $m+1$ and $k$ there is some $p \in \{ 1, \dots, m \}$
so that $a_\ell = a_p$. We then extend $G$ by letting $G_\ell = G_p$.

To prove the second part of the theorem, it is sufficient to show that every weak  \roms\ is a uniqueness variety, since if ${\mathcal U}$ 
is the uniqueness variety for some 
extremal Pick problem in $\D^k$, we can just append $\dimpol - k$ zeroes to each $\mu_j$ to get 
a problem on $\D^\dimpol$ with uniqueness set ${\mathcal U} \times \D^{\dimpol -k}$.

So pick a weak \roms\ $\LL = \{ x \in \D^\dimpol : x_j = \a_j(x_1, x_2) , \ 3 \le j \le \dimpol \} $
where $\a_j$ form a weak edge family. Applying a rotation automorphism to $\D^\dimpol$, we can assume that
$\a_j = [ 1; A_j, B_j]$ and that $A_j + B_j = 1$.

By \cite[Sec. 3.1]{kosthree}, for any $a_1, \dots, a_ \dimpol$ in $\D$ and any $\g$ in the open convex hull
of $\{ a_i \}$, there is a function $\Phi: \D^\dimpol \to \D$ satisfying
\be
\label{eqn2}
\Phi( \l \vp_{a_i}(\l) ) = \l \vp_\g(\l) \qquad \forall \l \in \D .
\ee
Since the right-hand side in \eqref{eqn2} is a Blaschke product of degree 2, it follows that for any distinct triple $\l_1, \l_2, \l_3$, if we set $\mu_n =( \l_n \vp_{a_i}(\l_n)) $ and $w_n = \l_n \vp_\g(\l_n)$, then the 3 point problem mapping $\mu_n$ to $w_n$ is extremal.
By the previous part, the uniqueness set is equal to a weak \roms\ in $\D^\dimpol$.

We want to choose $\l_1 = 0$ and $w_1 = 0$ and 
choose $a_j$ so that 
\beq
 \frac{a_2 - a_1}{\bar a_2 - \bar a_1} &= & 1 \\
  \frac{a_2 - a_j}{\bar a_2 - \bar a_1}  &=& A_j.
\eeq
Then, by Lemma \ref{corn1}
the uniqueness variety for this extremal problem will be $\LL$.

Let $r > 0$ be small. Choose $a_2$  small,  let $a_1 = a_2 - r$, and $a_j = a_2 - rA_j$.
For $r$ sufficiently small, all the points $a_j$ will be in $\D$.
If some  $A_j$ is  not real, then $a_1, a_2, a_j$ are not collinear,
and we can find some $\g$ in the open convex hull of all the $a_i$'s.
Therefore
 by Lemma~\ref{corn1}
the uniqueness variety for this extremal problem will be $\LL$.

If every $A_j$ is real, then each $\alpha_j$ would be in $ \overline{\M} \setminus \M$,
and by definition of \wrs\ at least one must be in $\M$.
\ep

\section{Conclusion and Open Questions}
\label{seco}

We have shown that, among algebraic subvarieties of the polydisk, the Carathéodory sets are precisely the totally geodesic sets. Up to automorphism, every such set is a graph over a product of disks and Roman surfaces.

{\roms}s are important in understanding the complex geometry of the polydisk.
We also showed that weak {\roms}s, possibly times a polydisk, arise as uniqueness sets for strictly three-dimensional nondegenerate three-point extremal Pick problems.

The proof of our main result, the classification of \car sets $V \subseteq \D^\dimpol$, has several stages, because 
several local and global phenomena have to be controlled simultaneously.
 The simplest case is when the dimension of $V$ is $1$, because then there are no {\roms}s. The next simplest case is when the codimension is $1$, because then there is no Step 4 in the proof.
The next case is when either the dimension or codimension is $2$, since then Step 4 is simpler.
Once one reaches $\D^6$, the argument there is basically the same as for all higher $N$.

Here are some questions that the authors were unable to answer.

\begin{enumerate}
\item
We made the standing assumption that $V$ was the intersection of an algebraic set with the polydisk.
We did this so that when we parametrized $V$ locally, we knew the functions were Nash functions, and this 
made them more tractable. But they actually turn out to be rational functions in $\M$, suggesting the assumption is not necessary. Is Theorem \ref{thma1} true if we merely require
$V$ to be a holomorphic subvariety of $\D^\dimpol$?

\item
{\roms}s of height $1$ are all equivalent to $\K$ under an automorphism of $\D^3$.
We know that in addition to being \car sets, they occur as the uniqueness sets of 3 point 
 interpolation problems.
The function theory of $\K$ is fascinating. Here are three questions, each one stronger than its predecessor.

\begin{enumerate}
\item Is $\D^3$ a 3-extension set for $\K$? This means if there are points $\l_1, \l_2, \l_3$ in $\K$, $w_1, w_2, w_3$ in $\D$, and a holomorphic $\phi : \K \to \D$ such that $\phi(\l_j) = w_j$ for each $j$, does there exist a holomorphic $\Phi : \D^3 \to \D$ that also maps each $\l_j$ to $w_j$?

\item
Is $\D^3$ an extension set for $\K$, i.e.  an $n$-extension set for every $n$?

\item
Is there an And\^o inequality for $\K$?
An And\^o inequality would assert: whenever $T_1, T_2, T_3$ are commuting operators on a Hilbert space, each with norm less than $1$, and satisfying the equation $T_1 + T_2 + T_3 = T_1 T_2 + T_1 T_3 + T_2 T_3$, then
$\| p(T) \| \le \|p\|_\K$ for any polynomial $p$.
\end{enumerate}

There is a dichotomy for algebraic subsets of $\D^3$: either every set with the extension property is a retract, or every \car set
has the extension property. Resolving (b) determines which of these two statements is true.

\item
Whatever the answers to (a)-(c) are for $\K$, do the same answers apply to every \roms\ in $\D^\dimpol$?

\end{enumerate}

\section{Index of Notation}

\begin{tabular}{ll}
  Symbol & Usage \\
$[\tau; A,B]$ & $(x_1, x_2) \mapsto \tau  \frac{A x_1 + Bx_2 + x_1 x_2}{1 + \overline{B} x_1 + \overline{A} x_2}$\\
$\M$ & $ \{ [\tau; A,B] : |\tau| = 1, \ (|A|, |B|, 1) \ {\rm are\ triangular} \}$\\
$\kappa$ & $[1;-1,-1]$\\
$\K$ & $\{ (x_1,x_2, x_3) \in \D^3 : x_1 + x_2 + x_3 = x_1 x_2 + x_1 x_3 + x_2 x_3 \}$\\
$\a_j^i$ &	$\frac{\partial \a_j}{\partial x_i} (0) $\\
$T(\a_j(x_1, \dots, x_\dimv))$ & $\{ \o \in \T^\dimv: | \sum_{i=1}^\dimv \a_j^i \o_i | < 1 \}$\\
$T(\a)$& $ \cap_{j= \dimv+1}^\dimpol T(\a_j) $\\
$\hatt(\a_j(x_1, x_2)) $ & $\{ \o_2: (1, \o_2) \in T(\a_j) \}$\\
$D_{\a_j}$& $\{ x \in \D^2: |\a_j(x)| < 1 \}$\\
$D(\a_j(x_1, x_2)) $ & $\{ \z \in \D : |\a_j^1 + \a_j^2 \z| < 1 \}$\\
$\vp_\mu : \D \to \D$& $z \mapsto \frac{z-\mu}{1 - \bar \mu z}$\\
$\Phi_a : \D^\dimpol \to \D^\dimpol$& $\Phi_a(x_1, \dots, x_\dimpol) = (\vp_{a_1}(x_1), \dots, \vp_{a_\dimpol} (x_\dimpol)) $\\
$\tilde \a_j(x_1, \dots, x_\dimv) $& $\vp_{ \a_j(a_1, \dots, a_\dimv)}\circ \a_j( \vp_{-a_1}(x_1), \dots, \vp_{-a_\dimv} x_\dimv)$\\

\end{tabular}

\section{AI disclosure}

AI was not used at all in initially proving the theorems, nor in the first drafts of this paper, except for drawing the pictures. We ran the paper through AI at the end, 
and it found a few gaps and
inconsistencies, which we have fixed. The only change that was more than a line or two was in the proof  in Section \ref{seck}, where we realized that Lemmas \ref{lemk0} and \ref{lemk02} gave a more direct proof of Lemma \ref{lemd6}  than our original argument.  We formulated Lemma \ref{lemk0} after we realized this, but then  we used  AI to complete the proof of Lemma \ref{lemk0}.

 The proof was too complicated for the AI model we used to say whether or not it thought the complete proof was correct. The authors believe the theorem is correct and that we have a complete proof, though we acknowledge that in any human-written paper of this length and complexity there may be small errors. 
We take full responsibility for everything in this paper.

\vskip 25pt
{\bf Acknowledgement:} The first author is grateful to Jakub Byszewski for useful discussions.


\begin{thebibliography}{10}

\bibitem{ag91}
Jim Agler.
\newblock Operator theory and the {Carath\'eodory metric}.
\newblock {\em Inventiones Math.}, 101:483--500, 1991.

\bibitem{akmc22}
Jim Agler, {\L}ukasz Kosi\'{n}ski, and John~E. M\raise.45ex\hbox{c}Carthy.
\newblock Norm preserving extensions of holomorphic functions defined on
  varieties in {$\mathbb{C}^n$}.
\newblock {\em J. Funct. Anal.}, 283(9):Paper No. 109636, 39, 2022.

\bibitem{akmc23}
Jim Agler, {\L}ukasz Kosi\'{n}ski, and John~E. M\raise.45ex\hbox{c}Carthy.
\newblock Complete norm-preserving extensions of holomorphic functions.
\newblock {\em Israel J. Math.}, 255(1):251--263, 2023.

\bibitem{aly19}
Jim Agler, Zinaida Lykova, and Nicholas Young.
\newblock Geodesics, retracts, and the norm-preserving extension property in
  the symmetrized bidisc.
\newblock {\em Mem. Amer. Math. Soc.}, 258(1242):vii+108, 2019.

\bibitem{agmc_dv}
Jim Agler and John~E. McCarthy.
\newblock Distinguished varieties.
\newblock {\em Acta Math.}, 194(2):133--153, 2005.

\bibitem{agmcvn}
Jim Agler and John~E. M\raise.45ex\hbox{c}Carthy.
\newblock Norm preserving extensions of holomorphic functions from subvarieties
  of the bidisk.
\newblock {\em Ann. of Math. (2)}, 157(1):289--312, 2003.

\bibitem{and63}
Tsuyoshi {And\^o}.
\newblock On a pair of commutative contractions.
\newblock {\em Acta Sci. Math. (Szeged)}, 24:88--90, 1963.

\bibitem{bds20}
Sibaprasad Barik, B.~Krishna Das, and Jaydeb Sarkar.
\newblock Isometric dilations and von {N}eumann inequality for finite rank
  commuting contractions.
\newblock {\em Bull. Sci. Math.}, 165:102915, 25, 2020.

\bibitem{bs18}
Tirthankar Bhattacharyya and Haripada Sau.
\newblock Holomorphic functions on the symmetrized bidisk---realization,
  interpolation and extension.
\newblock {\em J. Funct. Anal.}, 274(2):504--524, 2018.

\bibitem{brga18}
Filippo Bracci and Herv\'e Gaussier.
\newblock A proof of the {M}uir-{S}uffridge conjecture for convex maps of the
  unit ball in {$\Bbb C^n$}.
\newblock {\em Math. Ann.}, 372(1-2):845--858, 2018.

\bibitem{bgnt24}
Filippo Bracci, Herv\'e{} Gaussier, Nikolai Nikolov, and Pascal~J. Thomas.
\newblock Local and global visibility and {G}romov hyperbolicity of domains
  with respect to the {K}obayashi distance.
\newblock {\em Trans. Amer. Math. Soc.}, 377(1):471--493, 2024.

\bibitem{br20}
A.~Brudnyi.
\newblock Oka principle on the maximal ideal space of {$H^\infty$}.
\newblock {\em Algebra i Analiz}, 31(5):24--89, 2019.

\bibitem{cara27}
C.~Carath\'eodory.
\newblock {\"U}ber eine spezielle {M}etrik, die in der {T}heorie der
  analytischen {F}unktionen auftritt.
\newblock {\em Atti Pontifica Acad. Sc.}, Nuovi Lincei 80:135--141, 1927.

\bibitem{car51}
Henri Cartan.
\newblock {\em S\'eminaire {Henri Cartan} 1951/2}.
\newblock W.A. Benjamin, New York, 1967.

\bibitem{dks25}
B.~Krishna Das, Poornendu Kumar, and Haripada Sau.
\newblock Distinguished varieties and the {N}evanlinna-{P}ick problem on the
  symmetrized bidisk.
\newblock {\em Math. Z.}, 310(3):Paper No. 40, 29, 2025.

\bibitem{den08}
Zofia Denkowska and Jacek Stasica.
\newblock {\em Ensembles sous-analytiques à la polonaise}.
\newblock Number~69 in Travaux en Cours. Hermann, 2008.

\bibitem{dey17}
Michael~A. Dritschel, Daniel Est\'evez, and Dmitry Yakubovich.
\newblock Traces of analytic uniform algebras on subvarieties and test
  collections.
\newblock {\em J. Lond. Math. Soc. (2)}, 95(2):414--440, 2017.

\bibitem{gkvw09}
Anatolii Grinshpan, Dmitry~S. Kaliuzhnyi-Verbovetskyi, Victor Vinnikov, and
  Hugo~J. Woerdeman.
\newblock Classes of tuples of commuting contractions satisfying the
  multivariable von {N}eumann inequality.
\newblock {\em J. Funct. Anal.}, 256(9):3035--3054, 2009.

\bibitem{ghw08}
Kunyu Guo, Hansong Huang, and Kai Wang.
\newblock Retracts in polydisk and analytic varieties with the
  {$H^\infty$}-extension property.
\newblock {\em J. Geom. Anal.}, 18(1):148--171, 2008.

\bibitem{har25}
Michael Hartz.
\newblock On von {N}eumann's inequality on the polydisc.
\newblock {\em Math. Ann.}, 391(4):5235--5264, 2025.

\bibitem{hs81}
L.~F. Heath and T.~J. Suffridge.
\newblock Holomorphic retracts in complex {$n$}-space.
\newblock {\em Illinois J. Math.}, 25(1):125--135, 1981.

\bibitem{hu94}
Xiao~Jun Huang.
\newblock A preservation principle of extremal mappings near a strongly
  pseudoconvex point and its applications.
\newblock {\em Illinois J. Math.}, 38(2):283--302, 1994.

\bibitem{jp99}
Marek Jarnicki and Peter Pflug.
\newblock A remark on the product property for the {C}arath\'eodory
  pseudodistance.
\newblock volume~53, pages 89--94. 1999.
\newblock XII-th Conference on Analytic Functions (Lublin, 1998).

\bibitem{jp13}
Marek Jarnicki and Peter Pflug.
\newblock {\em Invariant distances and metrics in complex analysis}, volume~9
  of {\em De Gruyter Expositions in Mathematics}.
\newblock Walter de Gruyter GmbH \& Co. KG, Berlin, extended edition, 2013.

\bibitem{kn07b}
Greg Knese.
\newblock A {S}chwarz lemma on the polydisk.
\newblock {\em Proc. Amer. Math. Soc.}, 135:2759--2768, 2007.

\bibitem{kn11b}
Greg Knese.
\newblock A refined {A}gler decomposition and geometric applications.
\newblock {\em Indiana Univ. Math. J.}, 60(6):1831--1841, 2011.

\bibitem{ko67}
Shoshichi Kobayashi.
\newblock Invariant distances on complex manifolds and holomorphic mappings.
\newblock {\em J. Math. Soc. Japan}, 19:460--480, 1967.

\bibitem{kosthree}
{\L}ukasz Kosi\'nski.
\newblock Three-point {N}evanlinna-{P}ick problem in the polydisc.
\newblock {\em Proc. Lond. Math. Soc. (3)}, 111(4):887--910, 2015.

\bibitem{kmc19}
{\L}ukasz Kosi\'{n}ski and John~E. M\raise.45ex\hbox{c}Carthy.
\newblock Norm preserving extensions of bounded holomorphic functions.
\newblock {\em Trans. Amer. Math. Soc.}, 371(10):7243--7257, 2019.

\bibitem{kmc20}
{\L}ukasz Kosi\'{n}ski and John~E. M\raise.45ex\hbox{c}Carthy.
\newblock Extensions of bounded holomorphic functions on the tridisk.
\newblock {\em Rev. Mat. Iberoam.}, 36(3):791--816, 2020.

\bibitem{kz21}
{\L}ukasz Kosi\'{n}ski and Wlodzimierz Zwonek.
\newblock Extension property and universal sets.
\newblock {\em Canad. J. Math.}, 73(3):717--736, 2021.

\bibitem{lem81}
L\'aszl\'o Lempert.
\newblock La m\'etrique de {Kobayashi} et la repr\'esentation des domaines sur
  la boule.
\newblock {\em Bull. Soc. Math. France}, 109:427--484, 1981.

\bibitem{lem84}
L\'aszl\'o Lempert.
\newblock Intrinsic distances and holomorphic retracts.
\newblock In {\em Complex analysis and applications '81 ({V}arna, 1981)}, pages
  341--364. Publ. House Bulgar. Acad. Sci., Sofia, 1984.

\bibitem{loj91}
S.~{\L }ojasiewicz.
\newblock {\em Complex Analytic Geometry}.
\newblock Birkh{\"a}user, Boston, 1991.

\bibitem{maci19}
Krzysztof Maciaszek.
\newblock On polynomial extension property in {N}-disc.
\newblock {\em Complex Anal. Oper. Theory}, 13(6):2771--2780, 2019.

\bibitem{maci20}
Krzysztof Maciaszek.
\newblock Polynomial extension property in the classical {C}artan domain
  {${\mathcal R}_{II}$}.
\newblock {\em Arch. Math. (Basel)}, 115(6):657--666, 2020.

\bibitem{maci23}
Krzysztof Maciaszek.
\newblock Geometry of uniqueness varieties for a three-point {P}ick problem in
  {$\mathbb{D}^3$}.
\newblock {\em Bull. Sci. Math.}, 188:Paper No. 103337, 12, 2023.

\bibitem{mark18}
Vladimir Markovic.
\newblock Carath\'eodory's metrics on {T}eichm\"uller spaces and {$L$}-shaped
  pillowcases.
\newblock {\em Duke Math. J.}, 167(3):497--535, 2018.

\bibitem{mmw17}
Curtis~T. McMullen, Ronen~E. Mukamel, and Alex Wright.
\newblock Cubic curves and totally geodesic subvarieties of moduli space.
\newblock {\em Ann. of Math. (2)}, 185(3):957--990, 2017.

\bibitem{mum76}
David Mumford.
\newblock {\em Algebraic geometry. {I}}, volume No. 221 of {\em Grundlehren der
  Mathematischen Wissenschaften}.
\newblock Springer-Verlag, Berlin-New York, 1976.
\newblock Complex projective varieties.

\bibitem{nt26}
Nikolai Nikolov and Pascal~J. Thomas.
\newblock Quasi-triangle inequality for the {L}empert function.
\newblock {\em Complex Anal. Synerg.}, 12(1):Paper No. 3, 9, 2026.

\bibitem{ps14}
Sourav Pal and Orr~Moshe Shalit.
\newblock Spectral sets and distinguished varieties in the symmetrized bidisc.
\newblock {\em J. Funct. Anal.}, 266(9):5779--5800, 2014.

\bibitem{rud69}
Walter Rudin.
\newblock {\em Function theory in polydiscs}.
\newblock W. A. Benjamin, Inc., New York-Amsterdam, 1969.

\bibitem{si75}
Nessim Sibony.
\newblock Prolongement des fonctions holomorphes born\'ees et m\'etrique de
  {C}arath\'eodory.
\newblock {\em Invent. Math.}, 29(3):205--230, 1975.

\bibitem{tw90}
Piotr Tworzewski.
\newblock Intersections of analytic sets with linear subspaces.
\newblock {\em Ann. Scuola Norm. Sup. Pisa Cl. Sci. (4)}, 17(2):227--271, 1990.

\bibitem{vonN51}
John von Neumann.
\newblock {Eine Spektraltheorie f\"ur allgemeine Operatoren eines unit\"aren
  Raumes}.
\newblock {\em Math. Nachr.}, 4:258--281, 1951.

\end{thebibliography}
\end{document}